\documentclass[3p,times,11pt]{elsarticle}

\graphicspath{{./figures/TB_spline_ex/},{./figures/TB_spline_ex/tcheb_spline/},{./figures/examples/},{./figures/ex_quarter_annulus_adv_tan/},
	{./figures/ex_mickey_mouse/},{./figures/ex_adv_sqr/}, {./figures/ex_adv_radial/}, {./figures/tikz/}, {./figures/}}

\journal{Computer Methods in Applied Mechanics and Engineering}

\usepackage{amssymb, amsthm, amsmath, amsfonts, bm}
\usepackage{amsmath}

\usepackage{hyperref}
\usepackage{cleveref}

\usepackage{lineno}

\usepackage{array}
\usepackage{multirow}

\usepackage{tcolorbox}
\usepackage{pbox}

\usepackage[]{subcaption} 
\usepackage[font = small]{caption} 

\usepackage[toc,page]{appendix} 
\usepackage{mathtools} 

\usepackage{tikz}
\usepackage{xspace}
\usepackage{xcolor}

\theoremstyle{plain}
\newtheorem{theorem}{Theorem}
\newtheorem{proposition}[theorem]{Proposition}

\newtheorem{remark}[theorem]{Remark}
\newtheorem{example}[theorem]{Example}

\theoremstyle{definition}
\newtheorem{definition}[theorem]{Definition}

\newcommand{\peclet}{P\'eclet\xspace}

\newcommand{\angleSpace}[1]{\Big\langle {#1} \Big\rangle}

\newcommand{\myInterval}{J}
\newcommand{\bbR}{\mathbb{R}}
\newcommand{\bbN}{\mathbb{N}}
\newcommand{\bbZ}{\mathbb{Z}}
\newcommand{\bbT}{\mathbb{T}}
\newcommand{\bbS}{\mathbb{S}}
\newcommand{\bbP}{\mathbb{P}}

\newcommand{\bbV}{\mathbb{V}}
\newcommand{\bbW}{\mathbb{W}}
\newcommand{\xx}{\mathbf{x}}

\newcommand{\myi}{\mathrm{i}}
\newcommand{\rootVec}{{\scriptstyle{\mathcal{W}}}}

\newcommand{\smthVec}{\mathbf{r}}
\newcommand{\map}{\mathbf{G}}
\newcommand{\spaceT}{\mathbb{P}}

\usepackage{settobox} 
\newbox{\subbox}
\newlength{\xheight}
\newlength{\subheight}
	\newcommand\xheightsub[2]{%
    \savebox{\subbox}{\({}_{#2}\)}%
    \settoheight{\xheight}{\({}_x\)}%
    \settoboxheight{\subheight}{\subbox}%
    \addtolength{\subheight}{-\xheight}%
    {#1}{\raisebox{-\subheight}{\usebox{\subbox}}}%
 }	

\modulolinenumbers[5]
\usepackage{etoolbox}
\makeatletter
\patchcmd{\ps@pprintTitle}
  {Preprint submitted}
  {Published}
  {}{}
\makeatother
\begin{document}					

\begin{frontmatter}			
\title{Tchebycheffian B-splines in isogeometric Galerkin methods}
	
\corref{cor1}
\author[address1]{Krunal Raval}
\ead{raval@mat.uniroma2.it}	
	
\author[address1]{Carla Manni}
\ead{manni@mat.uniroma2.it}

\author[address1]{Hendrik Speleers}
\ead{speleers@mat.uniroma2.it}
	
\cortext[cor1]{Corresponding author}

\address[address1]{Department of Mathematics, University of Rome Tor Vergata, 00133 Rome, Italy}

\begin{abstract}
 Tchebycheffian splines are smooth piecewise functions whose pieces are drawn from (possibly different) Tchebycheff spaces, a natural generalization of algebraic polynomial spaces.
 They enjoy most of the properties known in the polynomial spline case.
 In particular, under suitable assumptions, Tchebycheffian splines admit a representation in terms of basis functions, called Tchebycheffian B-splines (TB-splines), completely analogous to polynomial B-splines.
 A particularly interesting subclass consists of Tchebycheffian splines with pieces belonging to null-spaces of constant-coefficient linear differential operators.
 They grant the freedom of combining polynomials with exponential and trigonometric functions with any number of individual shape parameters. Moreover, they have been recently equipped with efficient evaluation and manipulation procedures. 
 In this paper we consider the use of TB-splines with pieces belonging to null-spaces of constant-coefficient linear differential operators as an attractive substitute for standard polynomial B-splines and rational NURBS in isogeometric Galerkin methods. We discuss how to exploit the large flexibility of the geometrical and analytical features of the underlying Tchebycheff spaces according to problem-driven selection strategies. TB-splines offer a wide and robust environment for the isogeometric paradigm beyond the limits of the rational NURBS model.
\end{abstract}

\begin{keyword}
Isogeometric analysis, Tchebycheffian B-splines, B-splines with shape parameters
\end{keyword}
\end{frontmatter}

\section{Introduction}
Spline functions are ubiquitous in numerical methods. Besides their theoretical interest,
they have application in several branches of the sciences including geometric
modeling, signal processing, data analysis, visualization and numerical simulation just to mention a few.
Splines, in the broad sense of the term, are functions consisting of pieces of smooth
functions glued together in a certain smooth way. There is a large variety of spline species, often referred to as the zoo of splines. The most popular species is the one where the pieces
are algebraic polynomials of a given degree $p$ and inter-smoothness is imposed by means of equality of derivatives up to a certain order. Their popularity can be mainly
attributed to their representation in terms of the so-called B-splines. B-splines enjoy properties such as local linear independence, minimal support, non-negativity, and partition of unity; they can be computed through a stable recurrence relation and can even be seen as the geometrically
optimal basis for piecewise polynomial spaces.

Tchebycheffian splines are smooth piecewise functions whose pieces are drawn from (possibly different) extended Tchebycheff (ET-) spaces which are natural generalizations of algebraic polynomial spaces \cite{KarlinS1966,Mazure2011,Schumaker2007}. 
Any non-trivial element of an ET-space of dimension $p+1$ has at most $p$ zeros counting  multiplicity. We will refer to $p$ as the degree, in analogy with the polynomial case.
Extended complete Tchebycheff \mbox{(ECT-)} spaces are an important subclass that can be generated through a set of positive weight functions  \cite{Mazure2011b,Schumaker2007} and are spanned by generalized power functions \cite{lyche2019tchebycheffian}; the latter are the natural extension of the monomial basis functions for algebraic polynomials; see \Cref{sec-ECT-spaces}. Relevant examples are null-spaces of linear differential operators on suitable intervals \cite{Schumaker2007}; see \Cref{sec_null_spc}. On bounded and closed intervals the concepts of ET-space and ECT-space coincide \cite{Mazure2007}, so from now on we will focus on ECT-spaces.

Most of the results known for splines in the polynomial case extend in a natural way to the Tchebycheffian setting.
In particular, under suitable assumptions on the involved ECT-spaces, Tchebycheffian splines admit a representation in terms of basis functions, called Tchebycheffian B-splines (TB-splines), with similar properties to polynomial B-splines.
TB-splines were introduced in 1968 by Karlin \cite{Karlin1968} using generalized divided differences. 
There are several other ways to define them, including Hermite interpolation \cite{BuchwaldM2003,NurnbergerSSS1984}, de Boor-like recurrence relations \cite{DynR1988,Lyche1985}, integral recurrence relations \cite{BisterP1997}, and blossoming \cite{Mazure2011}; see also the historical notes in \cite[Chapters~9 and~11]{Schumaker2007} for further details. Each of these definitions has advantages according to the properties to be proved and lead to the same functions, up to a proper scaling.
 
Multivariate extensions of Tchebycheffian splines can be easily obtained by means of (local) tensor-product structures \cite{BraccoLMRS2016-gb,BraccoLMRS2016}.

Due to the richness of ECT-spaces, Tchebycheffian splines can find applications in several contexts including data approximation/interpolation, geometric modeling and numerical simulation; see \cite{lyche2019tchebycheffian} and references therein.
Thanks to their structural similarities, TB-splines are theoretically plug-to-plug compatible with classical polynomial B-splines, so they can be potentially easily incorporated in any software library supporting polynomial B-splines to enrich its capability.

Unfortunately, despite their theoretical interest and applicative potential, TB-splines have not gained much attention in practice so far. The reason behind this is that TB-splines are generally difficult to compute. The classical approaches mentioned above, based on generalized divided differences, Hermite interpolation, or repeated integration, are computationally expensive and/or numerically unstable.
An important step forward was recently made in \cite{hiemstra2020tchebycheffian} where the authors proposed a strategy that represents TB-splines as linear combinations of local Tchebycheffian Bernstein functions through a suitable extraction operator. The local Tchebycheffian Bernstein functions form a basis of the local ECT-spaces involved in the definition of the TB-splines. In the polynomial case, these are nothing but the classical Bernstein polynomial basis functions.
Following the approach in \cite{hiemstra2020tchebycheffian}, an object-oriented Matlab toolbox has been developed in \cite{speleers2022algorithm} for the construction and manipulation of TB-splines whenever they exist. The toolbox supports TB-splines whose pieces belong to ECT-spaces that are null-spaces of constant-coefficient linear differential operators, and is publicly available. Note that both \cite{hiemstra2020tchebycheffian} and \cite{speleers2022algorithm} address the more general setting of multi-degree TB-splines (MDTB-splines) where the local ECT-spaces of the Tchebycheffian spline space are not required to be of the same dimension. 
Having at our disposal a publicly available implementation for such a large class of TB-splines paves the path for their effective use in practical applications. An alternative strategy has been proposed in \cite{beccari2022} but no implementation is available.
 
Isogeometric analysis (IgA) is a well established methodology for the analysis of problems governed by partial differential equations (PDEs); see \cite{IGAbook,hughes2005}. 
IgA aims to simplify the interoperability between geometric modeling and numerical simulation by constructing a fully integrated framework for computer-aided design (CAD) and finite element analysis (FEA). In its original formulation, IgA uses B-splines, or their rational generalization (NURBS), both to represent the computational domain and to approximate the solution of the PDE that models the problem of interest. The isogeometric paradigm shows important advantages over classical FEA: complicated geometries can be represented more accurately and some common profiles, such as conic sections, are described exactly; the precise geometry description incorporated at the coarsest mesh level can be maintained during mesh refinement through the process of knot insertion; and the inherently high smoothness of B-splines/NURBS allows for a higher accuracy per degree of freedom.
Nevertheless, the above advantages are not a distinguishing property of B-splines/NURBS, and B-splines/NURBS are not a requisite ingredient in IgA.
 
An interesting subclass of TB-splines, the so-called generalized polynomial B-splines (GB-splines), has been proposed as an alternative to classical polynomial and rational splines in the context of IgA; see \cite{manni2011generalized,ManniRS2017} and references therein.
 GB-splines can be seen as the minimal extension of polynomial B-splines towards the wide variety of TB-splines: their pieces  belong to ECT-spaces obtained by enriching an algebraic polynomial space with a pair of functions, typically hyperbolic or trigonometric functions identified by a single shape parameter.   
 GB-splines overcome two issues of the rational NURBS model: they allow for an exact (almost) arc-length parameterization of conic sections and behave with respect to differentiation and integration as nicely as polynomial B-splines (for instance, the derivative of a trigonometric generalized spline with a given phase parameter $\beta$ and degree $p$ is again a trigonometric generalized spline with the same shape parameter $\beta$ and degree $p-1$).
 These properties make GB-splines an interesting tool to face both geometrical and analytical hurdles in IgA. The effectiveness of hyperbolic or trigonometric GB-splines in isogeometric Galerkin and collocation methods has been illustrated in a sequence of papers, where their properties have been exploited to obtain exact representations of common geometries \cite{Aimi,manni2011generalized,ManniPS2014} or to beneficially deal with advection-dominated problems \cite{manni2011tension,manni2015isogeometric}; see also \cite{ManniRS2017} for the related spectral properties. 
 
Being a so minimal extension of the polynomial setting, however, GB-splines are not always flexible enough for practical applications. In particular, in a tensor-product GB-spline space only two hyperbolic or trigonometric functions identified by the same shape parameter are added to polynomials along each parametric direction. Therefore, a given tensor-product GB-spline space does not allow for:
 \begin{itemize}
\item an exact representation of different arcs of conic sections at opposite sides;   	
\item a proper treatment of different analytic features (like sharp layers) along a given parametric direction;
\item a simultaneous treatment of geometrical and analytical features along the same parametric direction.
 \end{itemize} 

Tchebycheffian splines with pieces belonging to ECT-spaces that are null-spaces of constant-coefficient linear differential operators are an extension of hyperbolic and trigonometric generalized polynomial splines. They  enjoy all the geometrical and analytical features that motivate the interest in GB-splines without suffering from the above mentioned restrictions. They grant the freedom of combining polynomials with exponential and trigonometric functions with any number of individual shape parameters. For this class of Tchebycheffian splines, when the various pieces are drawn from a single ECT-space which contains constants, the existence of TB-splines is always ensured, possibly with some restriction on the partition; see \Cref{sec_TB_splines} for more details. Furthermore, they are supported by the Matlab toolbox available in \cite{speleers2022algorithm}. 
In summary, TB-splines with pieces belonging to ECT-spaces that are null-spaces of constant-coefficient linear differential operators offer a suitable balance between the immense variety of ECT-spaces and the practical needs of a problem-driven space selection and efficient evaluation procedures for the space elements. They can be seen as a flexible environment for the IgA paradigm, beyond the limits of the rational NURBS model.
In this paper we present and discuss the use of such TB-splines in Galerkin isogeometric methods.

The remainder of the paper consists of four sections. 
In the next section we introduce ECT-spaces and Tchebycheffian spline spaces. We also briefly discuss the construction and the properties of the corresponding TB-spline basis.
 \Cref{sec-TB-IGA} is devoted to the use of TB-splines in isogeometric analysis and to describe a general road map for the selection of suitable ECT-spaces in a problem-driven fashion. A bunch of case studies illustrating the performance of the proposed approach are presented in \Cref{sec-numerics}.
 \Cref{sec-conclusion} collects some concluding remarks and outlines possible future research lines.

\section{TB-spline theory}
In this section we outline the mathematical framework for the spaces we are going to use. We recall from the literature the definition of ECT-spaces and we summarize some of their main properties.
We also discuss a particular subclass of ECT-spaces identified by linear differential operators with real constant coefficients that will be our target class for isogeometric analysis.
Afterwards, we introduce Tchebycheffian splines -- piecewise functions where the pieces belong to given ECT-spaces glued together with certain smoothness -- and we briefly discuss the construction of TB-splines, which are B-spline like basis functions for Tchebycheffian splines. We refer the reader to \cite{lyche18,lyche2019tchebycheffian,Mazure2011,Schumaker2007} for in-depth information about Tchebycheffian splines and spline functions in general.

The presentation of this section is somehow technical for the sake of completeness. However, all the technical details can be skipped at a first reading, focusing on the features of the (spline) spaces we are going to describe.

\subsection{Extended complete Tchebycheff spaces}
\label{sec-ECT-spaces}

We consider a class of $(p+1)$-dimensional spaces that are a natural extension of the space of algebraic polynomials of degree less than or equal to $p$.

\begin{definition}[Extended Tchebycheff space]\label{def_et_space}
Given an interval $\myInterval$ and an integer $p\geq 0$, a space $\bbT_p(\myInterval) \subset C^{p}(\myInterval)$ of dimension $p+1$ is an extended Tchebycheff (ET-) space on $\myInterval$, if any Hermite interpolation problem with $p+1$ data on $\myInterval$ has a unique solution in $\bbT_p(\myInterval)$. More in details, let $z_1,z_2,\ldots,z_{\bar m}$ be distinct points on $\myInterval$, with any positive integer $\bar m$ and $d_1,d_2,\ldots,d_{\bar m}$ be non-negative integers such that $p+1 = \sum_{i=1}^{\bar m}(d_i+1)$. Then for any $f_{i,j} \in \bbR$, there exists a unique $g\in\bbT_p(\myInterval)$ such that
$$ D^jg(z_i) = f_{i,j}, \quad j = 1,\ldots,d_i, \quad i = 1,\ldots,\bar m.$$
\end{definition}

As an immediate consequence of the previous definition, we have that any non-trivial element of an ET-space on $\myInterval$ has at most $p$ roots in $\myInterval$, counting multiplicity. This property shows the strong interconnection between ET-spaces and the space of algebraic polynomials.

We are now going to identify an interesting subclass of ET-spaces that strengthens the similarity with algebraic polynomials even further. We do this by generalizing the monomial basis functions, which form the so-called power basis of the polynomial space. We denote by $\langle g_0,\ldots,g_k \rangle$ the linear space spanned by $\{ g_0,\ldots,g_k \}$.

\begin{definition}[Extended complete Tchebycheff space]\label{def_ect_space}
Given an interval $\myInterval$ and an integer $p\geq 0$, a space $\bbT_p(\myInterval) \subset C^{p}(\myInterval)$ of dimension $p+1$ is an extended complete Tchebycheff (ECT-) space on $\myInterval$, if there exist functions $g_0,\ldots, g_p$ such that $\bbT_p(\myInterval) =\langle g_0,\ldots,g_p\rangle $ and every subspace $\langle g_0,\ldots,g_k \rangle$ for $k=0,\ldots,p$ is an ET-space on $\myInterval$. The basis $\big\{g_0,\ldots,g_p\big\}$ is called an ECT-system.
\end{definition}

\begin{example}\label{ex_algbraic_pol}
The space of algebraic polynomials $\bbP_p=\langle 1,x,\ldots,x^p \rangle$ is an ECT-space on the real line. Indeed, with any fixed point $z\in\bbR $, it can be accounted as the span of the ECT-system 
\begin{equation}
	\label{eq-powers}
	\Bigg\{ 1,x-z,\frac{(x-z)^2}{2},\ldots,\frac{(x-z)^p}{p!}\Bigg\}.
\end{equation}
\end{example}
It is evident from \Cref{def_ect_space} that every ECT-space on the interval $\myInterval$ is an ET-space of the same dimension on the same interval. The contrary case is not always true. However, any ET-space on a bounded and closed interval is an ECT-space on the same interval \cite{Mazure2007}.
In the scope of this paper, we are only interested in bounded closed intervals. Hence, the notions of ET-space and ECT-space are interchangeable. Throughout the paper, we consider ECT-spaces to allow for a neat construction of a proper basis of related spline spaces. In this perspective, we recall some additional properties of ECT-spaces.

A space $\bbT_p(\myInterval) \subset C^{p}(\myInterval)$ of dimension $p+1$ is an ECT-space on $\myInterval$ if and only if there exist positive functions $w_{j} \in C^{p-j}(\myInterval)$, $j=0,1,\ldots, p$ such that  $\bbT_p(\myInterval)$ is spanned by the following functions \cite{lyche2019tchebycheffian,Mazure2007,Schumaker2007}:
\begin{equation}\label{eq_gen_powers}
\begin{cases}
g_0(x) \coloneqq w_0(x), \\
g_1(x) \coloneqq w_0(x) \int_{z}^x w_1(y_1) {\rm d}y_1, \\
\hspace*{1.1cm} \vdots \\
g_p(x) \coloneqq w_0(x) \int_{z}^x w_1(y_1) \int_{z}^{y_1} \cdots \int_{z}^{y_{p-1}} w_p(y_p) {\rm d}y_{p} {\rm d}y_{p-1}\cdots {\rm d}y_1,
\end{cases}
\end{equation}
where $z$ is any point in $\myInterval$.
The functions $g_0, \ldots, g_p$ are called generalized power functions, the functions $w_0, \ldots, w_p$ are usually referred to as weights, and the set $\{w_0,\ldots,w_p\}$ is called a weight system generating $\bbT_p(\myInterval)$.
 \begin{example}\label{ex_algbraic_pol_w} 
 	A weight system for the space of algebraic polynomials $\bbP_p$ is given by
	$ \{w_j=1 : j=0,\ldots,p\}$. Through \eqref{eq_gen_powers} these weights provide the basis \eqref{eq-powers}.
\end{example}
\begin{example}\label{ex_gb_exp_w}
	Let $1\leq\ell\leq p$ and $\alpha_1,\ldots,\alpha_\ell\in \mathbb{R}$ with $\alpha_i\neq\alpha_j$ for each $i\neq j$. Then, we consider the functions  
	\begin{equation}\label{eq-gb-exp-w}
	w_0(x)= \dots= w_{p-\ell}(x)=1,\quad  w_{p-\ell+1}(x)=e^{\alpha_1 x}, \quad w_{p-\ell+k}(x) =e^{(\alpha_k-\alpha_{k-1})x},\quad k=2,\ldots, \ell. 
	\end{equation}
	For any $z\in \mathbb{R}$, the identified generalized power functions in \eqref{eq_gen_powers} span the space
	$$
	\angleSpace{ 1,x,\ldots,x^{p-\ell},e^{\alpha_1 x},\ldots ,e^{\alpha_\ell x} }.
	$$
	Since the functions in \eqref{eq-gb-exp-w} are positive on $\mathbb{R}$, the above space is a $(p+1)$-dimensional ECT-space on $\mathbb{R}$ for any $p\geq 1$.	
\end{example}
\begin{example}\label{ex_gb_trig_w}
	For $\beta\in \mathbb{R}$, $\beta>0$, and $x\in\myInterval$ with
	$$\myInterval=\left(-\frac{\pi}{2\beta},\frac{\pi}{2\beta} \right),$$
	we consider the functions  
\begin{equation}\label{eq-gb-trig-w}
	w_0(x)= \dots= w_{p-2}(x)=1,\quad  w_{p-1}(x) =\cos(\beta x),\quad w_p(x)=\frac{1}{\cos^2(\beta x)}, \quad p\geq 1.
\end{equation}
	 For any $z\in \myInterval$, the identified generalized power functions in \eqref{eq_gen_powers} span the space
	\begin{equation}\label{eq-gb-trig-space}
	\angleSpace{ 1,x,\ldots,x^{p-2},\cos(\beta x), \sin(\beta x) }.
	\end{equation}
Since the functions in \eqref{eq-gb-trig-w} are positive on $\myInterval$, the above space is a $(p+1)$-dimensional ECT-space on $\myInterval$ for any $p\geq 1$. The space in \eqref{eq-gb-trig-space} is often referred to as cycloidal space.
Actually, for $p\geq 2$, it can be shown that the space in \eqref{eq-gb-trig-space} is an ECT-space on any interval of length less than $2\pi/\beta$; see \cite{carnicer2017critical} and also \cite[Section~2]{lyche2019tchebycheffian}.
\end{example}

\begin{remark}\label{rmk_weights_factor}
	A given ECT-space on an interval $\myInterval$ can be identified by different systems of weights; see \cite[Section~2]{lyche2019tchebycheffian} for details and examples. In particular, it is clear that
	the two weight systems
	\begin{equation*}
	w_0,\ldots,w_p \quad {\rm and} \quad K_0w_0,\ldots,K_pw_p,
	\end{equation*}
	where $K_0,\ldots, K_p$ are positive constants, identify the same ECT-space.
\end{remark}

\begin{remark}\label{rmk_gen-powers-der}
Let $\bbT_p$ be a space of dimension $p+1$ on $\myInterval$. 
	Assume that the derivative space of $\bbT_p$, i.e., the space spanned by the derivatives of functions in $\bbT_p$, is a $p$-dimensional ECT-space on the interval $\myInterval$ generated by the weights $\{w_1,\ldots, w_p\}$. Then, $\bbT_p$ is a  $(p+1)$-dimensional ECT-space on $\myInterval$ generated by the weights $\{1,w_1,\ldots,w_p\}$. 
\end{remark}

Many properties of polynomials involve working with derivatives. When dealing with ECT-spaces generated by the weight system $\{w_0,\ldots,w_p\}$ it is often convenient to replace the usual derivatives by some related differential operators. We define $D_0f \coloneqq f$ and
$$D_jf \coloneqq D\left(\frac{f}{w_{j-1}}\right), \quad j=1,\ldots,p+1,$$
and we set
\begin{equation}\label{eq_L}
L_j \coloneqq D_jD_{j-1}\cdots D_0,\quad j=0,\ldots,p+1.
\end{equation}
The operator $L_j$ can be seen as a natural substitute for $D^j$; see \cite[Section~9.1]{Schumaker2007}. In particular, we have 
$$
L_j g_k(x) = \begin{cases}
w_{j}(x) \int_{z}^x w_{j+1}(y_{j+1}) \int_{z}^{y_{j+1}} \cdots \int_{z}^{y_{k-1}} w_k(y_k) {\rm d}y_{k} {\rm d}y_{k-1}\cdots {\rm d}y_{j+1}, & k=j,\ldots,p,\\
0, & k=0,\ldots,j-1.
\end{cases}
$$
where $g_k$ are defined in \eqref{eq_gen_powers}. This confirms that the generalized power functions are the natural extension of the monomial basis functions for algebraic polynomials.

\subsection {ECT-spaces and null-spaces of linear differential operators} 
\label{sec_null_spc}

A large and interesting class of ECT-spaces is given by the null-spaces of linear differential operators with real constant coefficients. More precisely,
let us consider the differential operator defined by
\begin{equation}\label{eq_diff_op}
\mathcal{L}_pf \coloneqq 
D^{p+1}f +\sum_{j=0}^{p}a_jD^jf,\quad f\in C^{p+1}(\bbR),\quad a_j\in\bbR,\quad j=0,\ldots, p,
\end{equation}
and let
\begin{equation}\label{eq_char_pol}
\mathfrak{p}_p(\omega) \coloneqq \omega^{p+1}+\sum_{j=0}^p a_j\omega^j
\end{equation}
be its characteristic polynomial and $\bbN_p$ its null-space.

It is well known that the null-space of the differential operator \eqref{eq_diff_op}
is easily described by means of the roots of the characteristic polynomial \eqref{eq_char_pol}.
Let $\omega=\alpha+\myi\beta$ be a root of multiplicity $\mu\geq1$ of the polynomial in \eqref{eq_char_pol} for $\alpha,\beta\in\bbR$ and $\myi\coloneqq\sqrt{-1}$. Then, this root generates the following
fundamental subspace:
\begin{itemize}
	\item if $\beta=0$, then  $$\angleSpace{x^k e^{\alpha x} : k=0,\ldots,\mu-1} \subseteq \bbN_p;$$ 
	\item if $\beta\neq0$, then the complex conjugate of $\omega$ is also a root of multiplicity $\mu$, and $$\angleSpace{x^k e^{\alpha x}\cos(\beta x), x^k e^{\alpha x}\sin(\beta x): k=0,\ldots,\mu-1} \subseteq \bbN_p.$$ 	
\end{itemize}
The above subspaces for all roots together form the null-space $\bbN_p$.
To ensure that the constants belong to $\bbN_p$, we have to assume that $\omega=0$ is a root of the characteristic polynomial in \eqref{eq_char_pol}.

\begin{example}	Setting $a_j=0$, $ j=0,\ldots,p$ in \eqref{eq_diff_op} gives the linear differential operator $\mathcal{L}_p f = D^{p+1} f$. In this case, the null-space is the space of algebraic polynomials $\bbN_p=\bbP_p$; see also \Cref{ex_algbraic_pol}. 
\end{example}

\begin{remark}
	\label{rem_null_space}
	If the polynomial \eqref{eq_char_pol} has only real roots, then $\bbN_p$ is an ECT-space on the whole real line, and in particular on any bounded interval $[a,b]$. On the other hand, if the characteristic polynomial has also complex roots, then $\bbN_p$ is an ECT-space on sufficiently small intervals; see \cite{Schumaker2007}. More precisely, $\bbN_p$ is an ECT-space on any interval $[a,b]$ such that 
	$$ b-a<\mathfrak{l}_p. $$
	The value $\mathfrak{l}_p$ is referred to as critical length and it can be bounded from below as
	\begin{equation}
	\label{eq_critical_length}
	\mathfrak{l}_p\geq \pi/B>0,
	\end{equation}  
	where $B$ denotes the maximum imaginary part of all non-real roots of the characteristic polynomial. For a detailed study about the critical lengths of such ECT-spaces, the reader is referred to \cite{beccari2020critical} and \cite{carnicer2003critical}.
\end{remark}
In the rest of the paper we will mainly focus on $(p+1)$-dimensional ECT-spaces that are null-spaces of linear differential operators with real constant coefficients, as in \eqref{eq_diff_op}. The reference interval $\myInterval$ ensuring the ECT-property will be specified when necessary.
We will always assume that constants belong to the considered ECT-spaces. 
Let us now introduce the notation we are going to use for such spaces and give some examples.

Consider linear differential operators whose characteristic polynomial has roots
 $\omega_k = \alpha_k+\myi\beta_k$
 of multiplicity $\mu_k \geq 1$, with $k=0,\ldots,M$, such that $\sum_{k=0}^M \mu_k = p+1$, and $\omega_0=0$. Under this assumption, the corresponding null-space is uniquely characterized by the following vector with $p+1-\mu_0$ components:
\begin{equation}\label{eq_root_vec}
\rootVec \coloneqq \Big(  \underbrace{\omega_1,\,\ldots,\,\omega_1}_{\mu_1 \text{ times}},\,\underbrace{\omega_2,\,\ldots,\,\omega_2}_{\mu_2 \text{ times}},\, \ldots\,,\underbrace{\omega_M,\,\ldots,\,\omega_M}_{\mu_M \text{ times}} \Big).
\end{equation}
To highlight the dependence of the null-space on the roots \eqref{eq_root_vec}, we use the notation
\begin{equation}
\label{eq_null_space}
\spaceT_p^\rootVec =\spaceT_p^{( \alpha_1+\myi\beta_1,\ldots,\alpha_1+\myi\beta_1, \alpha_2+\myi\beta_2,\ldots,\alpha_2+\myi\beta_2, \ldots, \alpha_M+\myi\beta_M) }.
\end{equation}

The ECT-spaces constructed from the null-spaces of linear differential operators with real constant coefficients  offer great flexibility through the shape parameters to be chosen in \eqref{eq_root_vec}. A particular selection of these shape parameters results into familiar subclasses of ECT-spaces. We now list some popular examples. 
\begin{itemize}
\item The class of algebraic polynomial spaces is the most established subclass (see \Cref{ex_algbraic_pol}):
\begin{equation*} 
\spaceT_p = \angleSpace{1,x,\ldots,x^p}.
\end{equation*}
\item Generalized polynomial spaces, see \cite{lyche2019tchebycheffian,manni2011generalized}, enrich algebraic polynomial spaces by a pair of functions.  An important class of generalized polynomial spaces can be obtained as null-spaces of linear differential operators with constant coefficients. A first example is the algebraic polynomial space of degree $p-2$ enriched with two exponential functions. Taking $\omega_1=\alpha$, $\omega_2=-\alpha$, and $\mu_1=\mu_2=1$ gives the following null-space:
\begin{equation}\label{eq_gb_space_exp}
\spaceT_p^{(\alpha,-\alpha)} = \angleSpace{1,x,\ldots,x^{p-2},e^{\alpha x},e^{-\alpha x} } = \angleSpace{1,x,\ldots,x^{p-2},\cosh(\alpha x),\sinh(\alpha x) },\quad p\geq2;
\end{equation}
see also \Cref{ex_gb_exp_w}. This is an ECT-space on any interval according to \Cref{rem_null_space}.
Another instance of generalized polynomial spaces is the combination of the algebraic polynomial space of degree $p-2$ with two trigonometric functions sharing the same phase $\beta>0$. Taking $\omega_1=\myi\beta$, $ \omega_2=-\myi\beta$, and $\mu_1=\mu_2=1$ gives the following null-space:
\begin{equation}\label{eq_gb_space_trig}
\spaceT_p^{(\myi\beta,-\myi\beta)} = \angleSpace{ 1,x,\ldots,x^{p-2},\cos(\beta x), \sin(\beta x) },\quad p\geq2,
\end{equation}
which is a cycloidal space; see also \Cref{ex_gb_trig_w}. \Cref{rem_null_space} provides a lower bound on its critical length, while it can be bounded from above by 
$$\mathfrak{l}_p \leq\frac{2\pi}{\beta} \left\lfloor \frac{p}{2}\right \rfloor, \quad p\geq 2;$$
see \cite{carnicer2003critical}. A detailed and precise study on the critical length of cycloidal spaces can be found in \cite{carnicer2017critical}.

\item An interesting subclass of ECT-spaces obtained as null-spaces of linear differential operators are the following exponential and trigonometric spaces, both defined for $p=2q$:
\begin{align*}
\spaceT_p^{(\alpha,-\alpha,\ldots,q\alpha,-q\alpha)} &= \angleSpace{1,e^{\alpha x},e^{-\alpha x},e^{2\alpha x},e^{-2\alpha x},\ldots,e^{q\alpha x},e^{-q\alpha x}}\\
&= \angleSpace{1,\cosh(\alpha x),\sinh(\alpha x),\cosh(2\alpha x),\sinh(2\alpha x),\ldots,\cosh(q\alpha x),\sinh(q\alpha x)},
\end{align*}
and
\begin{equation*}
\spaceT_p^{(\myi\beta,-\myi\beta,\ldots,\myi q\beta,-\myi q\beta)}= \angleSpace{1,\cos(\beta x), \sin(\beta x), \ldots, \cos(q\beta x), \sin(q\beta x) },
\end{equation*}
identified by the vector \eqref{eq_root_vec} with $\omega_k=\pm k\alpha$ and $\omega_k=\pm \myi k\beta$, $k=1,\ldots,q$, respectively. In the former case, the space is an ECT-space on any interval, while in the latter case, the space is an ECT-space on any interval of length less than $\frac{\pi}{\beta}$ for $\beta>0$; see \cite{Mazure2005,Sanchez-R-1998}.
\end{itemize}
All the above examples are special instances of the following subclass of ECT-spaces that combines polynomial, exponential, and trigonometric functions. Let $0\leq \ell \leq p$ and $\alpha_1,\ldots,\alpha_\ell\in \mathbb{R}$ with $\alpha_i\neq\alpha_j$ for each $i\neq j$. Moreover, let $0\leq 2q\leq p-\ell$ and $\beta_1,\ldots,\beta_q\in \mathbb{R}$ with $\beta_i\neq\beta_j$ for each $i\neq j$.
Then, let us consider the space
\begin{equation}\label{eq_spaces} 
\spaceT_p^{(\alpha_1,\ldots,\alpha_\ell,\myi\beta_1,-\myi\beta_1,\ldots,\myi\beta_{q},-\myi \beta_{q})} 
 = \angleSpace{1,x,\ldots,x^{p-\ell-2q},e^{\alpha_1x},\ldots,e^{\alpha_\ell x},\cos(\beta_1 x), \sin(\beta_1 x),\ldots,\cos(\beta_{q} x), \sin(\beta_{q} x) },
\end{equation}
defined on a reference interval ensuring the ECT-property. Such space is invariant under translation of the interval. 
We are going to utilize this comprehensive subclass of ECT-spaces in the case studies presented in \Cref{sec-numerics}.
In the following we will refer to $p$ as the degree of the ECT-space. 

\begin{remark}\label{rem_ect_benefits}
The subclass of ECT-spaces in \eqref{eq_spaces} is a proper subclass of null-spaces of differential operators with real constant coefficients, see  \eqref{eq_null_space}, but still grants the freedom of combining polynomials with exponential and trigonometric functions with any number of individual shape parameters to be selected according to (automatic) problem-driven strategies. 
Moreover, derivatives and integrals of functions belonging to spaces of the form \eqref{eq_spaces} belong to spaces of the same form, possibly not including constants anymore.
\end{remark}

\subsection{Tchebycheffian spline spaces}
\label{sec-T-splines}
In a complete similarity with polynomial splines, Tchebycheffian splines are piecewise functions with pieces  belonging to ECT-spaces which can be glued together with certain prescribed smoothness.
More precisely, let $\Delta$ be a partition of the interval $[a,b]\subset\bbR$ specified in terms of a sequence of breakpoints,
$$\Delta \coloneqq \{a=: x_0<x_1<\dots<x_{m-1}<x_m\coloneqq b \}. $$
We set $\myInterval_i\coloneqq[x_{i-1},x_i)$, $i=1,\ldots,m-1$, and $\myInterval_m\coloneqq[x_{m-1},x_m]$. Furthermore, for $i=1,\ldots, m$ let $\bbT_{p,i}$ be an ECT-space of dimension $p+1$ on the closed interval $[x_{i-1},x_i]$. 
Finally, let us consider a sequence of $m-1$ integers
$$\smthVec \coloneqq \Big\{ r_i\in\bbZ: -1\leq r_i\leq p-1,\, i=1,\ldots,m-1\Big\}.$$ 
 The elements of the set   
\begin{equation}\begin{split}\label{eq_tcheb_sp} 
\bbS_p^{\smthVec}(\Delta) \coloneqq \Big\{f:[a,b]\rightarrow \bbR:\, &f|_{\myInterval_i}\in \bbT_{p,i},\ i=1,\ldots,m  ; \\& D^l_-f(x_i)=D^l_+f(x_i),\ l=0,\ldots,r_i,\ i=1,\ldots,m-1 \Big\}
\end{split}\end{equation}
are Tchebycheffian splines of degree $p$ and smoothness $\smthVec$ on the partition $\Delta$ of the interval $[a,b]$.
The dimension of the space $\bbS_p^{\smthVec}(\Delta)$ is given by
\begin{equation} 
\label{eq_dim_spline_space}
n \coloneqq p+1+\sum_{i=1}^{m-1}(p-r_i).
\end{equation}
The reader is referred to \cite{lyche2019tchebycheffian,Schumaker2007} for more details.

When building the spline space $\bbS_p^{\smthVec}(\Delta)$, all the pieces of the spline functions can be taken from the ``same'' ECT-space on the whole interval $[a,b]$, or from different ECT-spaces considered on the different intervals. The latter case provides a more general framework which allows us to locally exploit the full richness of ECT-spaces; however, it also entails theoretical difficulties, mainly concerning the existence of a basis possessing all the nice properties of classical polynomial B-splines. This will be discussed in the next subsection.

Being interested in the application of Tchebycheffian splines in IgA, it is important to understand the approximation power of the space in \eqref{eq_tcheb_sp}. It turns out that, in general, smooth functions can be approximated by Tchebycheffian splines with the same orders of approximation as in the case of classical polynomial splines. The reader is referred to \cite[Section~9.7]{Schumaker2007} for details and we state a simplified version of \cite[Theorem~9.38]{Schumaker2007} below. Let $\bbW^{r,q}([a,b])$ be the standard Sobolev space on $[a,b]$ and let $\|\cdot\|_q$ denote the standard $L^q$-norm on $[a,b]$ for $1\leq q\leq\infty$.

\begin{theorem}\label{thm_approx}
Let $\bbT_p$ be a $(p+1)$-dimensional ECT-space on the interval $[a, b]$ and 
let $\bbS_p^{\smthVec}(\Delta)$ be a Tchebycheffian spline space where $\bbT_{p,i}=\bbT_p$ for $i=1,\ldots,m$.
For any $f\in \bbW^{p+1,q}([a,b])$ there exists a spline $Qf\in\bbS_p^{\smthVec}(\Delta)$ such that
$$
\|f-Qf\|_q\leq C\, (h_\Delta)^{p+1}\,\|L_{p+1}f\|_q,
$$
where $L_{p+1}$ is defined in \eqref{eq_L},
$$
h_\Delta \coloneqq \max_{i=1,\ldots,m}(x_{i}-x_{i-1}),
$$
and $C$ is a constant independent of $f$ and $\Delta$.
\end{theorem}

\subsection{Tchebycheffian B-splines}
\label{sec_TB_splines}
We now introduce a basis of the space $\bbS_p^{\smthVec}(\Delta)$ that possesses all the characterizing properties of classical polynomial B-splines. In particular, the elements of this basis form a non-negative partition of unity and have minimal support. These basis functions will be referred to as Tchebycheffian B-splines (TB-splines).

The existence of TB-splines requires some assumptions on the sequence of ECT-spaces
$\bbT_{p,i}$, $i=1,\ldots,m$ 
in \eqref{eq_tcheb_sp} 
which can be expressed in terms of weight systems; see \Cref{sec-ECT-spaces}.
\begin{definition}[Admissible weights]\label{def_admissible_w} 
For $i=1,\ldots, m$, let $\{w_{j,i} : j=0,\ldots,p\}$ be a weight system generating the ECT-space $\bbT_{p,i}$ on the interval $[x_{i-1}, x_i]$. These weight systems are admissible for the space $\bbS_p^{\smthVec}(\Delta)$ in \eqref{eq_tcheb_sp}, if 
$$w_{0,i}=1, \quad i=1,\ldots,m,$$
and for $i=1,\ldots,m-1$ and $j=0,\ldots,r_i$, 
$$ D^l_{-}w_{j,i}(x_i) = D^l_{+}w_{j,i+1}(x_i),\quad l=0,\ldots,r_i-j.$$
\end{definition}
Given a set of admissible weight systems for the space $\bbS_p^{\smthVec}(\Delta) $, we define the global weight system $\{w_j : j=0,\ldots,p\}$ as
\begin{equation}
\label{eq-admissible-w}
w_j(x) \coloneqq  w_{j,i}(x), \quad  x\in \myInterval_i, \quad i=1,\ldots,m.
\end{equation}

Similar to the standard polynomial B-splines, Tchebycheffian B-splines can be defined using a vector of non-decreasing (open) knots, 
\begin{equation}\label{eq_knot_vec} 
\bm{\xi}\coloneqq (\xi_k)_{k=1}^{n+p+1} \coloneqq \Bigl(\underbrace{x_0,\,\ldots,\, x_0}_{p+1 \text{ times}},  \underbrace{x_1,\,\ldots,\, x_1}_{p-r_1 \text{ times}},\, \ldots\,,\underbrace{x_{m-1},\,\ldots,\,x_{m-1}}_{p-r_{m-1} \text{ times}},\underbrace{x_m,\,\ldots,\,x_m}_{p+1 \text{ times}} \Bigr),
\end{equation}
with $n$ being the dimension of the space $\bbS_p^{\smthVec}(\Delta)$; see \eqref{eq_dim_spline_space}.
\begin{definition}[TB-splines]	\label{def-TB-splines}
Assume there exist admissible weight systems as in \Cref{def_admissible_w} for the space $\bbS_p^{\smthVec}(\Delta) $.
Given the knot vector \eqref{eq_knot_vec} and the global weight system \eqref{eq-admissible-w}, the
TB-splines $\{N_{k,p} : k=1,\ldots,n\}$ are defined recursively as follows \cite{lyche2019tchebycheffian}: for $x\in[a, b]$, $q = 0,\ldots,p$ and $k=1,\ldots,n+p-q$,
\begin{equation*}
N_{k,0}(x)\coloneqq \begin{cases} w_p(x),\quad & x\in[\xi_{k},\xi_{k+1}),\\
0,\quad & \text{otherwise},
\end{cases}
\end{equation*}
and 
\begin{equation*}
N_{k,q}(x)\coloneqq w_{p-q}(x) \int_{a}^{x} \Bigg[\frac{N_{k,q-1}(y)}{d_{k,q-1}} - \frac{N_{k+1,q-1}(y)}{d_{k+1,q-1}}\Bigg]\, {\rm d}y,\quad q>0,
\end{equation*}
where
\begin{equation*}
d_{j,q-1}\coloneqq \int_a^b N_{j,q-1}(y)\, {\rm d}y,
\end{equation*}
and if $d_{j,q}=0$ then
\begin{equation*}
\int_{a}^{x} \frac{N_{j,q}(y)}{d_{j,q}}\, dy \coloneqq \begin{cases} 1,\quad  &x\geq \xi_{j+q+1},\\
0,\quad & \text{otherwise.}\end{cases}
\end{equation*}
\end{definition}

TB-splines exhibit a complete structural similarity with polynomial B-splines. In particular, they enjoy the nice properties collected in the following proposition; see \cite{lyche2019tchebycheffian,Mazure2011,NurnbergerSSS1984} and references therein.
\begin{proposition}\label{prop-TB-splines}
Let $\{N_{k,p} : k=1,\ldots,n\}$
be the functions in \Cref{def-TB-splines}. Then, 
$$\bbS_p^{\smthVec}(\Delta)=\langle N_{k,p}:k=1,\ldots,n\rangle,$$
and the following properties hold (see \Cref{{fig_tb_spline_properties}}).
 \begin{itemize}
	\item Non-negativity: $ N_{k,p}(x)>0$ for all $x\in(\xi_k,\xi_{k+p+1})$;
	\item Local support: $ N_{k,p}(x)=0$ for all $x\notin[\xi_k,\xi_{k+p+1}] $;
	\item Partition of unity: $ \sum_{k=1}^n N_{k,p}(x)=1$ for all $x\in[a,b]$;
	\item Local linear independence: the set $\{ N_{k,p} \}_{k=j-p}^j$ forms a basis of $\bbT_{p,i_j+1}$ on $\myInterval_{i_j+1}$ for all $p+1\leq j\leq n$, where $i_j$ denotes the index such that $\xi_j=x_{i_j}$;
	\item Interpolation at end-points: \begin{align*}
	N_{1,p}(a)&=1,\quad N_{k,p}(a)=0,\quad k=2,\ldots,n;\\
	N_{n,p}(b)&=1,\quad N_{k,p}(b)=0,\quad k=1,\ldots,n-1. 
	\end{align*}
\end{itemize}
\end{proposition}

\begin{figure}[t!]
\centering
\includegraphics[width=0.5\linewidth]{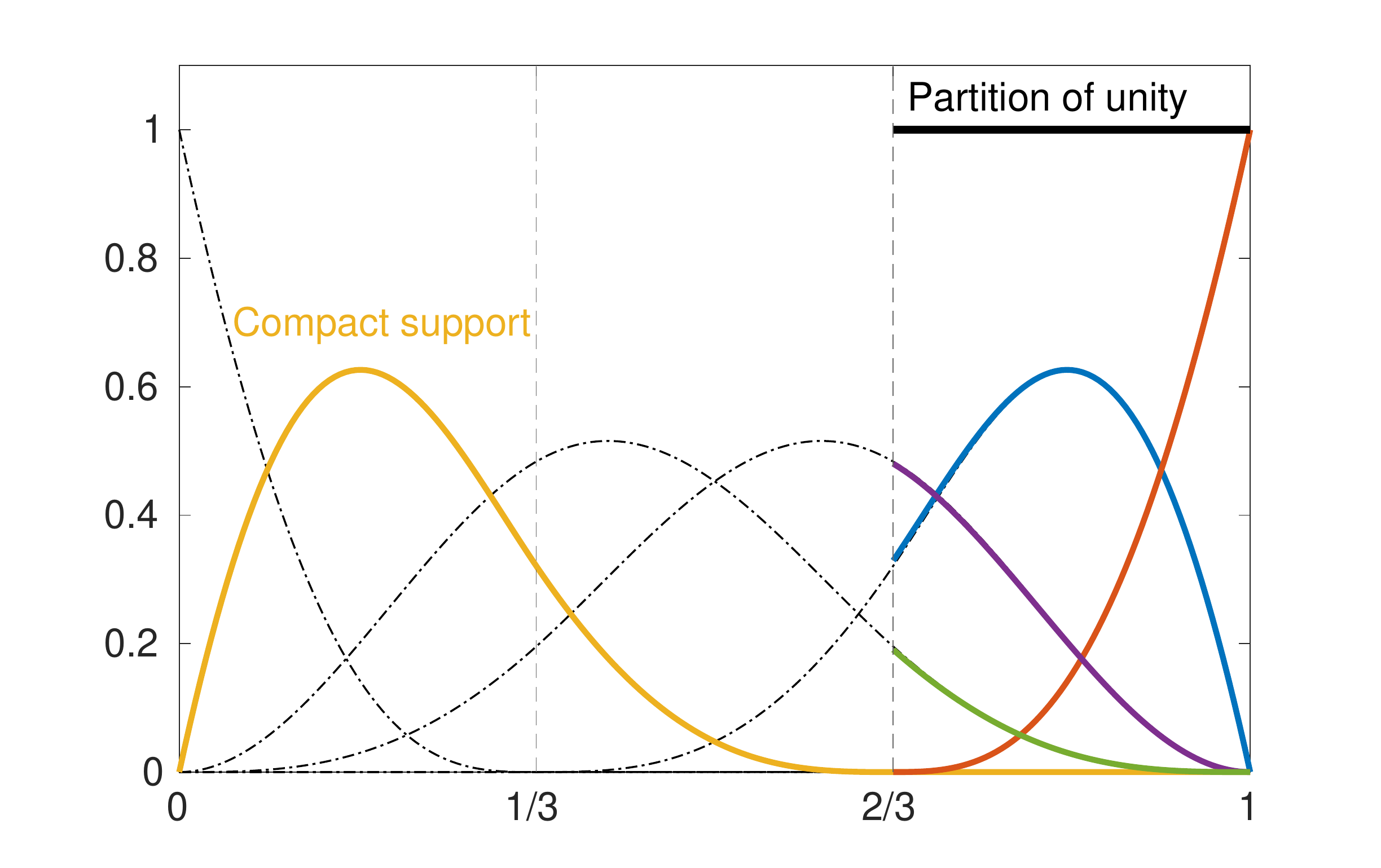}
\caption{Illustration of TB-spline properties.}
\label{fig_tb_spline_properties} 
\end{figure}

\begin{figure}[t!]
	\captionsetup[subfigure]{aboveskip=-1pt,belowskip=0cm}
	\centering
	\begin{subfigure}[t]{0.4\linewidth}
		\centering
		\includegraphics[width=\linewidth]{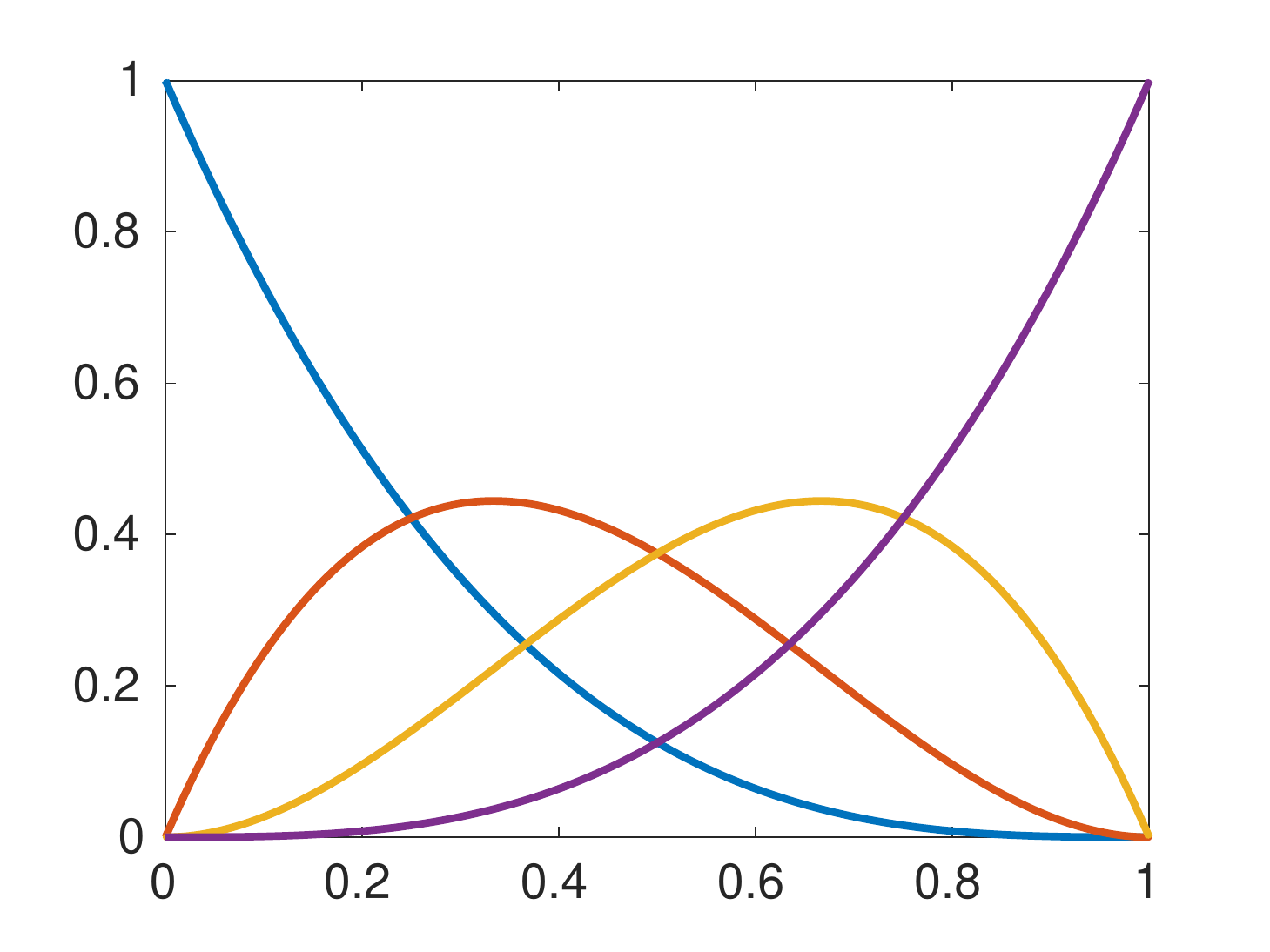} 
		\caption{$\spaceT_3 = \langle\, 1,x,x^2,x^3\, \rangle$} 
	\end{subfigure}
	\hspace{.75cm}
	\begin{subfigure}[t]{0.4\linewidth}
		\centering
		\includegraphics[width=\linewidth]{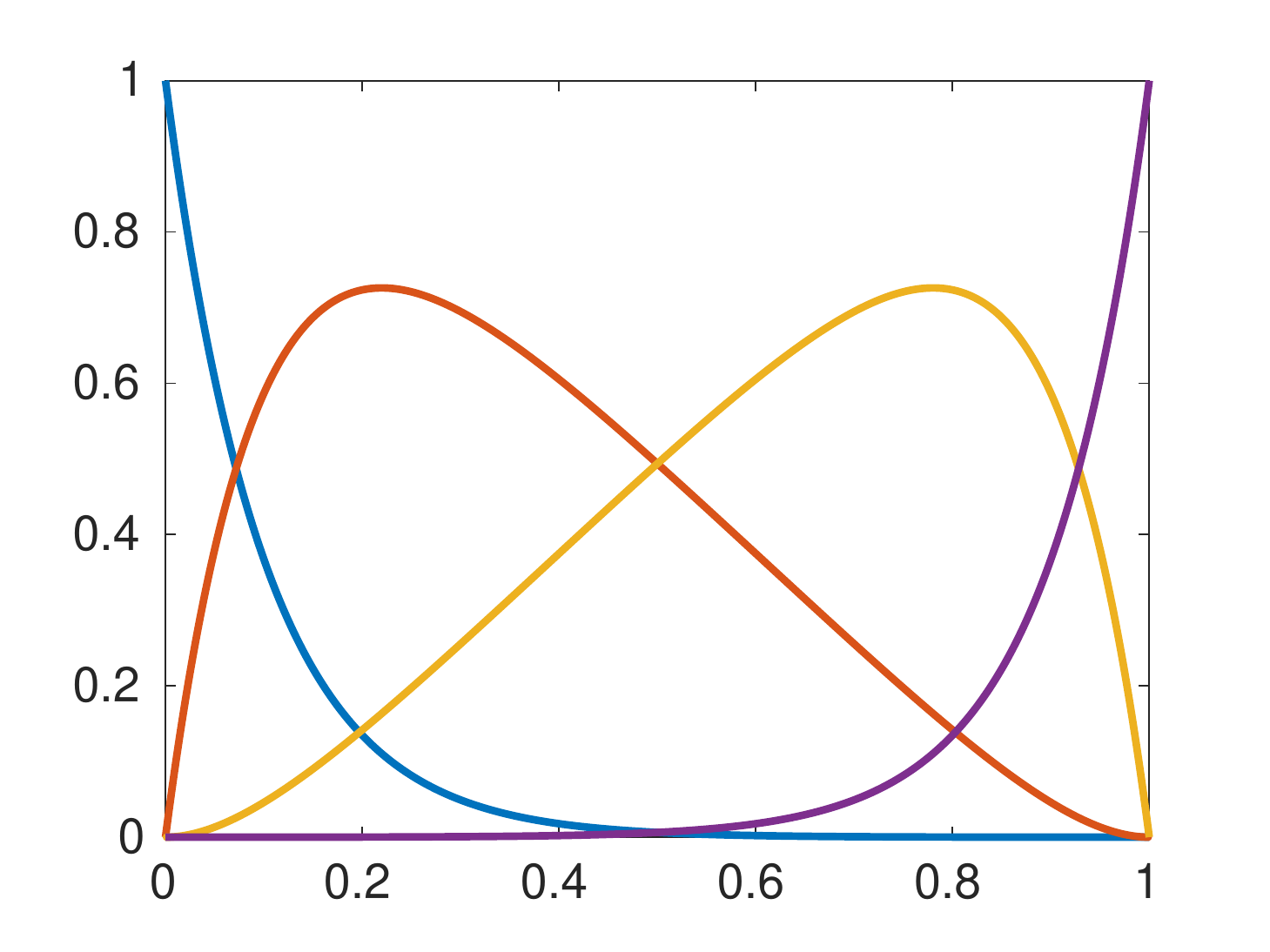} 
		\caption{$\spaceT_3^{(10,-10)} = \langle\, 1,x,e^{10x},e^{-10x}\, \rangle$} 
	\end{subfigure} 
	\\
	\begin{subfigure}[t]{0.4\linewidth}
		\centering
		\includegraphics[width=\linewidth]{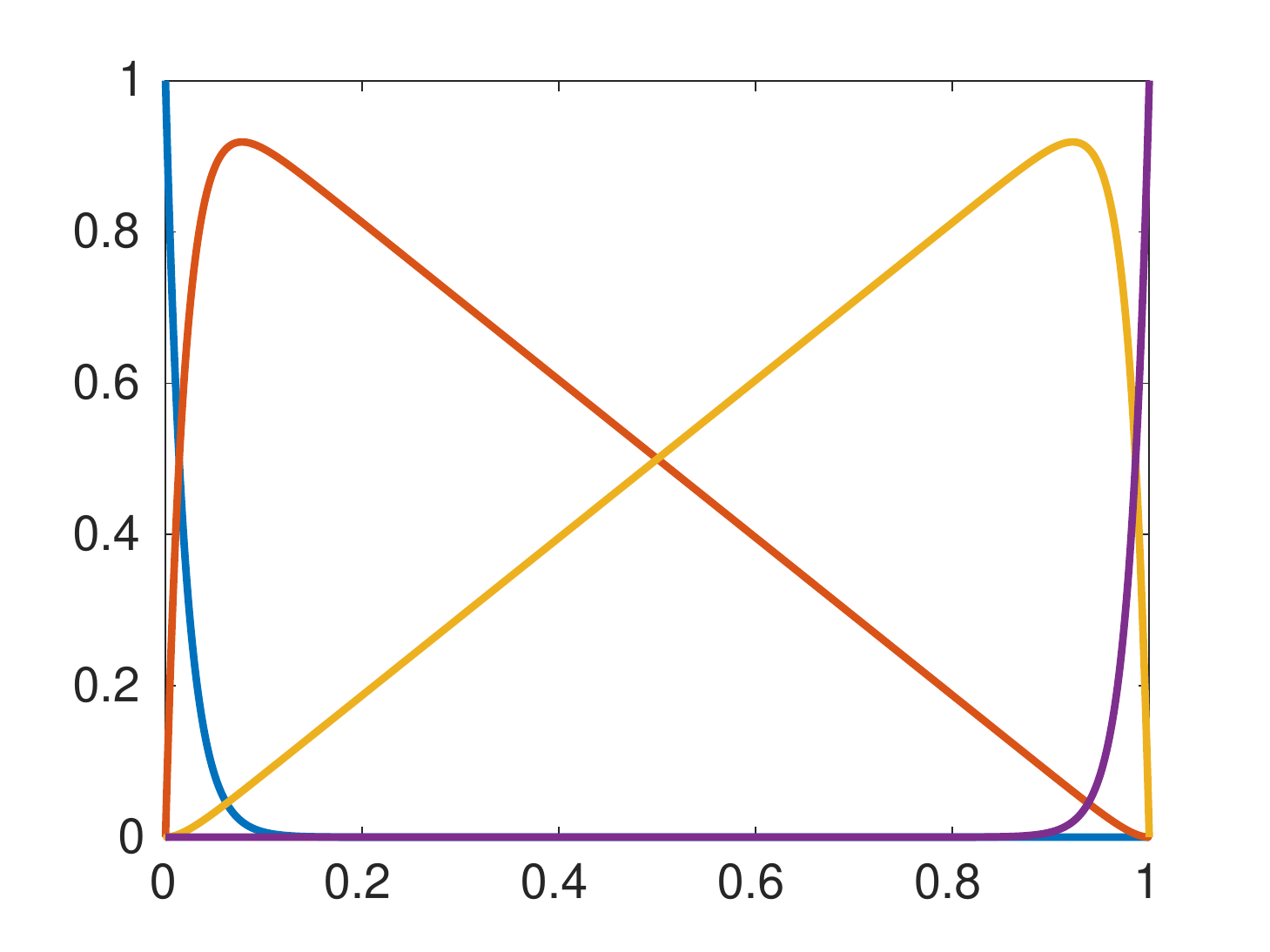} 
		\caption{$\spaceT_3^{(50,-50)} = \langle\, 1,x,e^{50x},e^{-50x}\, \rangle$} 
	\end{subfigure}
	\hspace{.75cm}
	\begin{subfigure}[t]{0.4\linewidth}
		\centering
		\includegraphics[width=\linewidth]{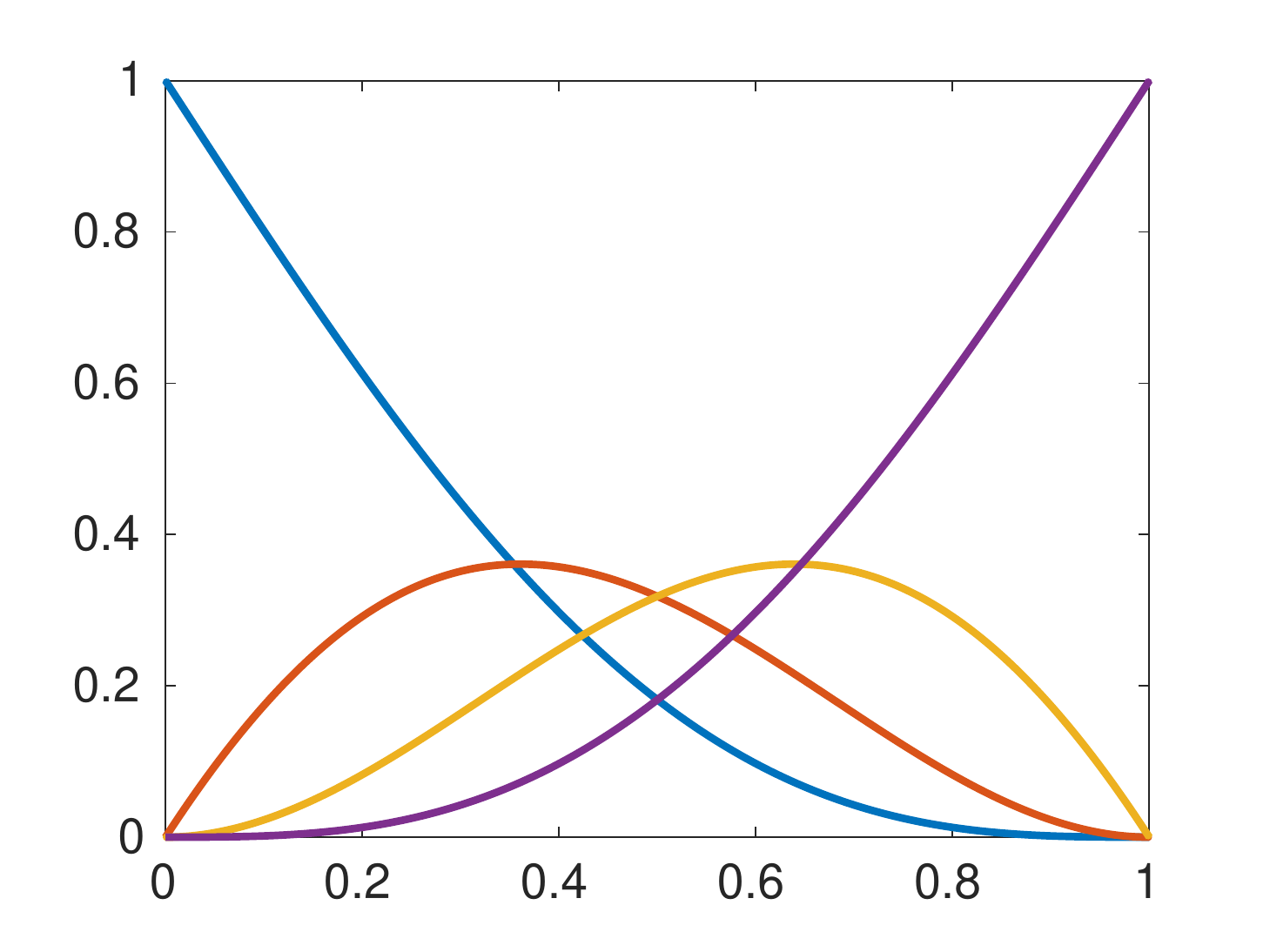} 
		\caption{$\spaceT_3^{(\myi\pi. -\myi\pi)} = \langle\, 1,x,\cos(\pi x),$ $ \sin(\pi x)\, \rangle$} 
	\end{subfigure} 
	\\
	\begin{subfigure}[t]{0.4\linewidth}
		\centering
		\includegraphics[width=\linewidth]{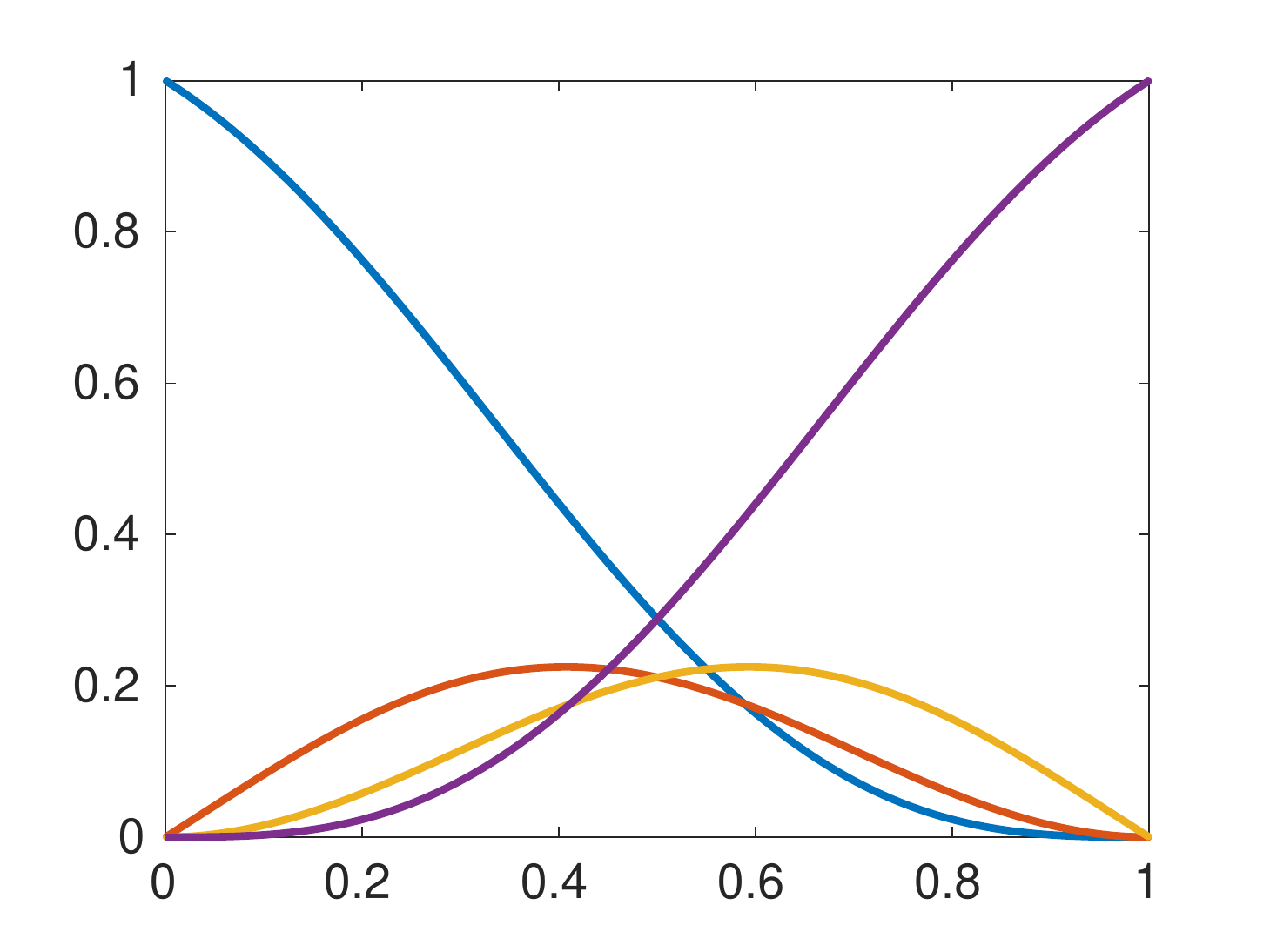} 
		\caption{$\spaceT_3^{(\myi \frac{3}{2}\pi. -\myi \frac{3}{2}\pi)} = \big\langle\, 1,x,\cos\Big(\frac{3\pi}{2} x\Big),\sin\Big(\frac{3\pi}{2} x\Big)\, \big\rangle$} 
	\end{subfigure} 
	\hspace{.75cm}
	\begin{subfigure}[t]{0.4\linewidth}
		\centering
		\includegraphics[width=\linewidth]{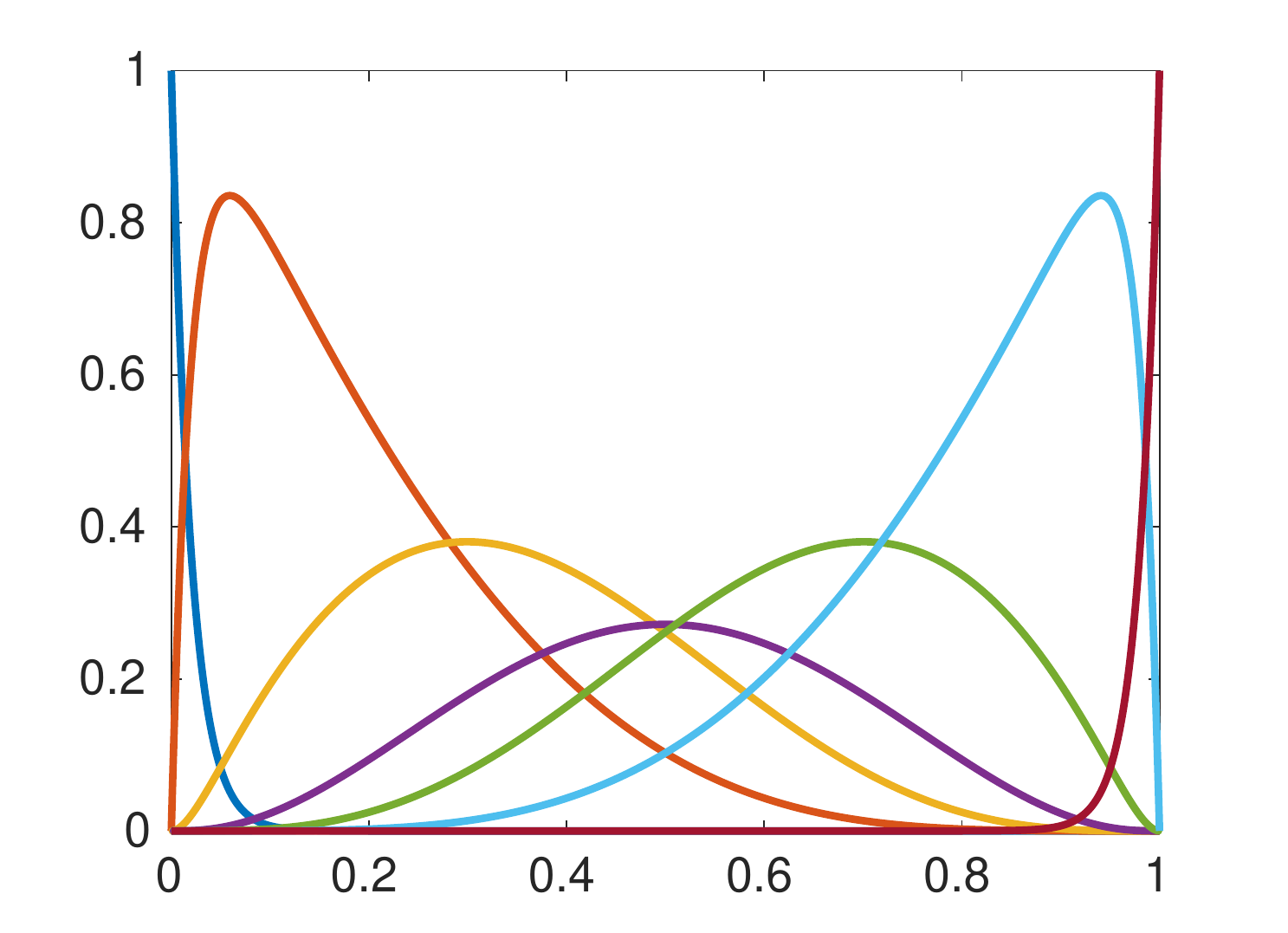} 
		\caption{\mbox{$\spaceT_6^{(50,-50,\myi \frac{3}{2}\pi. -\myi \frac{3}{2}\pi)} = \big\langle\, 1,x,x^2,e^{50x},e^{-50x},$  $\cos\Big(\frac{3\pi}{2} x\Big),\sin\Big(\frac{3\pi}{2} x\Big)\, \big\rangle$} }
	\end{subfigure} 
	\caption{Tchebycheffian Bernstein basis on the interval $[0,1]$ for different ECT-spaces.}
	\label{fig_brnstn}
\end{figure}
\begin{figure}[t!]
	\captionsetup[subfigure]{aboveskip=-.9cm,belowskip=0cm}
	\centering
	\begin{subfigure}[t]{0.49\linewidth}
		\centering
		\includegraphics[width=\linewidth]{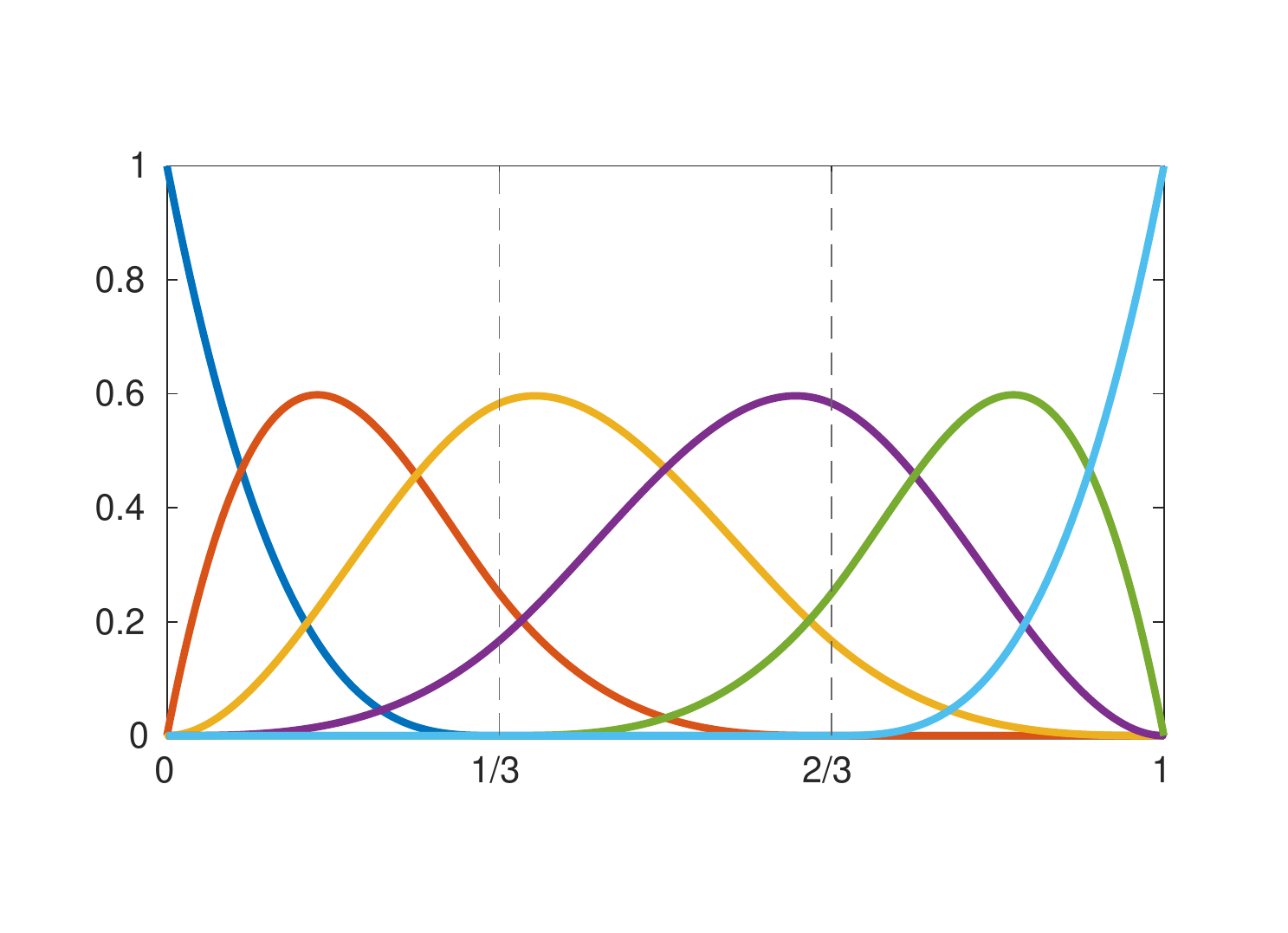} 
		\caption{$\spaceT_3 = \langle\, 1,x,x^2,x^3\, \rangle$} 
		\label{fig_sp_bas_pol} 
	\end{subfigure}
	\hfill
	\begin{subfigure}[t]{0.49\linewidth}
		\centering
		\includegraphics[width=\linewidth]{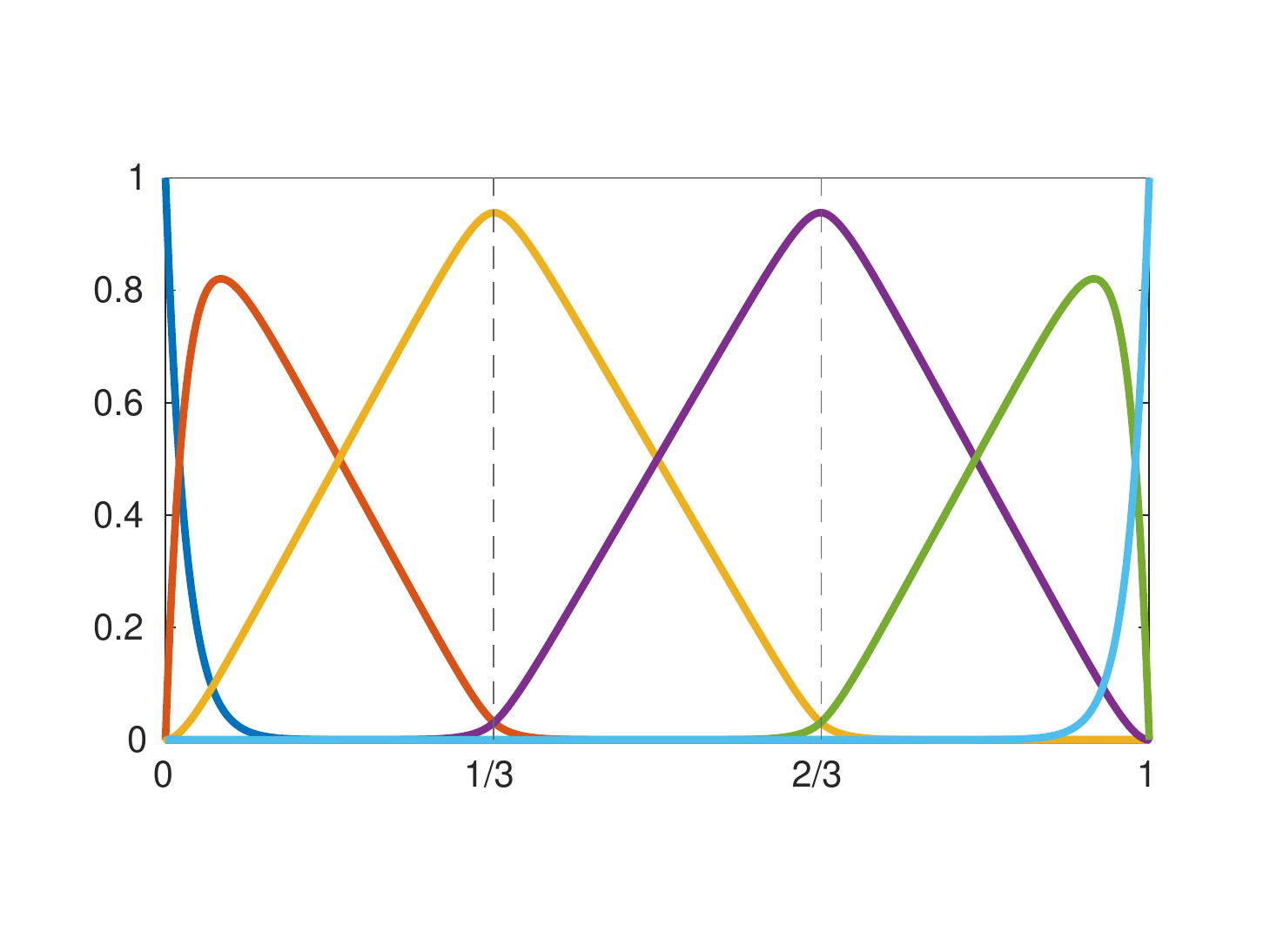} 
		\caption{$\spaceT_3^{(50,-50)} = \langle\, 1,x,e^{50x},e^{-50x}\, \rangle$} 
		\label{fig_sp_bas_exp} 
	\end{subfigure} 
	\hfill
	\vspace{-.4cm}
	\begin{subfigure}[t]{0.49\linewidth}
		\centering
		\includegraphics[width=\linewidth]{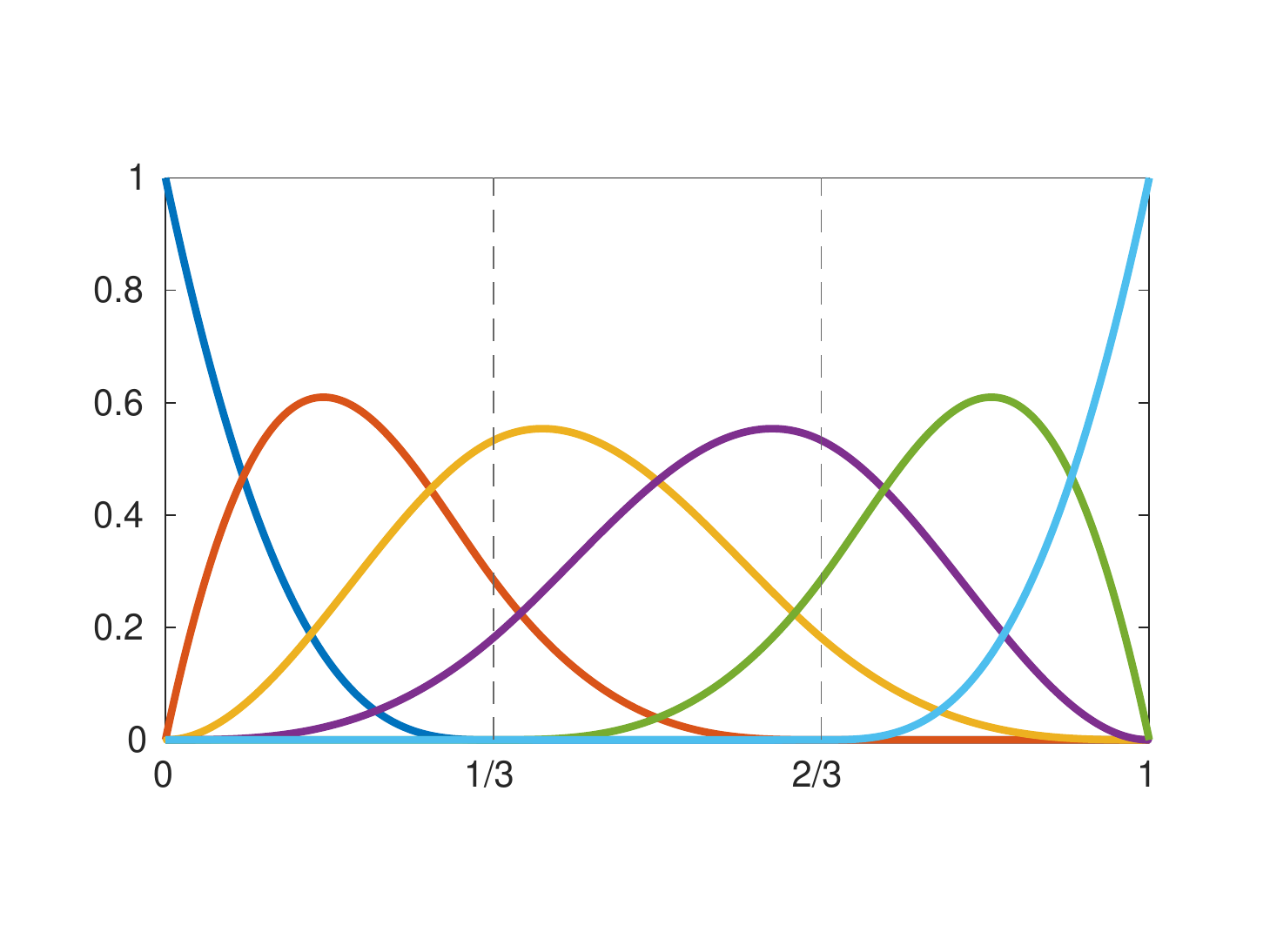} 
		\caption{$\spaceT_3^{(\myi \frac{3}{2}\pi. -\myi \frac{3}{2}\pi)} = \big\langle\, 1,x,\cos\Big(\frac{3\pi}{2} x\Big),$ $ \sin\Big(\frac{3\pi}{2} x\Big)\, \big\rangle$} 
		\label{fig_sp_bas_trig} 
	\end{subfigure}
	\hfill
	\begin{subfigure}[t]{0.49\linewidth}
		\centering
		\includegraphics[width=\linewidth]{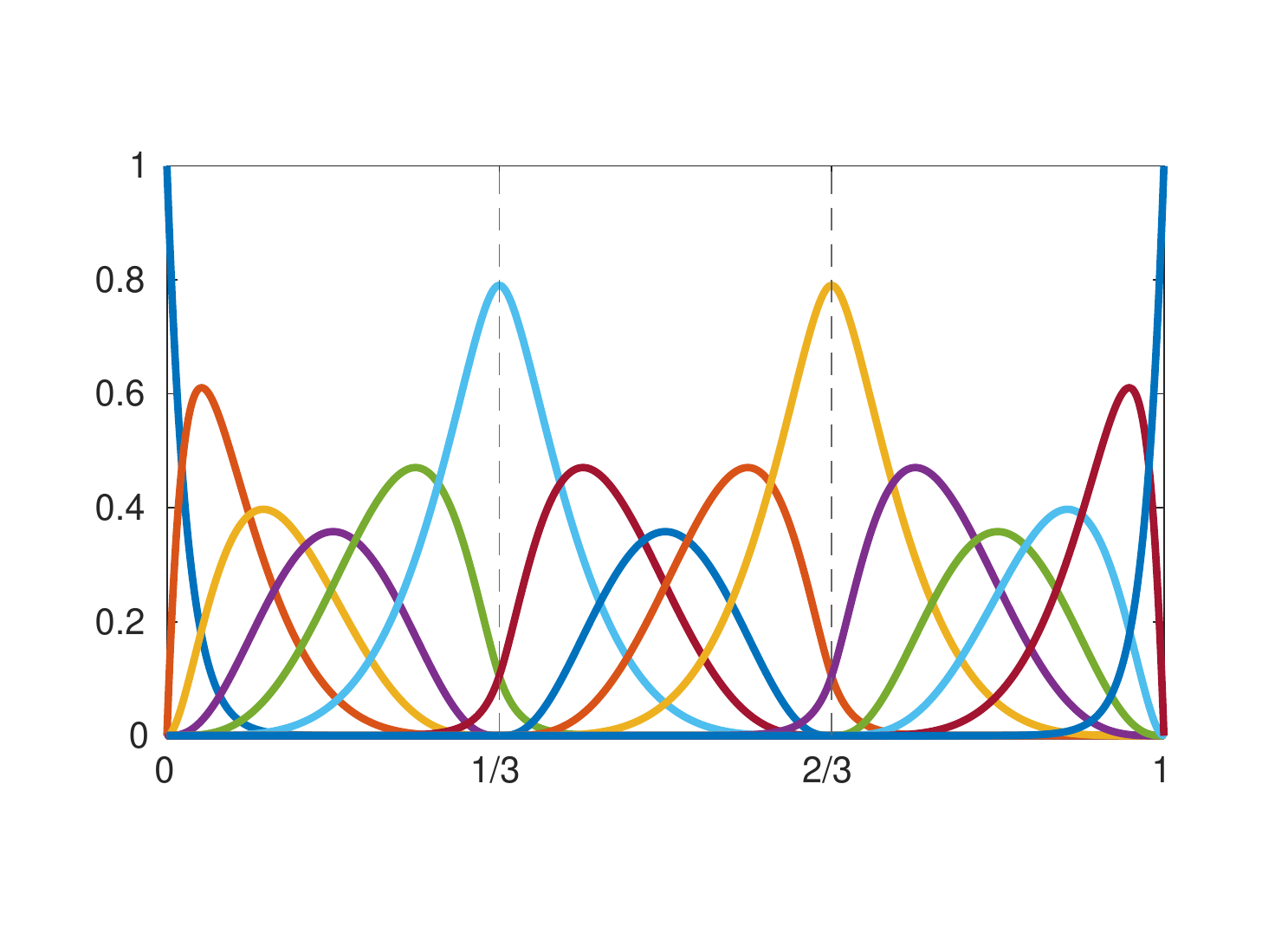} 
		\caption{$\spaceT_6^{(50,-50,\myi \frac{3}{2}\pi. -\myi \frac{3}{2}\pi)} = \big\langle\, 1,x,x^2,e^{50x},e^{-50x},\cos\Big(\frac{3\pi}{2} x\Big),$  $\sin\Big(\frac{3\pi}{2} x\Big)\, \big\rangle$} 
		\label{fig_sp_bas_mix} 
	\end{subfigure} 
	\caption{TB-splines on the partition $\Delta=\{0,1/3,2/3,1\}$ with $\smthVec=\{2,2\}$ for different local ECT-spaces.}
	\label{fig_tb_spline}
\end{figure}

When $m=1$, so that the partition $\Delta $ has no interior breakpoints, the space $\bbS_p^{\smthVec}(\Delta)$ reduces to the ECT-space $\bbT_p=\bbT_{p,1}$. The corresponding TB-splines are called Tchebycheffian Bernstein functions, and play the Tchebycheffian counterpart of Bernstein polynomials in case of classical algebraic polynomials. Tchebycheffian Bernstein functions provide an interesting alternative basis to the generalized power functions of the space $\bbT_p$. 

Some examples of the Tchebycheffian Bernstein basis for different ECT-spaces are depicted in \Cref{fig_brnstn}, and considering a partition consisting of three intervals, the corresponding $C^2$ TB-splines based on the same ECT-spaces are given in \Cref{fig_tb_spline}.

\Cref{def-TB-splines} allows for the use of different ECT-spaces on different intervals. However, the requirement of admissible weight systems states that the various ECT-spaces we are dealing with should be identified by a sequence of weights which connect smoothly across the different segments. In general, the construction of such weights is complicated; see \cite{beccari2019design, Mazure2018}. 
Some comments about the existence of a TB-spline basis of the space $\bbS_p^{\smthVec}(\Delta) $ are in order. 
\begin{itemize}
	\item Assume that $\bbT_p$ is a $(p+1)$-dimensional ECT-space on the interval $[a,b]$ and it contains constants. Moreover, assume that the derivative space of $\bbT_p$ is a $p$-dimensional ECT-space on $[a,b]$ generated by the weights $\{ w_1,\ldots,w_p\}$.
	If all the spaces $\bbT_{p,i}$, $i=1,\ldots, m$ in \eqref{eq_tcheb_sp} are selected as the restriction of $\bbT_p$, then from \Cref{rmk_gen-powers-der} it follows that the weights $\{1,w_1,\ldots,w_p\}$
	restricted to the intervals $[x_{i-1}, x_i]$, $i=1,\ldots, m$ form admissible weight systems for the space $\bbS_p^{\smthVec}(\Delta)$. 
	In other words, if all the pieces of the Tchebycheffian spline functions are taken from a single space $\bbT_p$ which is an ECT-space on $[a,b]$ containing constants, then the space $\bbS_p^{\smthVec}(\Delta)$ admits a TB-spline basis provided that the derivative space of $\bbT_p$ is also an ECT-space on $[a,b]$; see \cite{lyche2019tchebycheffian,Mazure2005,Schumaker2007} and references therein.
	\item Assume now that $\bbT_p$ is a $(p+1)$-dimensional ECT-space on each interval $[x_{i-1}, x_i]$ separately, but not necessarily on the entire interval $[a,b]$. 
	Moreover, assume that $\bbT_p$ is the null-space of a linear differential operator as in \eqref{eq_diff_op} and it contains constants.
	If all the pieces of the Tchebycheffian spline functions are taken from such space $\bbT_p$, then the space $\bbS_p^{\smthVec}(\Delta)$ admits a TB-spline basis on condition that 
	$h_\Delta \coloneqq \max_{i=1,\ldots,m}(x_{i}-x_{i-1})$ is sufficiently small. In particular, it admits a TB-spline basis when the derivative space of $\bbT_p$ is an ECT-space on each interval $[\xi_{k+1},\xi_{k+p}]$, where $\xi_{k+1}<\xi_{k+p}$ and $1\leq k\leq n$. We refer the reader to \cite{BrilleaudM2014} for details (and also \cite{LycheS1996}). Note that these intervals are closely related to the supports of the TB-splines; see Proposition~\ref{prop-TB-splines}.

	\item Admissible weight systems can be easily obtained for the interesting subclass of generalized  polynomial splines, i.e., splines where the pieces are taken from generalized polynomial spaces (such as \eqref{eq_gb_space_exp} and \eqref{eq_gb_space_trig}), even considering different spaces on different intervals. This ensures the existence of TB-splines, called generalized polynomial B-splines (GB-splines) in this setting, under very mild restrictions. In particular, the space $\bbS_p^{\smthVec}(\Delta)$ admits a TB-spline basis when the $(p-1)$-th order derivative space of $\bbT_{p,i}$ in \eqref{eq_tcheb_sp} is an ECT-space on $[x_{i-1},x_{i}]$, $i=1,\ldots, m$. We refer the reader to \cite[Section~4]{lyche2019tchebycheffian} for details.
\end{itemize}

\begin{example}\label{ex_TB_exp}
	Let $\bbS_p^{\smthVec}(\Delta)$ be a spline space as in \eqref{eq_tcheb_sp} where
	$$\bbT_{p,i}=\spaceT_p^{(\alpha_1,\ldots,\alpha_\ell)}=\angleSpace{1,x,\ldots,x^{p-\ell},e^{\alpha_1 x},\ldots,e^{\alpha_\ell x}}, \quad p\geq\ell\geq0,$$
	for $i=1,\ldots,m$.
	From \Cref{rem_null_space} we deduce that $\spaceT_p^{(\alpha_1,\ldots,\alpha_\ell)}$ is an ECT-space on the whole real line, and in particular on the interval $[a,b]$. Moreover, the derivative space of $\spaceT_p^{(\alpha_1,\ldots,\alpha_\ell)}$ is $\spaceT_{p-1}^{(\alpha_1,\ldots,\alpha_\ell)}$ (a space without constants when $p=\ell$), which is also an ECT-space on the whole real line.
	We conclude that the corresponding spline space $\bbS_p^{\smthVec}(\Delta)$ admits a TB-spline basis on any partition of the interval $[a,b]$, without any restrictions.
\end{example}
\begin{example}\label{ex_TB_null-space}
	Let $\bbS_p^{\smthVec}(\Delta)$ be a spline space as in \eqref{eq_tcheb_sp} where
	$$\bbT_{p,i} 
	=\angleSpace{1,x,\ldots,x^{p-\ell-2q},e^{\alpha_1x},\ldots,e^{\alpha_\ell x},\cos(\beta_1 x), \sin(\beta_1 x),\ldots,\cos(\beta_{q} x), \sin(\beta_{q} x) },$$
	for $i=1,\ldots,m$, i.e., the space considered in \eqref{eq_spaces}. 
	The case $q=0$ is covered in \Cref{ex_TB_exp}. Hence, let us here focus on the case $q>0$ and without loss of generality we assume $\beta_1>0,\ldots,\beta_q>0$.
	Since the pieces of the spline functions belong to the null-space of a linear differential operator and taking into account \Cref{rem_null_space}, we conclude that
	the space $\bbS_p^{\smthVec}(\Delta)$ admits a TB-spline basis on the interval $[a,b]$ if the knot sequence $\bm{\xi}$ in \eqref{eq_knot_vec} satisfies
	$$
	\xi_{k+p}-\xi_{k+1}<\frac{\pi}{\max\bigl\{\beta_1,\ldots,\beta_q\bigr\}}, \quad k=1,\ldots,n.
	$$
\end{example}
\begin{example}\label{ex_ppoly}
	Let $\bbS_p^{\smthVec}(\Delta)$ be a spline space as in \eqref{eq_tcheb_sp} where
    $$\bbT_{p,i}=\angleSpace{1,\cos(\beta x), \sin(\beta x), \ldots, \cos(q\beta x), \sin(q\beta x) },
    \quad \beta>0, \quad p=2q\geq0, $$
    for $i=1,\ldots,m$. As discussed in the previous example, 
	such spline space admits a TB-spline basis on the interval $[a,b]$ if the knot sequence $\bm{\xi}$ in \eqref{eq_knot_vec} satisfies $\xi_{k+p}-\xi_{k+1}<{\pi}/{(q\beta)}$, $k=1,\ldots,n$.
    Actually, this bound can be improved: the space $\bbS_{p}^{\smthVec}(\Delta)$ admits a TB-spline basis if 
    $\xi_{k+p}-\xi_{k+1}<{\pi}/{\beta}$, $k=1,\ldots,n$;
    see \cite{Mazure2005,Sanchez-R-1998}.
    \end{example} 
\begin{example}\label{ex_TB_trig} 
	Let $\bbS_p^{\smthVec}(\Delta)$ be a spline space as in \eqref{eq_tcheb_sp} where
	$$\bbT_{p,i}=\angleSpace{1,x,\ldots, x^{p-2},\cos(\beta x), \sin(\beta x)}, \quad \beta>0, \quad p\geq 2,$$
	for $i=1,\ldots,m$. As discussed in \Cref{ex_TB_null-space}, 
	such spline space admits a TB-spline basis on the interval $[a,b]$ if the knot sequence $\bm{\xi}$ in \eqref{eq_knot_vec} satisfies $\xi_{k+p}-\xi_{k+1}<{\pi}/{\beta}$, $k=1,\ldots,n$; see also \Cref{ex_gb_trig_w}. Actually, it is a generalized polynomial spline space and it can be shown to have a GB-spline basis if $x_i-x_{i-1}<{\pi}/{\beta}$, $i=1,\ldots,m$; see \cite[Section~4]{lyche2019tchebycheffian}. 
\end{example}

In the rest of the paper we will only consider Tchebycheffian spline spaces as in \eqref{eq_tcheb_sp} obtained by considering the same ECT-space $\spaceT_p^\rootVec$, defined in \eqref{eq_spaces}, on each interval of the partition $\Delta$. In view of the above discussion, this setting ensures the existence of a TB-spline basis whenever the derivative space of $\spaceT_p^\rootVec$ is an ECT-space on $[a,b]$. If this is not the case, a TB-spline basis can be obtained when a sufficiently fine partition of the given interval is considered. 
For notational convenience, we denote the corresponding TB-splines as
$$
N_{p,k}^\rootVec, \quad k=1,\ldots, n,
$$
in order to stress the selected ECT-space.
As numerically illustrated in \Cref{sec-numerics}, this simplified Tchebycheffian setting allows already for such a wide variety of ECT-spaces and related shape parameters to be extremely profitable in the context of IgA.

\begin{remark}
The recurrence relation in \Cref{def-TB-splines} is the main theoretical tool to define and analyze TB-splines. 
 However, it is of limited practical interest because it does not lead to an efficient evaluation strategy. 
 Fortunately, this deficiency has been recently overcome: efficient routines to construct and manipulate TB-splines with pieces belonging to ECT-spaces as in \Cref{sec_null_spc} have been provided in \cite{speleers2022algorithm} by implementing the approach described in \cite{hiemstra2020tchebycheffian}. This approach starts from the Tchebycheffian Bernstein functions on each interval $[x_{i-1},x_i]$, $i=1,\ldots,m$, and, by using the smoothness constraints at all breakpoints, generates an extraction matrix expressing the TB-spline basis locally in terms of the Tchebycheffian Bernstein functions. The latter can be computed directly through Hermite interpolation (see \Cref{ex_TB_trig2} for an illustration). No weights are required along the process.
Although the evaluation of TB-splines for extreme values of the shape parameters (i.e., the roots of the characteristic polynomial \eqref{eq_char_pol}), high degrees, high smoothness, and highly non-uniform partitions still faces numerical challenges \cite{speleers2022algorithm}, the above mentioned routines provide reliable tools for a large class of TB-splines of interest in practical applications. A Matlab implementation is publicly available through the ACM Collected Algorithms (CALGO) library \cite {speleers2022algorithm}.
\end{remark}
 \begin{remark} \label{rem-no-weights}
  The TB-spline construction procedure developed in \cite{hiemstra2020tchebycheffian,speleers2022algorithm} does not require any knowledge of weights. In case the space $\bbS_p^{\smthVec}(\Delta)$ does not possess admissible weight systems, the procedure simply returns a sequence of functions that might not enjoy the properties listed in \Cref{prop-TB-splines}. In this perspective, the computational process can also be used in practice to numerically test the (non-)existence of a TB-spline basis for a given Tchebycheffian spline space.
\end{remark}

\begin{figure}[t!]
	\centering	
		\includegraphics[width=.45\linewidth]{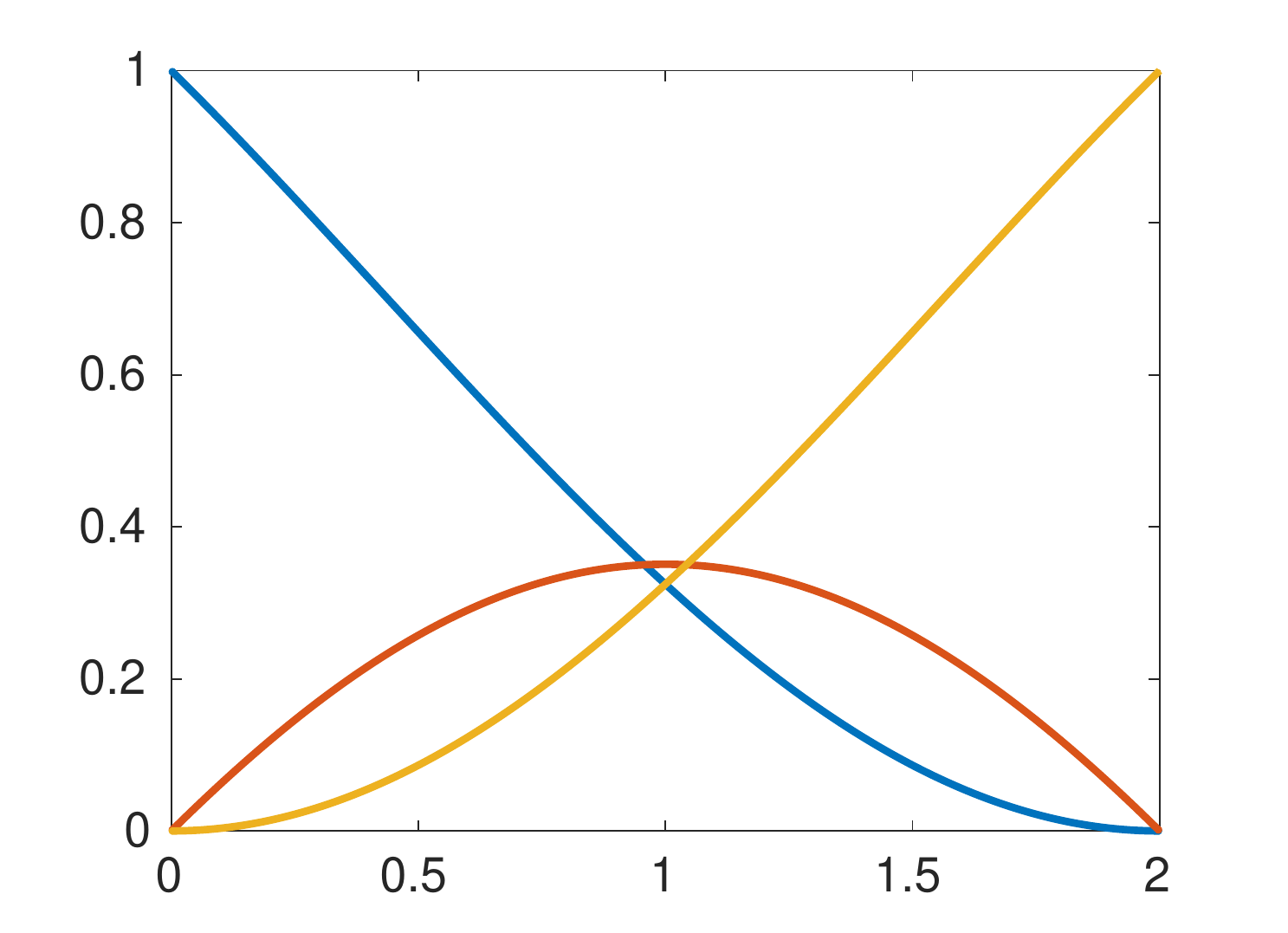} 
	\hspace{.75cm}
		\includegraphics[width=.45\linewidth]{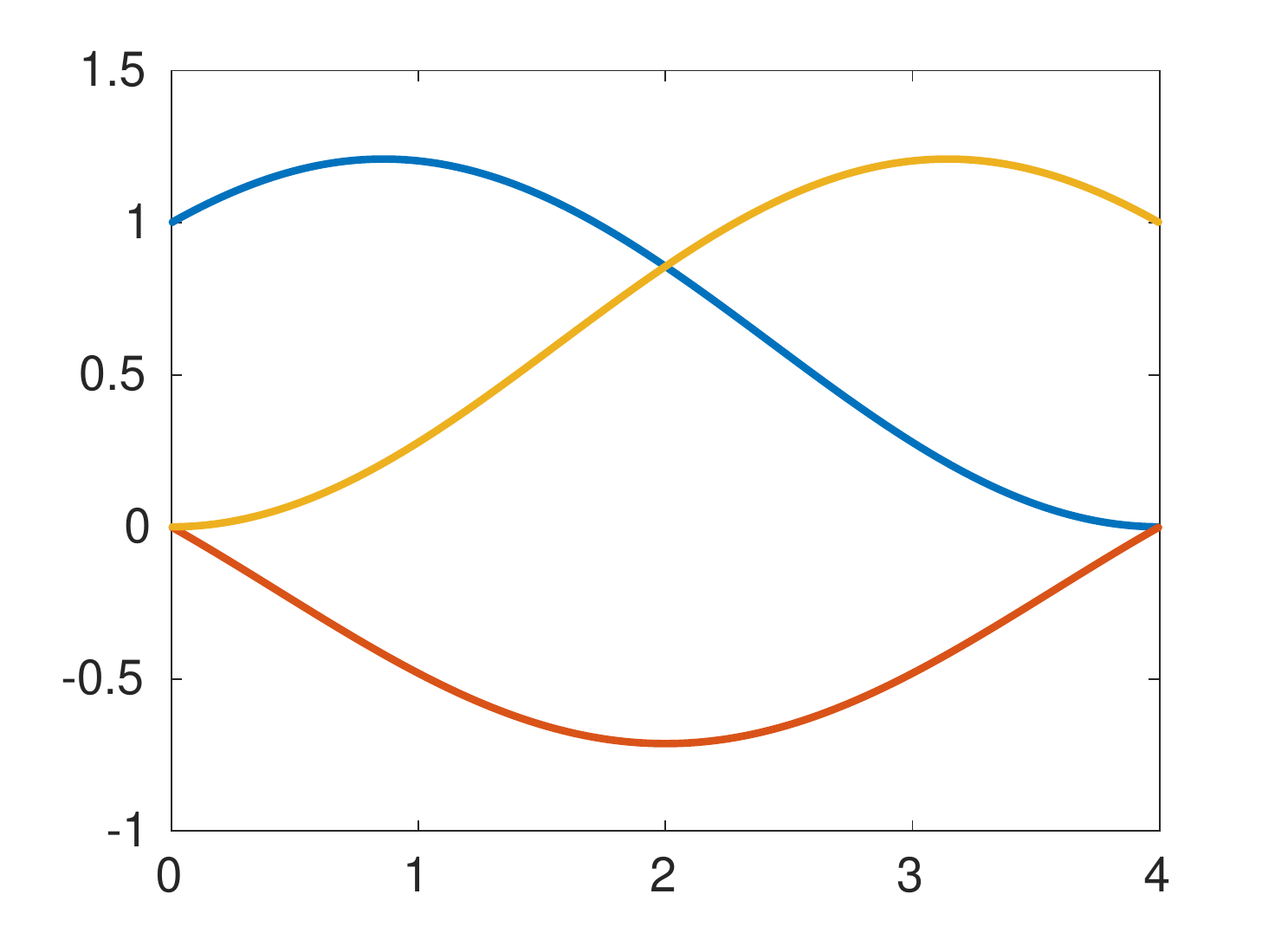} 	
	\caption{The functions satisfying \eqref{eq:B0}--\eqref{eq:B2} in the space $\spaceT_2^{(\myi,-\myi)}=\angleSpace{1,\cos(x),\sin(x)}$. 
	Left: interval $[0,2]$, Right: interval $[0,4]$.}
	\label{fig-Bernstein}
\end{figure}

\begin{example}
	\label{ex_TB_trig2} 
	The space
	$$\spaceT_2^{(\myi,-\myi)} = \angleSpace{1,\cos(x), \sin(x)}$$ 
	is an ECT-space of dimension $3$ on any interval of length less than $2\pi$; see \Cref{ex_gb_trig_w}.
	However, the space only admits a Tchebycheffian Bernstein basis on intervals of smaller length, namely less than $\pi$; see \Cref{ex_ppoly,ex_TB_trig}. If such a basis exists, then the corresponding basis elements $B_{j,2}$, $j=0,1,2$, can be computed by solving the following Hermite interpolation problems:
	\begin{align}
	 B_{0,2}(a) &=1,\quad B_{0,2}(b)=0,\quad  DB_{0,2}(b)=0; \label{eq:B0}\\
	 B_{1,2}(a) &=0,\quad D B_{1,2}(a) = - D B_{0,2}(a), \quad  B_{1,2}(b)=0; \label{eq:B1}\\
	 B_{2,2}(a) &=0,\quad  DB_{2,2}(a)=0,\quad B_{2,2}(b) =1. \label{eq:B2}
	\end{align}
	These interpolation problems have a unique solution whenever $b-a<2\pi$, since $\spaceT_2^{(\myi,-\myi)}$ is an ECT-space then. However, the unisolvency of the above interpolation problems does not ensure that the corresponding functions $B_{j,2}$, $j=0,1,2$ enjoy all the properties listed in \Cref{prop-TB-splines}, in particular their non-negativity.  This is illustrated in \Cref{fig-Bernstein} for the interval $[a,b] = [0,2]$ (so $b-a < \pi$) and the interval $[a,b] = [0,4]$ (so $\pi < b-a < 2\pi$).
\end{example}

%


\section{Isogeometric analysis with Tchebycheffian B-splines}
\label{sec-TB-IGA}
In this section we outline the isogeometric Galerkin method based on TB-splines.
For the sake of simplicity, we focus on second-order elliptic differential problems with homogeneous Dirichlet boundary conditions.

\subsection{Isogeometric Tchebycheffian Galerkin methods}\label{sec-Galerkin}

Let ${\mathfrak L}$ be a linear second-order elliptic differential operator on the domain
$\Omega\subset \bbR^d$ with Lipschitz boundary $\partial \Omega$.
We consider the differential problem
\begin{equation}\label{eq:PDE}
\begin{cases}
{\mathfrak L} u={\rm f}, & \text{in\ } \Omega,\\
u=0, & \text{on\ }  \partial\Omega,
\end{cases}
\end{equation}
whose weak form reads as follow: 
\begin{equation}\label{eq:weak-formulation}
\text{find } u\in\bbV \quad\text{ such that } \quad  a(u,v)=F(v), \quad \forall v\in \bbV,
\end{equation}
where $\bbV$ is a suitable function space,
$
a:\bbV\times \bbV\rightarrow\bbR $ is a bilinear form depending on $ {\mathfrak L}$,
and $
F:\bbV\rightarrow\bbR $ is a linear form depending on ${\rm f}$.

The Galerkin approach to approximate the solution of \eqref{eq:PDE} is based on the weak form \eqref{eq:weak-formulation}. We select a finite-dimensional approximation space on $\Omega$,
\begin{equation}\label{eq:Vh}
\bbW \coloneqq \bigl\langle \varphi_1,\varphi_2,\ldots,\varphi_{n_\bbW}\bigr\rangle\subset \bbV,
\quad \dim(\bbW)=n_{\bbW},
\end{equation}
and we look for 
\begin{equation*}
u_{\bbW}\in \bbW \quad \text { such that }\quad a(u_{\bbW},w)=F(w), \quad \forall w\in \bbW.
\end{equation*}
Taking
$$u_{\bbW}=\sum_{i=1}^{n_\bbW} c_{i}\varphi_i$$
gives rise to a linear system
$A\mathbf{c}=\mathbf{F}$, where the matrix $A$ and the vector $\mathbf{F}$ are defined as 
\begin{equation}\label{eq-matrix-rhs}
A_{i,j} \coloneqq a(\varphi_j,\varphi_i), \quad i,j=1,\ldots,n_\bbW,
\qquad
F_i \coloneqq F(\varphi_i),\quad i=1,\ldots,n_\bbW.
\end{equation}
Different Galerkin methods correspond to different choices of the subspace $\bbW$.

\begin{figure}[t!]
	\centering
	\includegraphics[width=.9\linewidth]{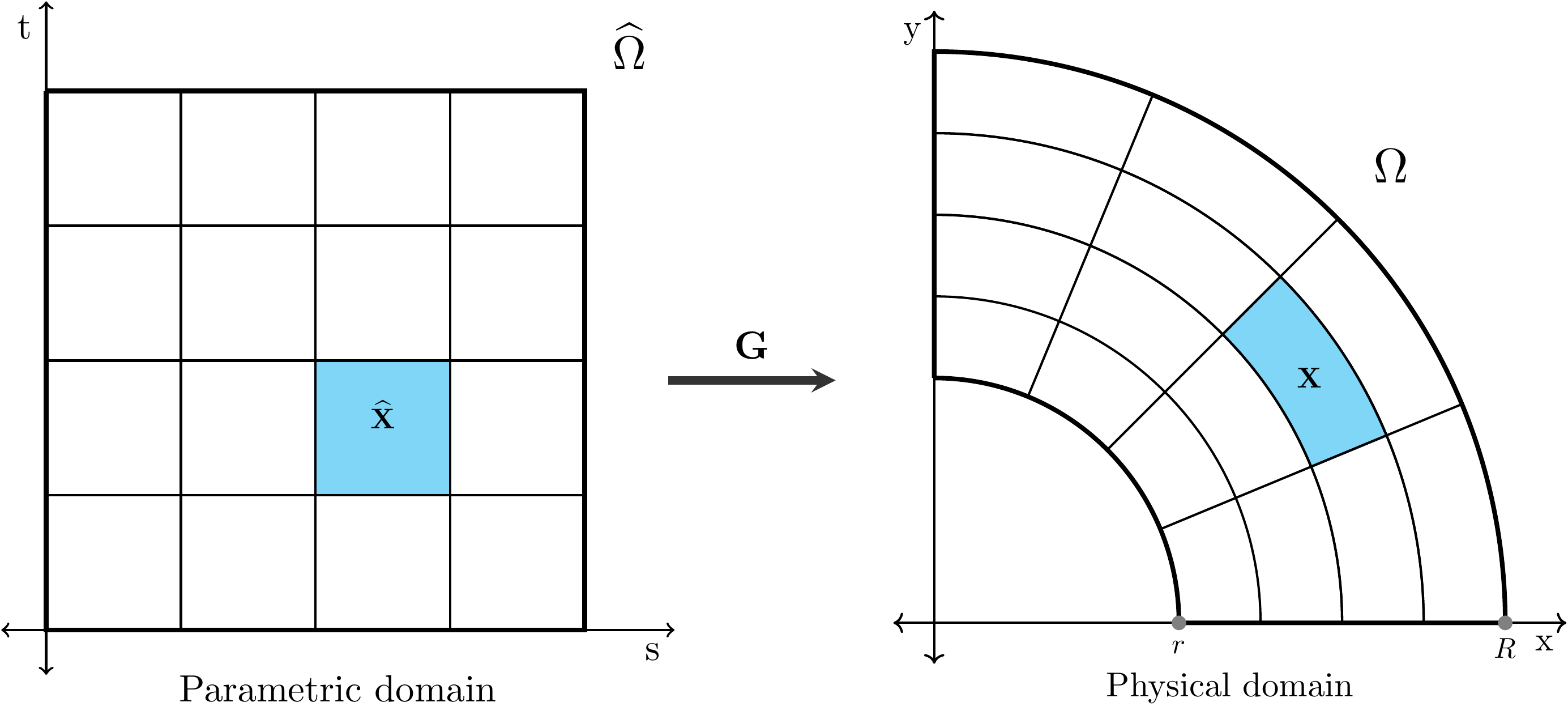} 
	\caption{Isogeometric Galerkin method. Description of the physical domain $\Omega$ by means of a global geometry map~$\map$.}
	\label{fig-IgA-Galerkin}
\end{figure}
%
In the standard formulation of IgA, the physical domain $\Omega$ is represented by means of an {invertible global geometry map}. We define the geometry map $\map$ from the closure of the parametric domain $\widehat \Omega \coloneqq (0,1)^d$ to the closure of the physical domain $\Omega$ as (see Figure~\ref{fig-IgA-Galerkin})
\begin{equation}\label{geometry-map}
\map(\hat\xx) \coloneqq \sum_{i=1}^{n_{\map}}{\mathbf{P}}_i  \hat \varphi_i(\hat\xx),
\quad {\mathbf{P}}_i\in \bbR^d, 
\end{equation}
where the basis functions
\begin{equation} \label{basis}
\bigl\{\hat\varphi_1,\ldots,\hat\varphi_{n_\map}\bigr\}
\end{equation}
have to be selected so as to produce an exact (or at least very accurate) representation of the geometry (and to allow imposition of homogeneous boundary conditions).
Following the isoparametric approach, the fields of interest are described with the same basis functions as the geometry in IgA. The space $\bbW$ in \eqref{eq:Vh}, which incorporates the homogeneous boundary conditions, is then spanned by the isogeometric functions
\begin{equation}
\label{basis-phi}
\varphi_i(\xx) \coloneqq \hat\varphi_{j_i}\circ \map^{-1}(\xx)=\hat\varphi_{j_i}(\hat\xx),\quad i=1,
\ldots, n_\bbW, \quad \ \xx=\map(\hat\xx),
\end{equation}
where $n_\bbW< n_\map$ and $j_i:\{1,\dots,n_\map\}\rightarrow\{1,\dots, n_\map\}$ are a suitable rearrangement of the indices.
In the most common isogeometric formulation, the functions in \eqref{basis} are chosen to be tensor-product B-splines or NURBS. Nevertheless, B-splines/NURBS are not a requisite ingredient in the isogeometric paradigm and several alternatives have been presented in the literature.

The properties described in \Cref{prop-TB-splines} show that TB-splines are plug-to-plug compatible with classical (polynomial) B-splines, hence they can be used as an alternative basis on the parametric domain for building isogeometric discretization spaces according to \eqref{basis-phi}. In this paper, we focus on TB-splines from ECT-spaces of the form \eqref{eq_spaces}. 

We now provide some notation that will be useful in the case studies reported in \Cref{sec-numerics}. We deal with two-dimensional problems, where the bivariate TB-splines on the parametric domain are constructed by taking the tensor product of univariate TB-splines considering a single ECT-space on the interval $[0,1]$. Let $\Delta_1 \times \Delta_2$ be a rectangular grid in $[0,1]^2$, 
let $\smthVec_{\Delta_1},\smthVec_{\Delta_2}$ be two smoothness vectors at the corresponding breakpoints, let $p_1,p_2\geq1$ be two degrees, and let $\rootVec_1,\rootVec_2$ be two root vectors as in \eqref{eq_root_vec}. Then, the tensor-product Tchebycheffian spline space is given by
\begin{equation}  \label{eq_space_tensor}
\Big\langle N_{p_1,k}^{\rootVec_1}\, N_{p_2,l}^{\rootVec_2} : k=1,\dots,n_1,\ l=1,\dots,n_2 \Big\rangle,
\end{equation}
where $n_1$ and $n_2$ are the dimensions of the two univariate Tchebycheffian spline spaces.
For notational convenience, the basis functions are also shortly indexed as
\begin{equation*}
N_i \coloneqq N_{p_1,k}^{\rootVec_1}\, N_{p_2,l}^{\rootVec_2},\quad i\coloneqq (l-1)n_1+k, \quad k=1,\dots,n_1, \quad l=1,\ldots,n_2.
\end{equation*}
A point in (the closure of) the physical domain $\Omega$ is denoted as $\mathbf{x} = (x,y)$ and a point in (the closure of) the parametric domain $\widehat\Omega$ as $\hat{\mathbf{x}} = (s,t)$. The geometry map in \eqref{geometry-map} is given by
\begin{equation} \label{eq_geo_map} 
\mathbf{x} = \map(\hat{\mathbf{x}}) = \sum_{i=1}^{n_\map} \mathbf{P}_iN_i(\hat{\mathbf{x}}),\quad \mathbf{P}_i\in \bbR^2, \quad n_\map=n_1n_2,
\end{equation}
while the basis functions in \eqref{basis-phi} are
$$ 
\phi_i(\mathbf{x}) \coloneqq N_{j_i} \circ \map^{-1}(\mathbf{x}) = N_{j_i}(\hat{\mathbf{x}}),
$$
where, in order to incorporate the homogeneous boundary conditions, the tensor-product TB-splines  
$$ 
N_{p_1,k}^{\rootVec_1}\, N_{p_2,l}^{\rootVec_2},\quad k=1,n_1, \quad l=1,\dots,n_2, 
\qquad N_{p_1,k}^{\rootVec_1}\, N_{p_2,l}^{\rootVec_2},\quad k=1,\dots,n_1, \quad l=1,n_2,
$$
are excluded according to the last property in \Cref{prop-TB-splines}. Thus we have $n_\bbW=(n_1-2)(n_2-2)$.

\subsection{Selection of the ECT-spaces}
As mentioned above, we are interested in testing the performance of isogeometric Galerkin methods based on TB-splines from ECT-spaces of the form \eqref{eq_spaces}.
It is clear that such a class provides a large variety of combinations of polynomial, exponential, and trigonometric functions equipped with a wide spectrum of shape parameters. Of course, it is crucial to properly choose the structure of the space (including different kinds of functions) and the related shape parameters to fully exploit the functionality of ECT-spaces. 
Particularly, there are two aspects that should be considered in such a selection to profit from a non-polynomial structure in the reference ECT-spaces.
\begin{itemize}
\item The first is the geometric aspect of the problem based on the geometry mapping at hand. Often, the boundary of the physical domain $\Omega$ is built from (arcs of) conic sections. Therefore, it can be exactly described in terms of polynomials, exponential, and trigonometric functions with suitable shape parameters. In the pure isogeometric philosophy, the geometry of the physical domain is exact from the beginning of the procedure and does not change along the process when possible refinements are required. In this perspective, it is then natural to select ECT-spaces \eqref{eq_spaces} that allow for an exact representation of the geometry.  
In its classical formulation, IgA is based on B-splines or NURBS. Even though (arcs of) conic sections can be exactly represented in terms of NURBS, it is not an arc-length parameterization and smooth representations of (closed) conic sections require high degrees. As an example, a $C^r$ NURBS parameterization of a circle requires at least degree $2(r+1)$; see \cite{Bangert1997}. Moreover, NURBS behave poorly when dealing with differentiation and integration. On the other hand, spaces of the form \eqref{eq_spaces} allow for natural smooth parameterizations of conic sections and their ``structure'' does not change when dealing with differentiation and integration, in a complete similarity with polynomial spaces; see also \Cref{rem_ect_benefits}.
\item The second, of interest in the numerical treatment of PDEs, is the fact that functions in ECT-spaces of the form \eqref{eq_spaces} are fundamental solutions of important ordinary differential operators, meaning that the shape parameters are immediately connected to the associated differential problems.
Although the extension to the multivariate case is not straightforward, the use of TB-splines with pieces drawn from suitable ECT-spaces can still be beneficial. As an example, they can lead to a significant reduction, or even a complete removal, of extraneous oscillations in advection-dominated problems, without the need of stabilization techniques. Once again, a proper selection of the shape parameters is crucial but can be driven by the problem setting. 
\end{itemize}

In the next section we present some case studies that show how TB-splines can outperform classical polynomial B-splines in isogeometric Galerkin discretizations of Poisson and advection-diffusion problems. As discussed above, the selection of the underlying ECT-spaces has to be done according to a problem-driven strategy; this will be detailed in each of the case studies.

\section{Case studies}
\label{sec-numerics}

In this section we discuss the performance of the isogeometric Galerkin method based on TB-splines through different case studies.
In all the numerical experiments we use the following setup. 

Bivariate tensor-product Tchebycheffian spline spaces of the form \eqref{eq_space_tensor} are taken as discretization spaces in the isogeometric Galerkin method.
The corresponding ECT-spaces in the two directions will be specified in each experiment by writing
$$
\spaceT_{p_1}^{\rootVec_1}\otimes\spaceT_{p_2}^{\rootVec_2},
$$ 
where $p_1,p_2$ are the degrees and $\rootVec_1, \rootVec_2$ are the root vectors as in \eqref{eq_root_vec}. 
We only consider uniform partitions $\Delta_1=\Delta_2$ with breakpoints 
$$x_i=\frac{i}{m}, \quad i=0,\dots, m,$$
and denote by $h$ the (uniform) distance between two consecutive breakpoints, i.e.,
$h=\frac{1}{m}$.
We assume maximal smoothness at the breakpoints, i.e.,
$$\smthVec_{\Delta_1}= \{p_1-1,\dots, p_1-1\}, \quad\smthVec_{\Delta_2}=\{p_2-1,\dots, p_2-1\}.$$ 
This amounts to univariate Tchebycheffian spline spaces of dimension $ n_1=m+p_1$ and $n_2=m+p_2$, respectively.

We focus on the general advection-diffusion problem of the form
\begin{equation}\label{eq_adv_diff} 
\begin{split}  \begin{cases}
-\nabla \cdot ({\kappa} \nabla u) + \textbf{a} \cdot \nabla u = {\rm f},\quad &\text{in} \ \Omega = \map\big((0,1)^2\big),\\
u = {\rm g},\quad &\text{in}\ \partial \Omega,
\end{cases} \end{split}
\end{equation}
where $\mathbf{a}$ is the 
advection flow velocity, 
$\kappa$ the diffusivity, ${\rm f}$ the prescribed source function, and ${\rm g}$   the Dirichlet boundary data. The physical domain $\Omega$ is described by the geometry map $\map$. 
The dominance of advection over diffusion is defined by the global \peclet number 
\begin{equation}\label{eq_peclet}
\bm{Pe}_g \coloneqq \frac{\|\bm{a}\|
}{\kappa}.
\end{equation}
For simplicity, we fix the diffusion coefficient $\kappa=1$ in all the case studies.

As usual, homogeneous boundary conditions are satisfied pointwise exactly. In the non-homogeneous case, the boundary values are represented in the underlying Tchebycheffian spline space by suitable approximation strategies (e.g., least squares approximation) and subsequently the reduction to the homogeneous case is considered, so dealing again with a special instance of the problem \eqref{eq:PDE}.

For an accurate numerical approximation of the integrals needed in the construction of the matrix and the vector in \eqref{eq-matrix-rhs}, we use elementwise Gaussian quadrature rules. The selection of quadrature points for discretizations based on TB-splines is not trivial and  quadrature rules of higher order compared to the algebraic polynomial splines are often required. In the presented case studies we use $3p$ quadrature points in each element for TB-splines only involving trigonometric functions (see \Cref{sec_Example1,sec_Example2}), and $5p$ quadrature points in each element when also exponential functions with large shape parameters are considered (see \Cref{sec_Example3,sec_Example4,sec_Example5}).

\subsection{Case study 1: Poisson problem on a domain bounded by circular arcs}
\label{sec_Example1}

\begin{figure}[t!]
\captionsetup[subfigure]{aboveskip=-1pt,belowskip=0cm}
\centering
 \begin{subfigure}[t]{0.55\linewidth}
    \centering
    \includegraphics[width=\linewidth]{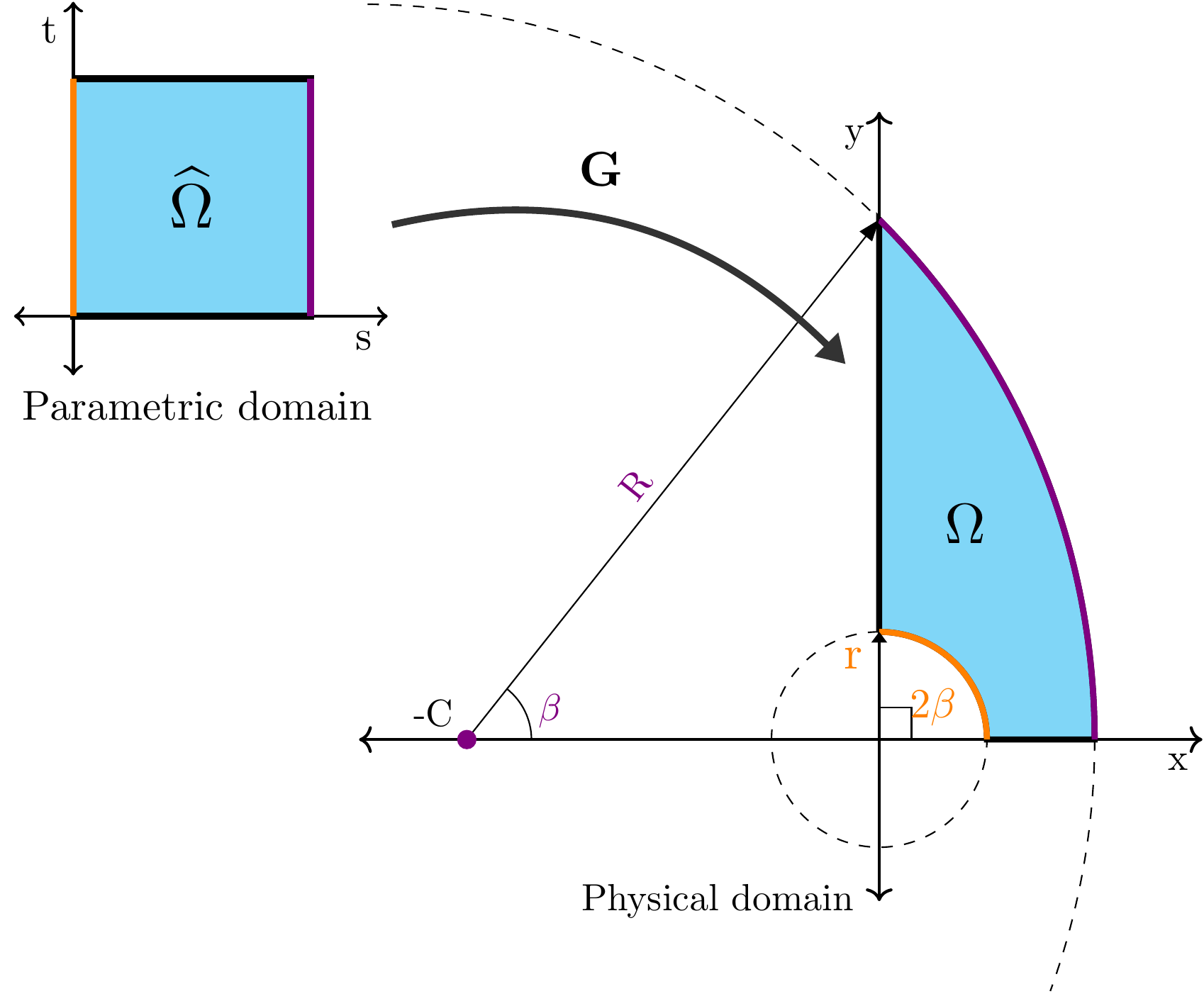} 
  \end{subfigure}
 \hfill 
   \begin{subfigure}[t]{0.35\linewidth}
    \hspace{-2cm} \centering
    \includegraphics[width=\linewidth]{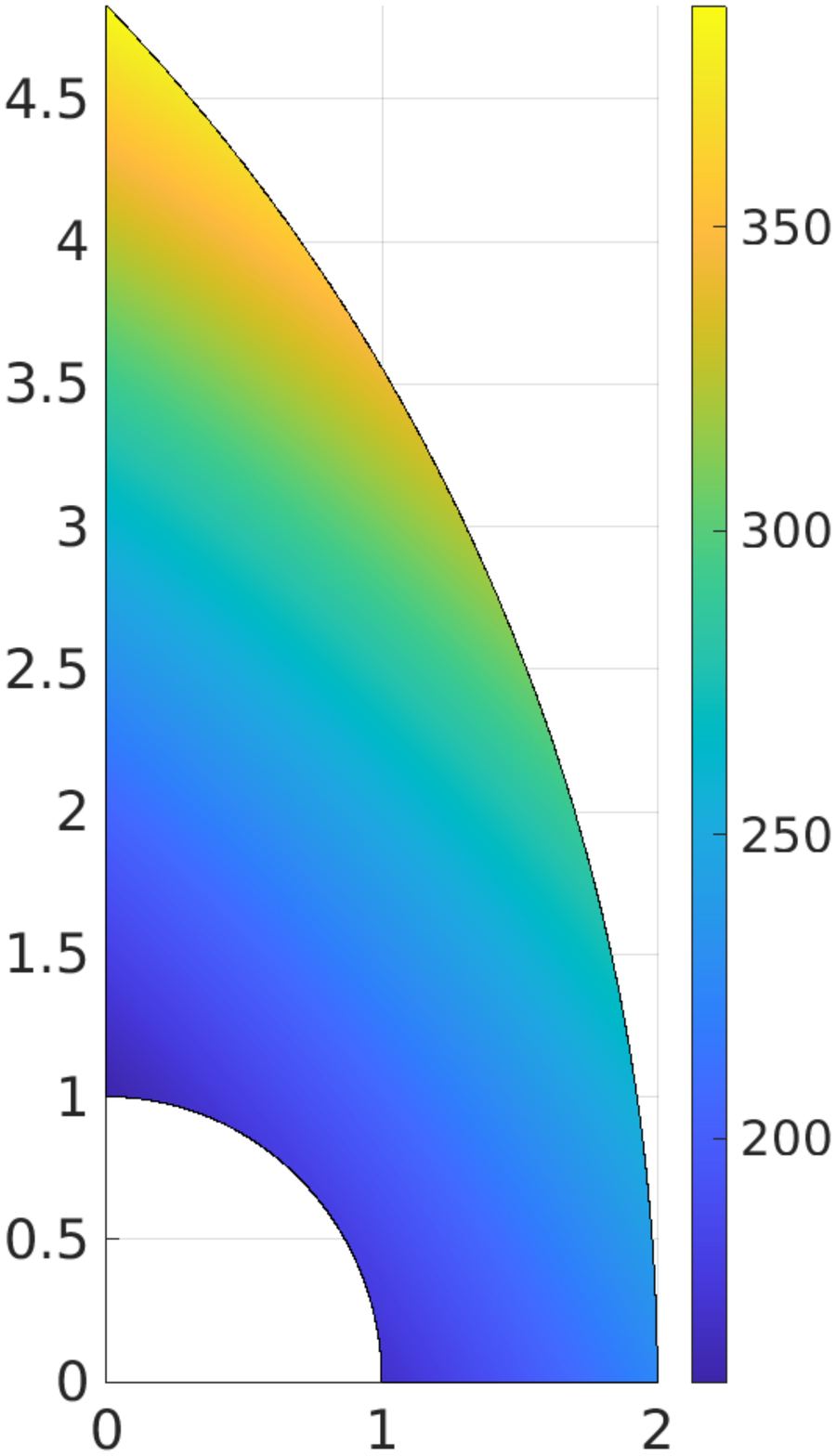} 
  \end{subfigure} 
  \caption{Case study 1. Left: Geometry map presented in \eqref{eq_map_ecce_annulus} as a combination of curves with phases $\beta$ and $2\beta$. Right: Plot of the approximate solution obtained by using tensor-product TB-splines identified by $\xheightsub{\spaceT}{4} \otimes \spaceT_{4}^{(\myi\frac{\pi}{4},-\myi\frac{\pi}{4}, \myi \frac{\pi}{2}, -\myi \frac{\pi}{2})}$ with $m=8$ and dof $=144$.}
  \label{fig_ex_ecce_annulus} 
\end{figure}

\begin{figure}[t!] 
\centering
\includegraphics[width=0.5\linewidth]{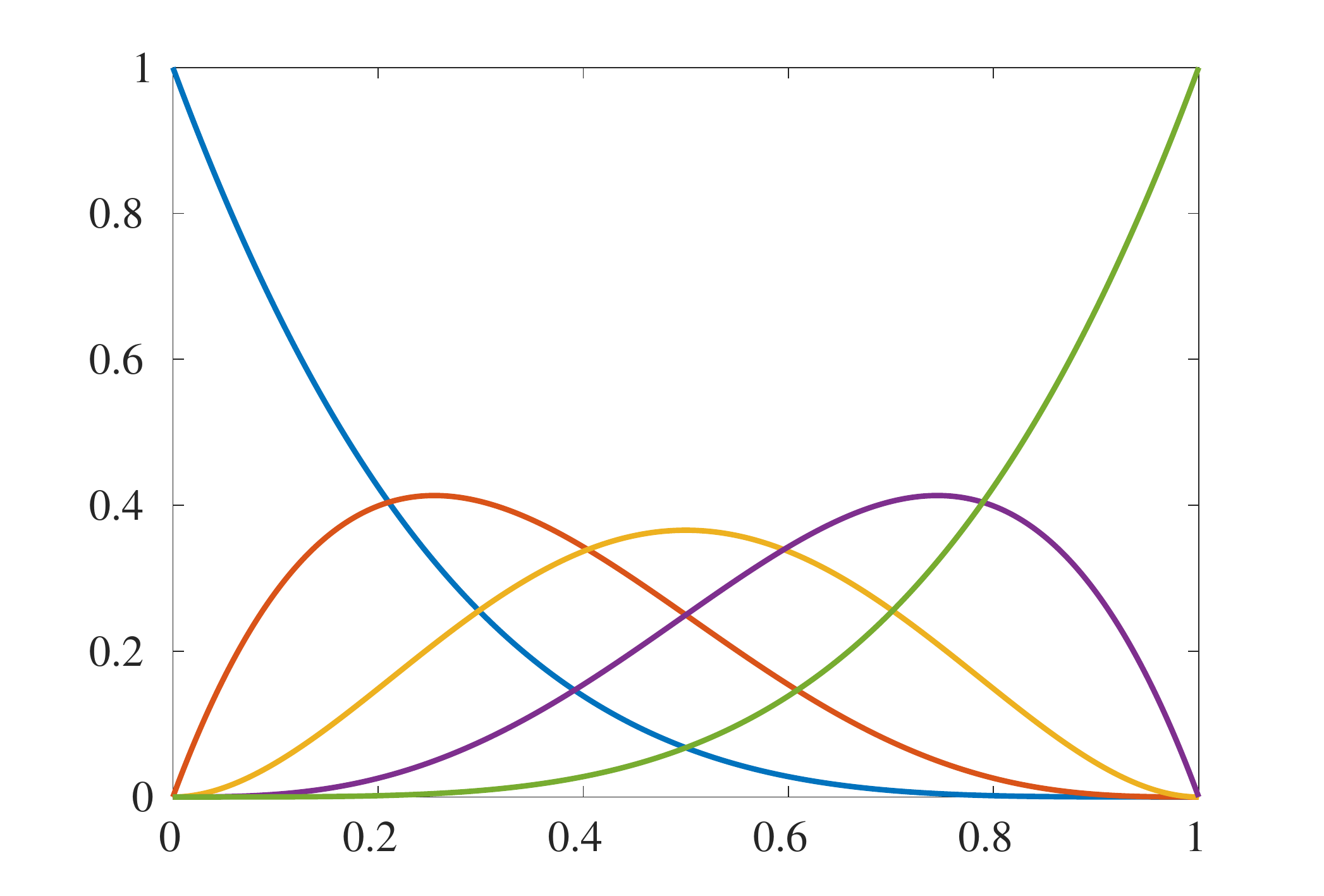}
\caption{Case study 1. Tchebycheffian Bernstein basis of $\spaceT_4^{(\myi\beta,-\myi\beta, \myi 2\beta, -\myi 2\beta)}$ on the interval $[0,1]$.}
\label{fig_ex1_benstn} 
\end{figure}

In this case study we address a Poisson problem by considering \eqref{eq_adv_diff} with $\kappa=1$, $\textbf{a}=\textbf{0}$,
and ${\rm f}$ obtained from the exact solution 
$$u(x,y) = 5 ((x + 4)^2 + (y + 3)^2 ) + xy.$$
The parametric domain $\widehat\Omega= (0, 1)^2$ is mapped to the physical domain $\Omega$ through the geometry map given by
\begin{equation} \label{eq_map_ecce_annulus} 
\begin{pmatrix} x\\ y \end{pmatrix} = \textbf{G}
\begin{pmatrix} s\\ t \end{pmatrix} = (1 - s)r \,
\begin{pmatrix} \cos(2\beta t)\\ \sin(2\beta t) \end{pmatrix} +\ s\, 
\begin{pmatrix}
-C + R\cos(\beta t)\\ R\sin(\beta t)\end{pmatrix}, \end{equation}
with
$$r = 1,\quad \beta = \frac{\pi}{4},\quad C = \frac{2r\cos(\beta)}{1-\cos(\beta)},\quad R = C+2r,$$
which leads to the domain depicted in \Cref{fig_ex_ecce_annulus}, left. 
It is evident that the geometry can be exactly represented in the form \eqref{eq_geo_map} by considering in the parametric $t$-direction Tchebycheffian splines identified by an ECT-space that contains the trigonometric functions with phases $\beta$ and $2\beta$. Hence we consider in the $t$-direction the space of degree $p_2 = 4$, 
\begin{equation*}
\spaceT_4^{(\myi\beta,-\myi\beta, \myi 2\beta, -\myi 2\beta)} = 
\Biggl\langle\, 1,\cos\biggl(\frac{\pi}{4}t\biggr), \sin\biggl(\frac{\pi}{4}t\biggr), \cos\biggl(\frac{\pi}{2}t\biggr), \sin\biggl(\frac{\pi}{2}t\biggr) \, \Biggr\rangle.
\end{equation*}
The Tchebycheffian Bernstein basis of this ECT-space is illustrated in \Cref{fig_ex1_benstn}. 
From \Cref{ex_ppoly} we know that the corresponding Tchebycheffian spline space admits a TB-spline basis on any partition of the interval $[0, 1]$ since $\beta<\pi$.
In the $s$-direction, we simply consider algebraic polynomials of degree $p_1 = 4$, i.e.,
$$\xheightsub{\spaceT}{4} = \angleSpace{ 1,s,s^2,s^3,s^4 }.$$
This gives the tensor-product space 
$$\xheightsub{\spaceT}{4} \otimes \spaceT_4^{(\myi\frac{\pi}{4},-\myi\frac{\pi}{4}, \myi \frac{\pi}{2}, -\myi \frac{\pi}{2})}. $$ 
To obtain the control coefficients $\mathbf{P}_i$ of the geometry in the form \eqref{eq_geo_map} we simply do interpolation on the grid obtained by the tensor product of the Greville points of polynomial B-splines of the same degree $4$. We remark that this geometry cannot be exactly represented by using generalized polynomial B-splines (see \cite{manni2011generalized}) with pieces in ECT-spaces of the form \eqref{eq_gb_space_trig}.

We approximate the $L^\infty$ norm of the error by sampling the approximate and exact solutions on a uniform grid in the parametric domain consisting of $501$ points along each direction.
 \Cref{tab_ecce_annulus} shows the $L^\infty$ error obtained for different levels of refinements, starting with a very coarse mesh of only one interval in each parametric direction, comparing the TB-splines against the classical B-splines for spaces of the same dimension. We see an improvement of about a factor $10$ that we achieve just by exactly representing the geometry using TB-splines. The contour plot of the approximate solution  partitioning the domain in $8\times 8$ elements is illustrated in \Cref{fig_ex_ecce_annulus}, right.  

\begin{table}[t!]
\centering
\renewcommand{\arraystretch}{1.1}
  \begin{tabular} { >{\centering}p{1cm}  >{\centering}p{1cm} >{\centering}p{3.5cm} >{\centering\arraybackslash}p{2.5cm}}
  \hline
$m$ & dof & $ \xheightsub{\spaceT}{4} \otimes \spaceT_{4}^{(\myi\frac{\pi}{4},-\myi\frac{\pi}{4}, \myi \frac{\pi}{2}, -\myi \frac{\pi}{2})}$ & $\spaceT_{4} \otimes \spaceT_{4}$\\
\hline \hline
1        & 25     & 7.9294e-04      & 5.0629e-03 \\
2        & 36     & 3.4854e-04      & 5.9308e-04 \\
4        & 64     & 2.9655e-05      & 8.4760e-05 \\
8        & 144    & 8.3556e-07      & 2.7812e-06 \\
16       & 400    & 2.6336e-08      & 9.2305e-08 \\
32       & 1296   & 8.4881e-10      & 2.9782e-09 \\ 
\hline
\end{tabular}\caption{Case study 1. Comparison of the $L^\infty$ error of tensor-product TB-splines and B-splines of degrees $p_1=p_2=4$ for different numbers of intervals $m$ in each direction.} 
\renewcommand{\arraystretch}{1.0}
\label{tab_ecce_annulus}
\end{table}
%

\subsection{Case study 2: Poisson problem on a symmetric domain bounded by circular arcs}
\label{sec_Example2}
\begin{figure}[t!]
\captionsetup[subfigure]{aboveskip=-1pt,belowskip=0cm}
\centering
  \hspace{8mm}
 \begin{subfigure}[t]{0.42\linewidth}
    \centering
    \includegraphics[width=\textwidth]{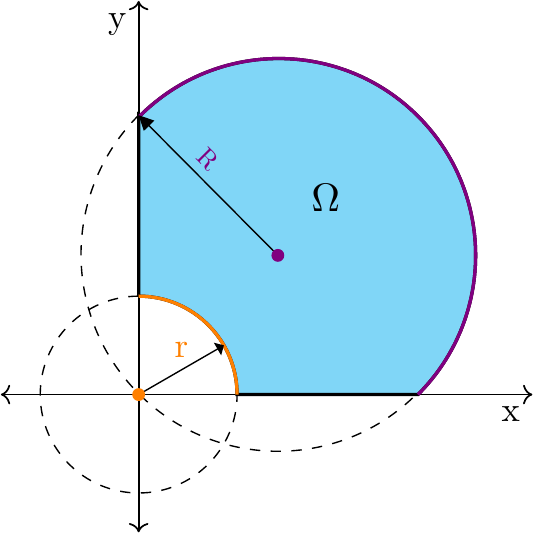}
  \end{subfigure}
 \hfill
   \begin{subfigure}[t]{0.35\linewidth}
    \centering
    \includegraphics[width=\textwidth]{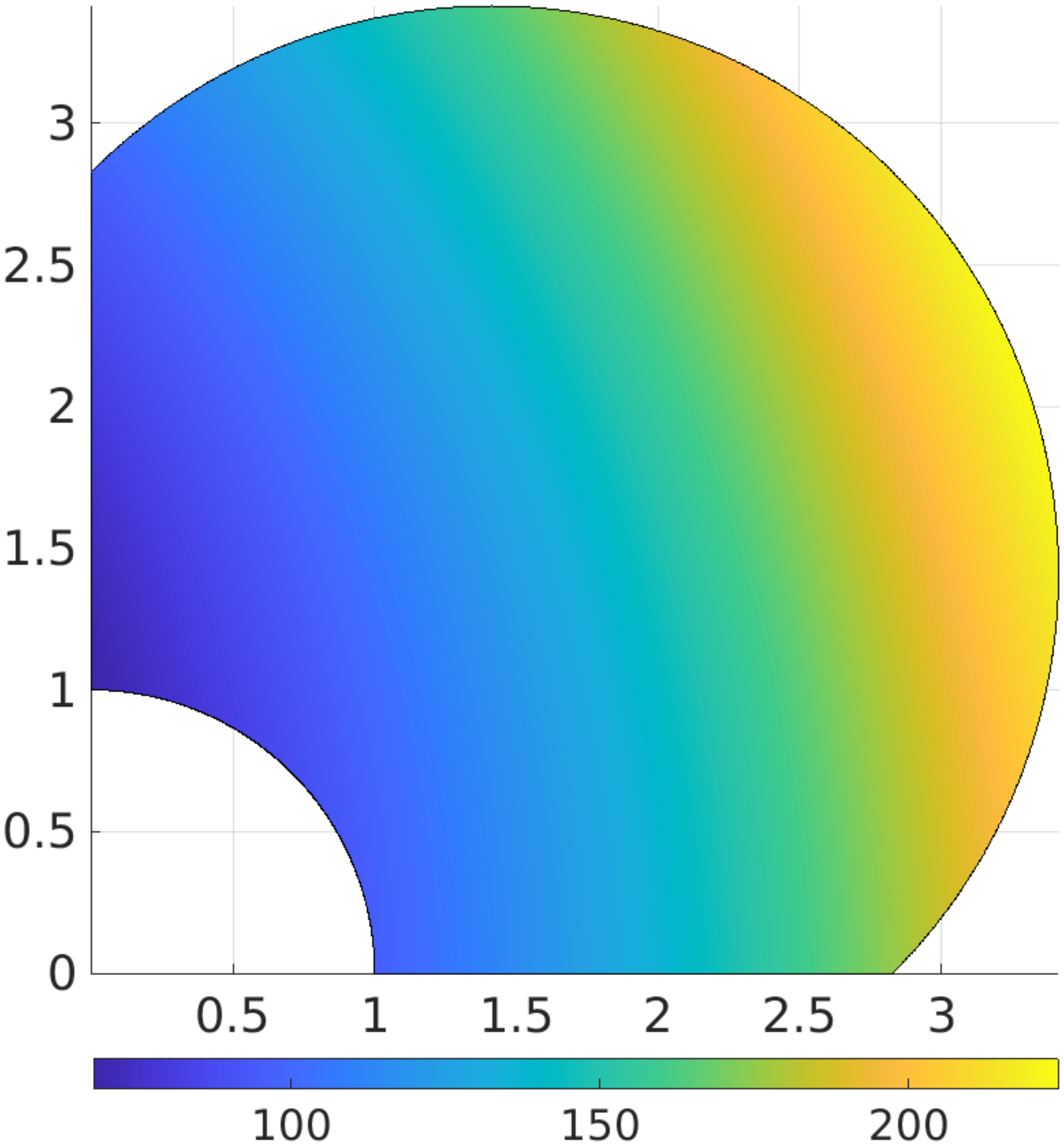}
  \end{subfigure} 
  \hspace{8mm}
  \caption{Case study 2. Left: Geometry map presented in \eqref{eq_map_mickey} as a combination of curves with phases $\beta$ and $2\beta$. Right: Plot of the approximate solution obtained by using tensor-product TB-splines identified by $\xheightsub{\spaceT}{6}\otimes\spaceT_6^{(\myi \beta,-\myi \beta, \myi 2\beta, -\myi 2\beta,\myi{3\beta},-\myi{3\beta})} $, $\beta=\frac{\pi}{2}$ with $m=8$ and dof $=196$.}
  \label{fig_ex_mickey} 
\end{figure}

\begin{figure}[t!] \centering
\includegraphics[width=0.5\linewidth]{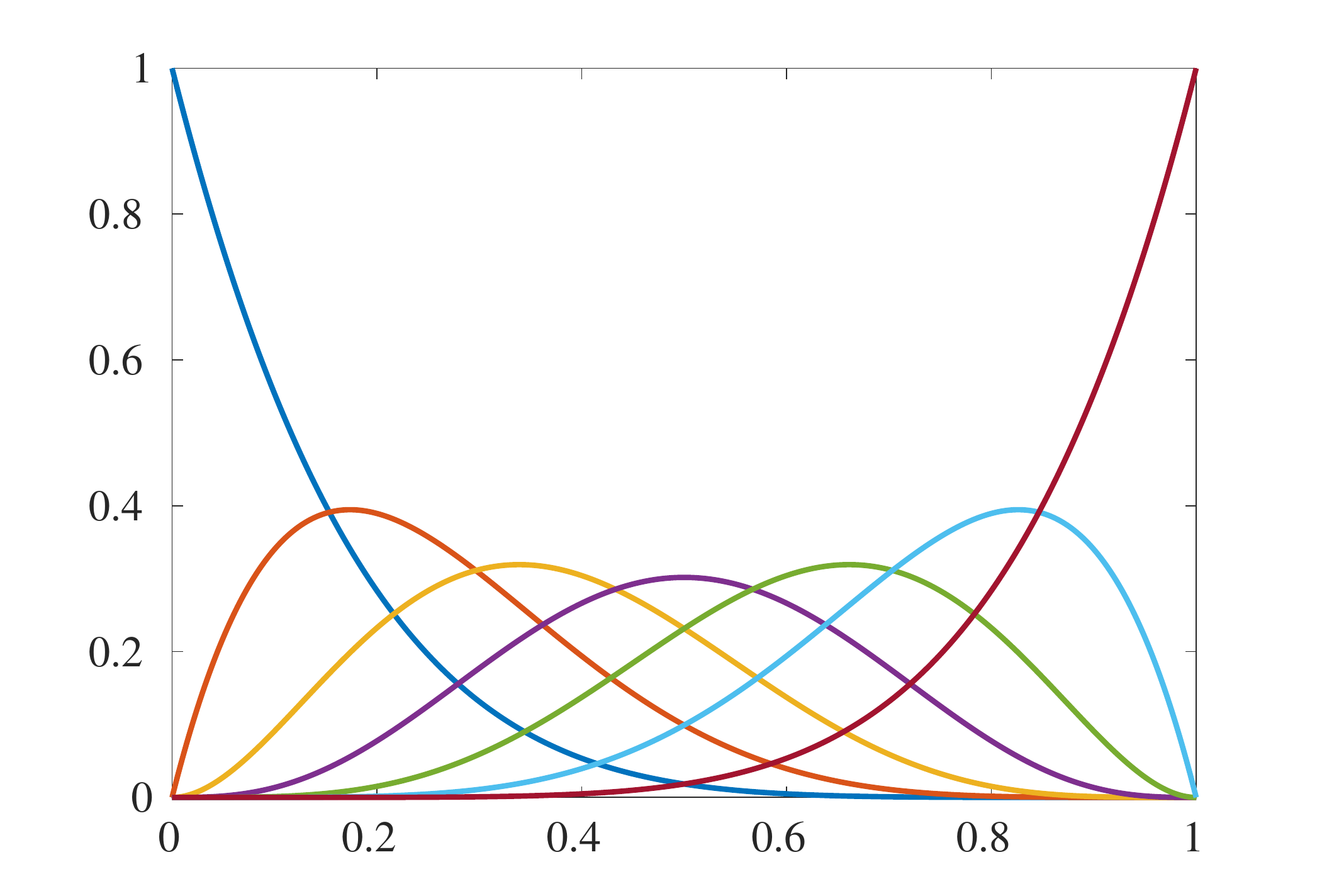}
\caption{Case study 2. Tchebycheffian Bernstein basis of $\spaceT_6^{(\myi \beta,-\myi \beta, \myi 2\beta, -\myi 2 \beta,\myi 3\beta,-\myi 3\beta)}$, $\beta=\frac{\pi}{2}$, on the interval $[0,1]$.}
\label{fig_ex2_benstn} 
\end{figure}

Here we focus again on a Poisson problem by considering \eqref{eq_adv_diff} with $\kappa=1$, $\textbf{a}=\textbf{0}$, and ${\rm f}$ obtained from the exact solution $$u(x,y) = 4((x+4)^2+y^2) + x(y-4),$$  to exhibit the evolution of error under uniform refinement and verify the consistency of the TB-splines for optimal convergence. The geometry map from the parametric domain $\widehat\Omega= (0, 1)^2$ to the physical domain $\Omega$ is given by
\begin{equation} \label{eq_map_mickey}
\begin{pmatrix} x\\ y \end{pmatrix} = \textbf{G}
\begin{pmatrix} s\\ t \end{pmatrix} = (1 - s)r \, \begin{pmatrix} \cos(\beta t)\\ \sin(\beta t) \end{pmatrix}
 +\  sR\, 
\begin{pmatrix} C + \cos(2\beta t-\gamma)\\ C + \sin(2\beta t-\gamma) \end{pmatrix},
\end{equation}
with
 $$r = 1,\quad R = 2,\quad C = 1/\sqrt{2},\quad \gamma = \frac{\pi}{4},\quad \beta = \frac{\pi}{2},$$
which leads to the domain depicted in \Cref{fig_ex_mickey}, left. Similar to the previous case study, the geometry consists of two arcs with phases $\beta$ and $2\beta$, hence a suitable selection of the ECT-space in the $t$-direction is
\begin{equation}\label{eq_space_mickey_t}
\spaceT_{p_2}^{(\myi \beta,-\myi \beta,\dots,\myi q\beta,-\myi q\beta)} = \angleSpace{ 1,\cos(\beta t), \sin(\beta t),\dots, \cos(q\,\beta t), \sin(q\, \beta t) },\quad p_2=2q\geq 4.
\end{equation}
The Tchebycheffian Bernstein basis of this ECT-space is illustrated in \Cref{fig_ex2_benstn} for $p_2 =6$.
From \Cref{ex_ppoly} we know that the corresponding Tchebycheffian spline space admits a TB-spline basis on any partition of the interval $[0, 1]$ since $\beta<\pi$. 
In the $s$-direction, we simply consider algebraic polynomials of the same degree as in the $t$-direction.
This gives the tensor-product space 
$$\spaceT_{p_1} \otimes \spaceT_{p_2}^{(\myi\frac{\pi}{2},-\myi\frac{\pi}{2},\ldots, \myi q\frac{\pi}{2}, -\myi q\frac{\pi}{2})},\quad p_1=p_2=p=2q. $$ 
Again, this geometry cannot be exactly represented by using generalized polynomial B-splines (see \cite{manni2011generalized}) with pieces in ECT-spaces of the form \eqref{eq_gb_space_trig}.

\begin{figure}[t!]
\captionsetup[subfigure]{aboveskip=-1pt,belowskip=0cm}
\centering
\includegraphics[width=.5\textwidth]{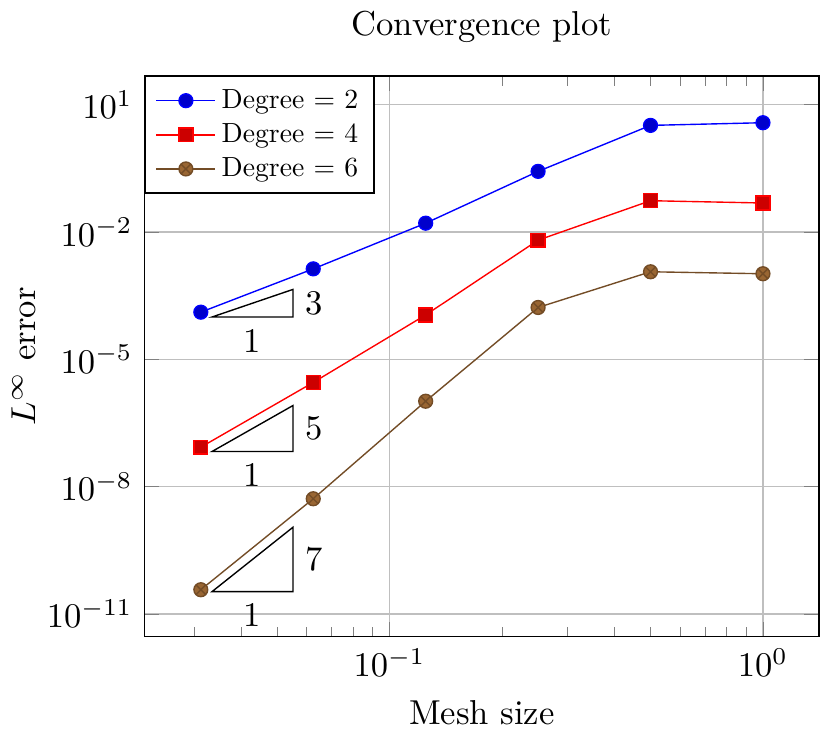}
  \caption{Case study 2. Convergence plot of the $L^\infty$ error of tensor-product TB-splines identified by $\xheightsub{\spaceT}{p}\otimes\spaceT_p^{(\myi \beta,-\myi \beta, \myi 2\beta, -\myi 2\beta,\dots,\myi{q\beta},-\myi{q\beta})} $, $\beta=\frac{\pi}{2}$,  for $p=2q$, $q=\{1,2,3\}$, and different resolution of the mesh $h=1,\frac{1}{2},\frac{1}{2^2},\dots,\frac{1}{2^5}$.}
  \label{fig_ex_mickey_conv} 
\end{figure}

\Cref{fig_ex_mickey_conv} illustrates the convergence of the error in $L^\infty$ norm,
which has been approximated by sampling the approximate and exact solutions on a uniform grid in the parametric domain consisting of $501$ points along each direction.
The plot clearly shows that the approximate solution obtained by TB-splines exhibits the optimal convergence of order $p+1$ as classical polynomial B-splines in agreement with the theoretical expectations; see \Cref{thm_approx}. This establishes the reliability of the Tchebycheffian spline spaces from the accuracy point of view. 
The contour plot of the approximate solution with $p=6$ on a mesh consisting of $8\times 8$ elements is presented in \Cref{fig_ex_mickey}, right.
 
%
\subsection{Case study 3: Advection-diffusion problem with radial advection flow}
\label{sec_Example3}

In the experiments presented so far, the selection of the ECT-spaces has been driven by their geometrical features. Now we will also take into account some analytical features of the problem. Let us consider the advection-diffusion problem \eqref{eq_adv_diff} with 
\begin{equation}
\label{eq_ex_adv_radial_data}
 \kappa = 1,\quad \mathbf{a}\,(x, y) \coloneqq\mathfrak{a} \begin{pmatrix}
\frac{x}{\sqrt{x^2+y^2}}\\
\frac{y}{\sqrt{x^2+y^2}},\end{pmatrix},
\quad \mathfrak{a}=100,\quad {\rm f}(x,y)=1, 
\end{equation}
and homogeneous boundary conditions, on the same geometry used in the previous case study given by the geometry map \eqref{eq_map_mickey}.

It is well known that for problems with high \peclet number $(\bm{Pe}_g\gg1)$ approximate solutions belonging to piecewise polynomial spaces tend to exhibit spurious oscillations unless the discretization is fine enough. 
The above mentioned issue can be efficiently addressed by turning to Tchebycheffian splines built from exponential functions with suitable parameters. Their use for solving advection-dominated advection-diffusion problems reduces, and sometimes even eliminates, the spurious oscillations while still providing an accurate localization of the sharp layer(s) without demanding extremely fine meshes and so keeping the number of degrees of freedom (dof) reasonably low.

\begin{figure}[t!] \centering
	\includegraphics[width=.4\linewidth]{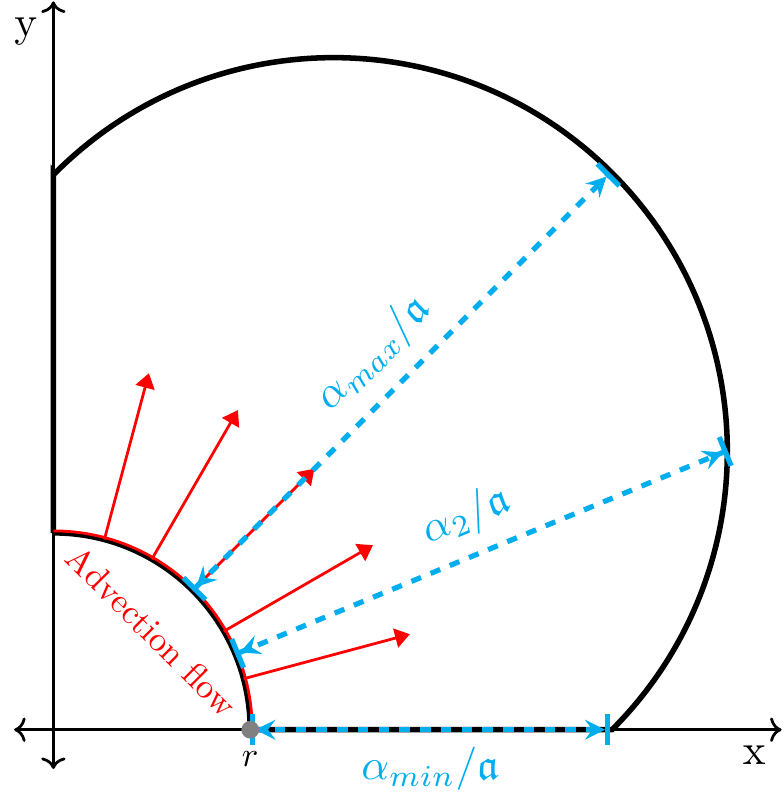}
	\caption{Case study 3. Illustration of the selection strategy of the shape parameters $\{\alpha_i: i=1,\ldots,\ell\}$ of the exponential functions to accommodate for different distances traveled by the advection flow for $\ell=3$ in \eqref{eq_ex_exp_adv_radial}. The extremes of the chosen shape parameters $\alpha_1 = \alpha_{min} = \mathfrak{a}\,(\sqrt{2}R - r)$ and $\alpha_3=\alpha_{max} = \mathfrak{a}\,(2R - r)$ are evidently related to the maximum and minimum radial distances of the domain, where $\mathfrak{a}$ is the modulus of the advection coefficient.}
	\label{fig_ex_adv_radial_description}
\end{figure}  

\begin{figure}[t!] \centering
	\includegraphics[width=0.75\linewidth]{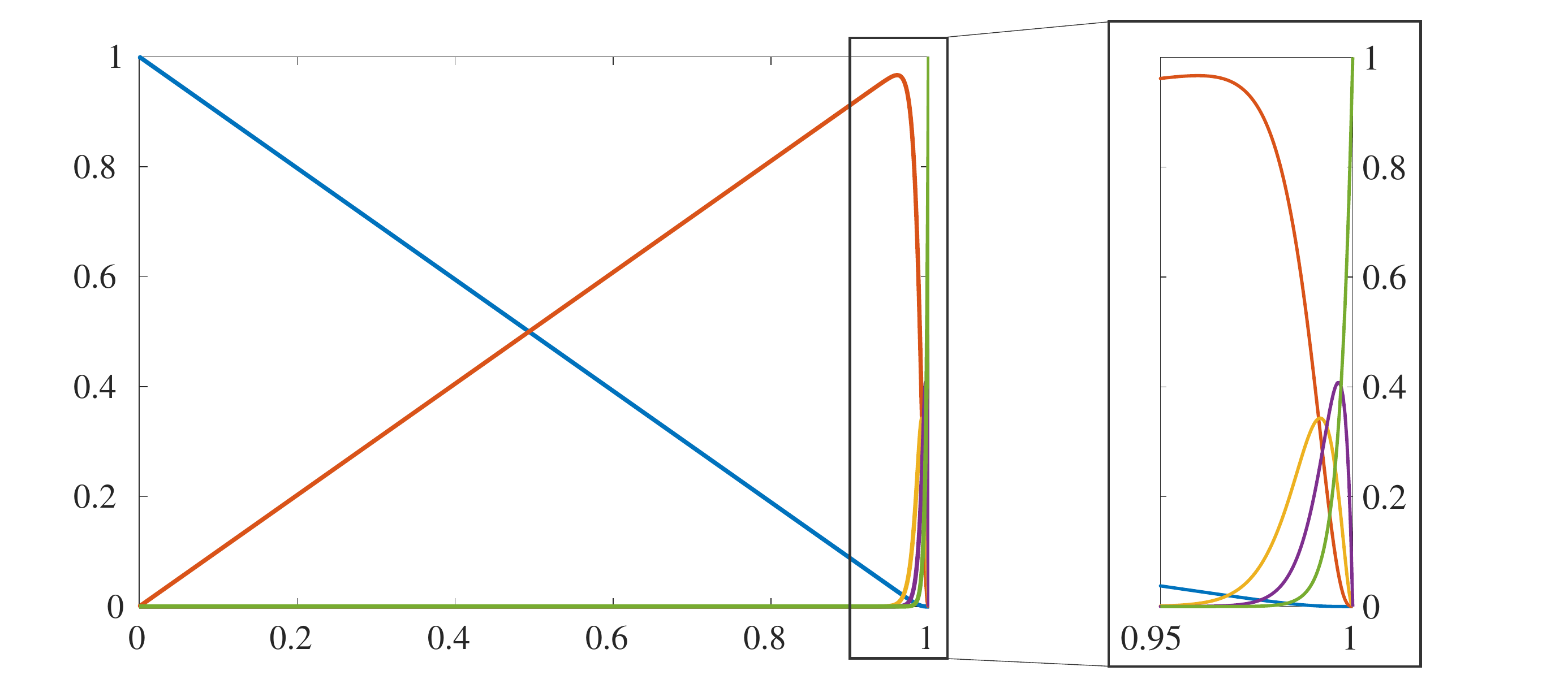}
	\caption{Case study 3. Tchebycheffian Bernstein basis of $\spaceT_4^{(\alpha_1,\alpha_2,\alpha_3)}$ on the interval $[0,1]$, with a magnified view of the functions on the interval $[0.95,1]$.}
	\label{fig_ex_adv_radial_brnstn}
\end{figure}

For an advection-dominated system with advection flow in the radial direction as in \eqref{eq_ex_adv_radial_data}, we can exploit exponential functions in the parametric $s$-direction with suitable shape parameters to represent the sharp layer generated near the homogeneous boundary.    
To construct effective Tchebycheffian splines with exponential functions in the $s$-direction, we consider the ECT-space
\begin{equation*}
\spaceT_{p_1}^{(\alpha_1,\alpha_2,\ldots,\alpha_\ell)} = \angleSpace{ 1, s,\ldots, s^{p_1-\ell}, e^{\alpha_1 s},\dots e^{\alpha_\ell s}},\quad \quad p_1\geq\ell\geq1,
\end{equation*}
where for $i = 1,\ldots,\ell,$
\begin{equation}\label{eq_ex_exp_adv_radial}
\alpha_i \coloneqq
\begin{cases} 
\frac{1}{2}(\alpha_{max}+\alpha_{min}), & \ell=1,
\\[0.1cm]
\alpha_{min} + \frac{i-1}{\ell-1}(\alpha_{max}-\alpha_{min}),& \ell\geq2,
\end{cases}
\end{equation}
and
$$ \alpha_{min}\coloneqq \mathfrak{a}\,(\sqrt{2}R - r)\,,\quad \alpha_{max}\coloneqq \mathfrak{a}\,(2R - r). $$
The governing factors in the selection of the shape parameters here are the global \peclet number \eqref{eq_peclet}, which is $\bm{Pe}_g = \mathfrak{a}$, and the range of distances traveled by the advection flow in the radial direction. \Cref{fig_ex_adv_radial_description} illustrates this approach for the selection of the shape parameters by considering $\ell=3$ different radial distances traveled by the advection flow. 
The Tchebycheffian Bernstein basis of this ECT-space is illustrated in \Cref{fig_ex_adv_radial_brnstn} for $p_1 =4$.
From \Cref{ex_TB_exp} we know that the corresponding Tchebycheffian spline space admits a TB-spline basis on any partition of the interval $[0, 1]$.

As discussed in \Cref{sec_Example2}, to exactly reproduce the geometry, we consider in the $t$-direction the ECT-space \eqref{eq_space_mickey_t} of degree $p_2=2q\geq 4$,
which allows for a Tchebycheffian spline space equipped with a TB-spline basis on the interval $[0, 1]$. 
The resulting tensor-product space is
$$\spaceT_{p_1}^{(\alpha_1,\alpha_2,\ldots,\alpha_{\ell})} \otimes \spaceT_{p_2}^{(\myi\frac{\pi}{2},-\myi\frac{\pi}{2},\ldots, \myi q\frac{\pi}{2}, -\myi q\frac{\pi}{2})},\quad p_2=2q,$$ 
with both spaces of the form \eqref{eq_spaces}.

\begin{figure}[t!]
	\begin{subfigure}{0.5\textwidth}
		\includegraphics[width=\linewidth]{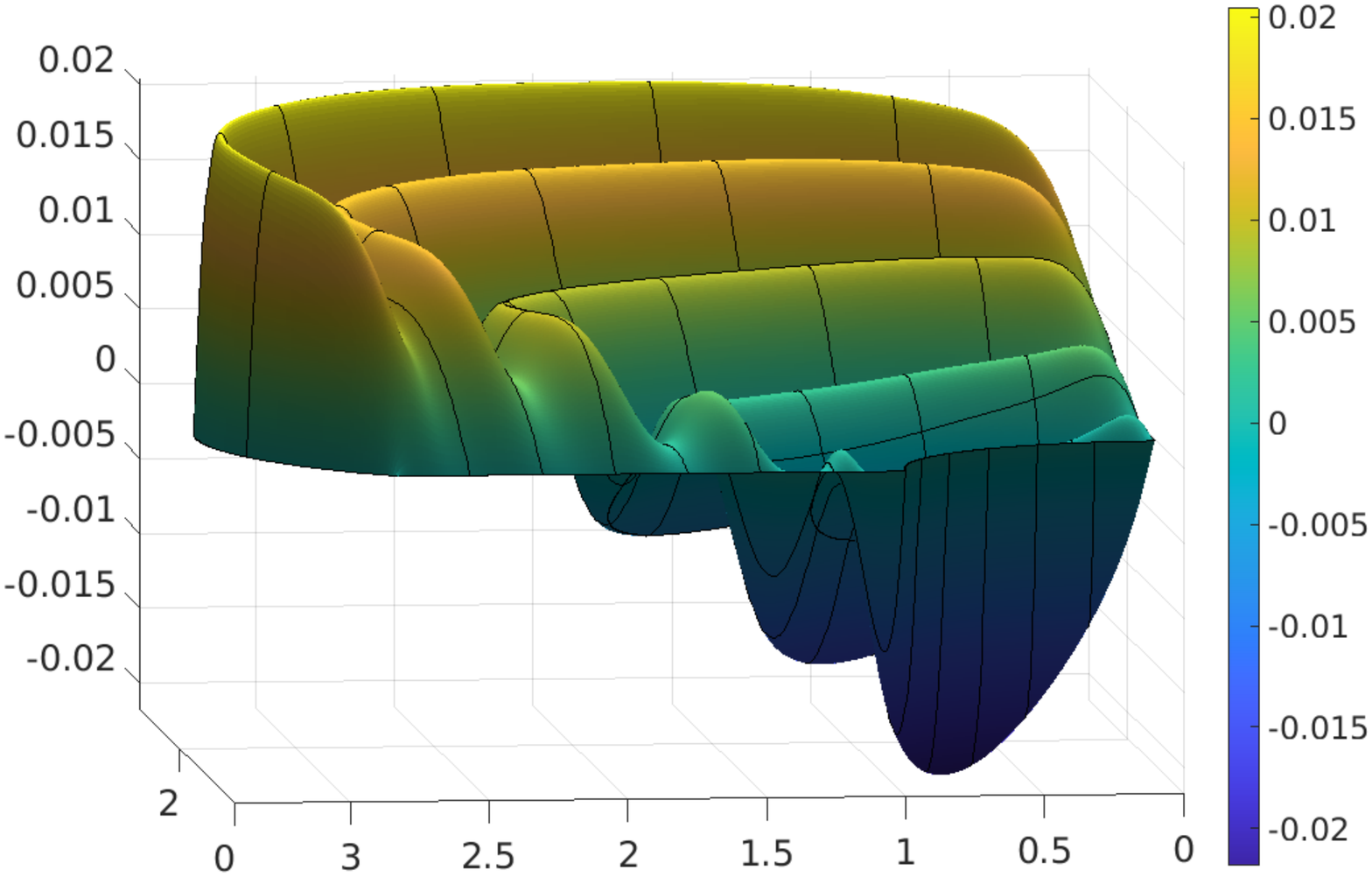}
		\caption{TB-spline space: $\xheightsub{\spaceT}{4} \otimes \spaceT_4^{(\myi \frac{\pi}{2},-\myi \frac{\pi}{2}, \myi \pi, -\myi \pi)}$} \label{fig_ex_adv_rad_exp0}
	\end{subfigure}\hspace*{\fill}
	\begin{subfigure}{0.5\textwidth}
		\includegraphics[width=\linewidth]{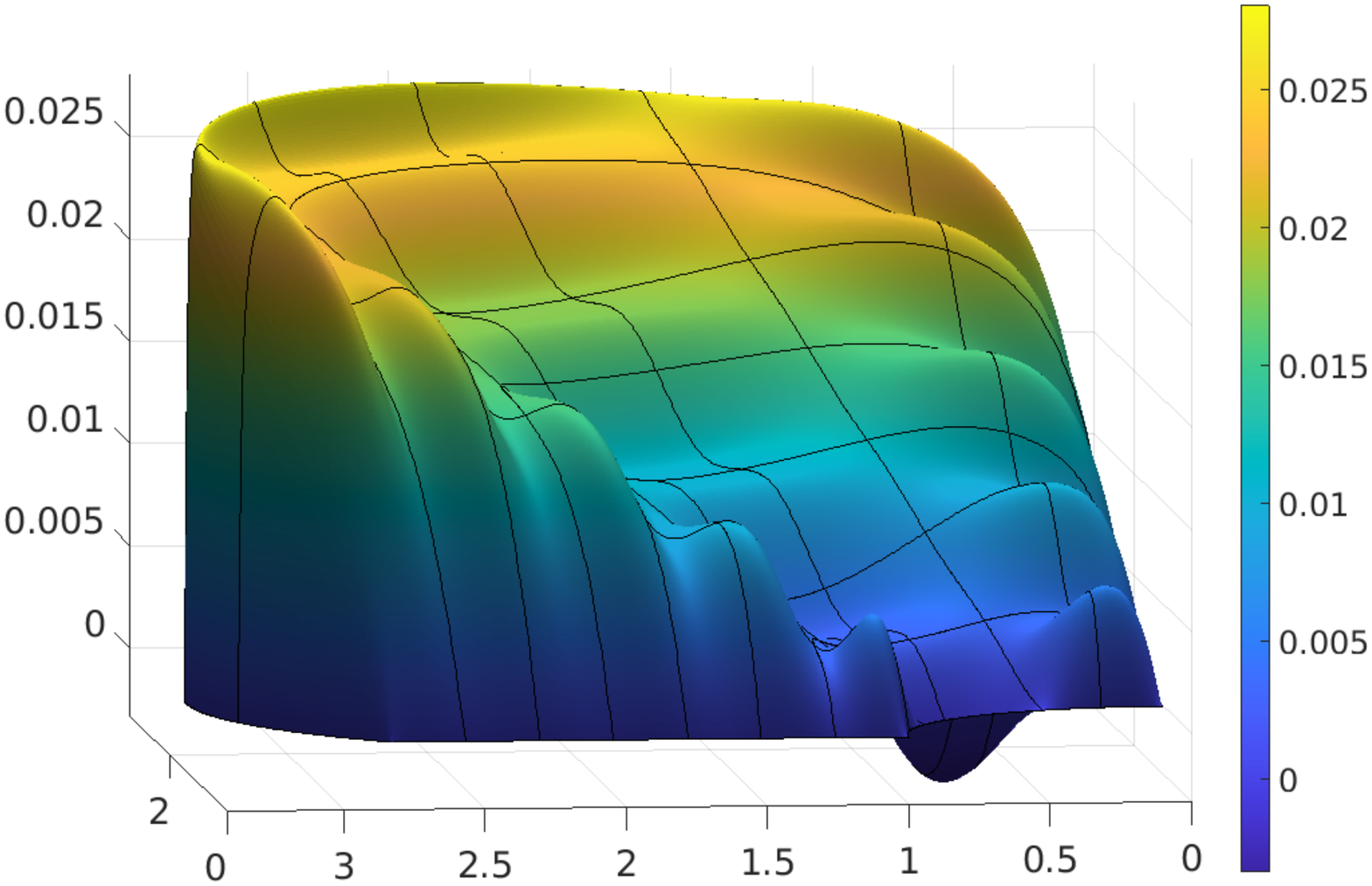}
		\caption{TB-spline space: $\spaceT_{4}^{(\alpha_1)} \otimes \spaceT_4^{(\myi \frac{\pi}{2},-\myi \frac{\pi}{2}, \myi \pi, -\myi \pi)}$} \label{fig_ex_adv_rad_exp1}
	\end{subfigure}
	\medskip
	\begin{subfigure}{0.5\textwidth}
		\includegraphics[width=\linewidth]{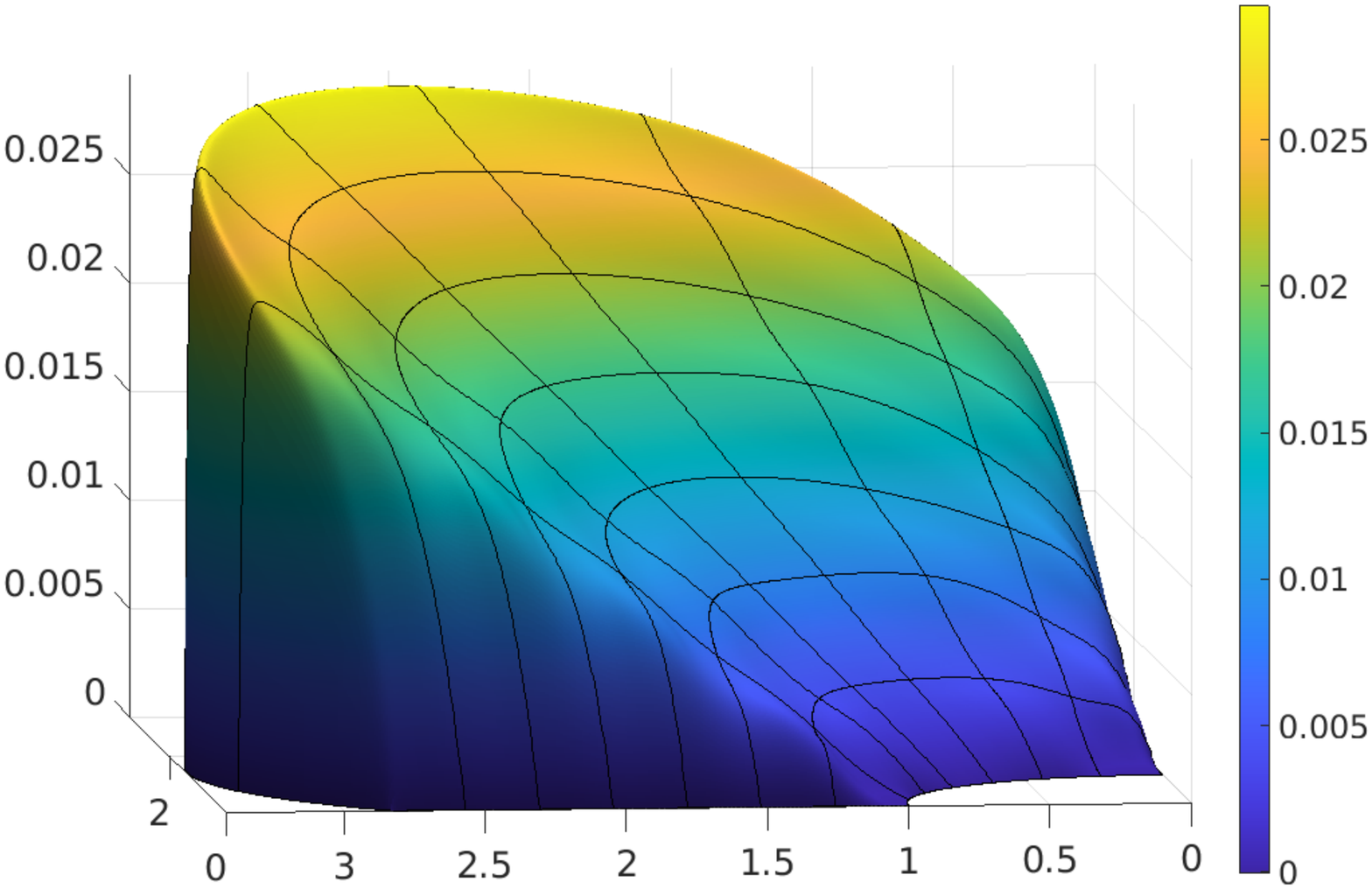}
		\caption{TB-spline space: $\spaceT_{4}^{(\alpha_1, \alpha_2)} \otimes \spaceT_4^{(\myi \frac{\pi}{2},-\myi \frac{\pi}{2}, \myi \pi, -\myi \pi)}$} \label{fig_ex_adv_rad_exp2}
	\end{subfigure}\hspace*{\fill}
	\begin{subfigure}{0.5\textwidth}
		\includegraphics[width=\linewidth]{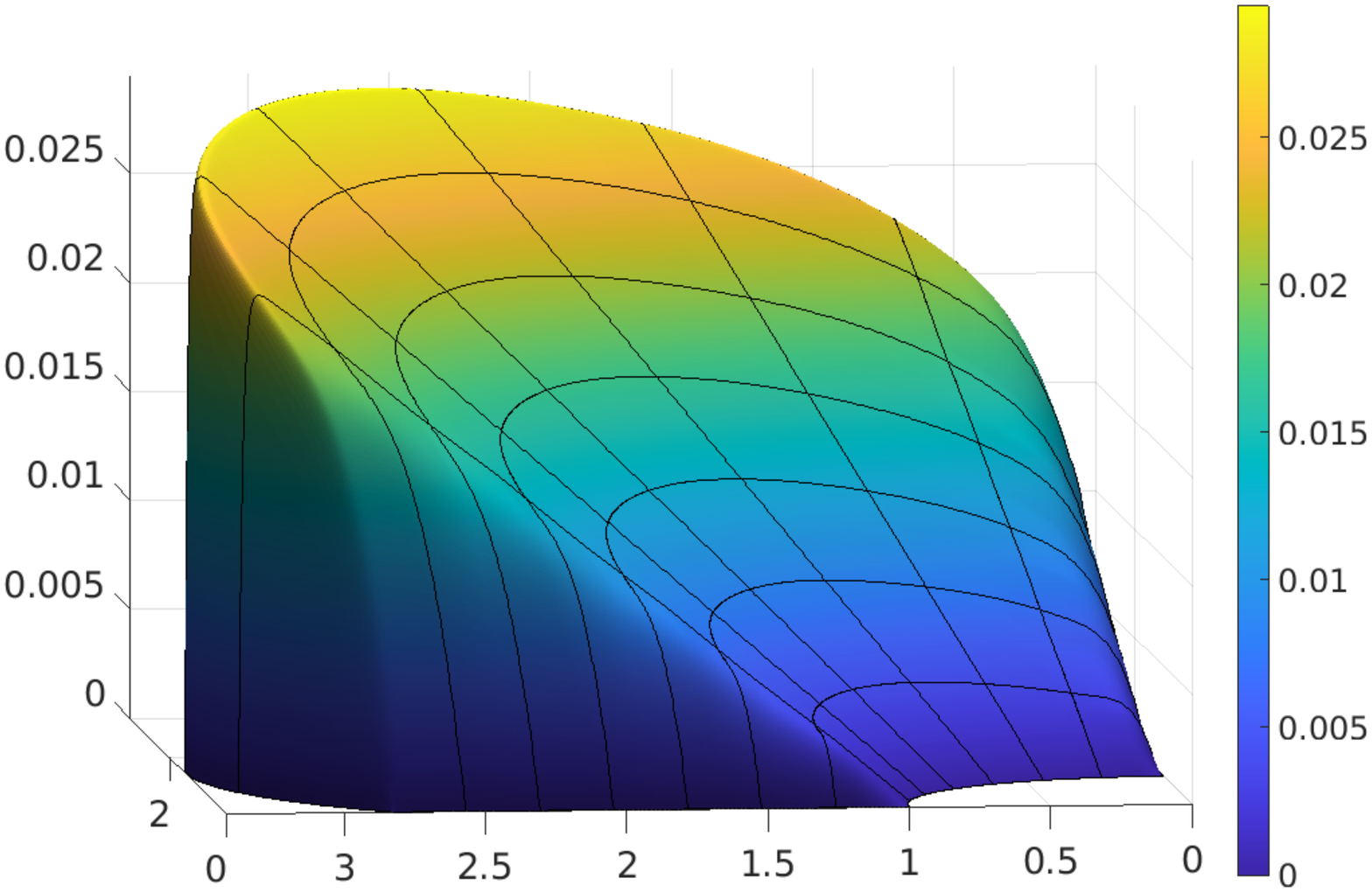}
		\caption{TB-spline space: $\spaceT_{4}^{(\alpha_1, \alpha_2, \alpha_3)} \otimes \spaceT_4^{(\myi \frac{\pi}{2},-\myi \frac{\pi}{2}, \myi \pi, -\myi \pi)}$} \label{fig_ex_adv_rad_exp3}
	\end{subfigure}
	\caption{Case study 3. Plots of the approximate solution to the advection-diffusion problem with global \peclet number $\bm{Pe}_g= 100$, obtained by using tensor-product TB-splines  with $p_1= p_2=4$, $m= 8$, dof = $144$. Increasing the number of exponential functions with shape parameters $\{\alpha_i: i=1,\ldots,\ell\}$ according to \eqref{eq_ex_exp_adv_radial} in the TB-spline space results into a reduction of oscillations in the solution.} \label{fig_ex_adv_radial}
\end{figure}

In \Cref{fig_ex_adv_radial} we present some results obtained by using tensor-product TB-splines of degrees $p_1=p_2=4$ on a mesh consisting of $8\times8$ elements. 
To exhibit the effect of the presence of exponential functions in the TB-spline space for addressing the advection-dominant problem \eqref{eq_adv_diff} with \eqref{eq_ex_adv_radial_data}, we vary the number of shape parameters in the selection from \eqref{eq_ex_exp_adv_radial}.
The outcome of using only algebraic polynomial splines (so without any exponential functions) in the $s$-direction is shown in \Cref{fig_ex_adv_rad_exp0}.
Even though the geometry is described exactly with the trigonometric functions present in the $t$-direction, the absence of the exponential functions results in spurious oscillations in the solution. However, taking one exponential function in the TB-spline space somewhat reduces the oscillations as shown in \Cref{fig_ex_adv_rad_exp1}, and taking two exponential functions in the space almost removes all the over- and undershoots as shown in \Cref{fig_ex_adv_rad_exp2}. Following the trend, with three exponential functions in the TB-spline space, the oscillations are completely eliminated from the solution (at least visually) as illustrated in \Cref{fig_ex_adv_rad_exp3}.

%
\subsection{Case study 4: Advection-diffusion problem with tangential advection flow}
\label{sec_Example4}
The ECT-spaces used in the previous case studies do not exploit the full functionality of mixing exponential and trigonometric functions with polynomials, as in \eqref{eq_spaces}. Hence, now we consider a problem with analytical and geometrical features in the same parametric direction. Let us consider the advection-diffusion problem \eqref{eq_adv_diff} with 
\begin{equation}\label{eq_ex_adv_tan_data}
    \kappa = 1,\quad \mathbf{a}\,(x, y) \coloneqq\mathfrak{a} \begin{pmatrix}
\frac{-y}{\sqrt{x^2+y^2}}\\
 \frac{x}{\sqrt{x^2+y^2}},\end{pmatrix},
 \quad \mathfrak{a}=100,\quad {\rm f}(x,y)=1, 
\end{equation}
 and homogeneous boundary conditions, on a classical quarter of an annulus given by the geometry map
\begin{equation}\label{eq_map_annulus}
\begin{pmatrix} x\\ y \end{pmatrix} = \textbf{G}
\begin{pmatrix} s\\ t \end{pmatrix} = ((1 - s)r  + sR) 
\begin{pmatrix} \cos(\frac{\pi}{2} t)\\ \sin(\frac{\pi}{2} t) \end{pmatrix},
\quad r = 1,\quad R = 2. \end{equation}
which leads to the domain depicted in \Cref{fig-IgA-Galerkin}, right.

\begin{figure}[t!] \centering
    \includegraphics[width=.4\linewidth]{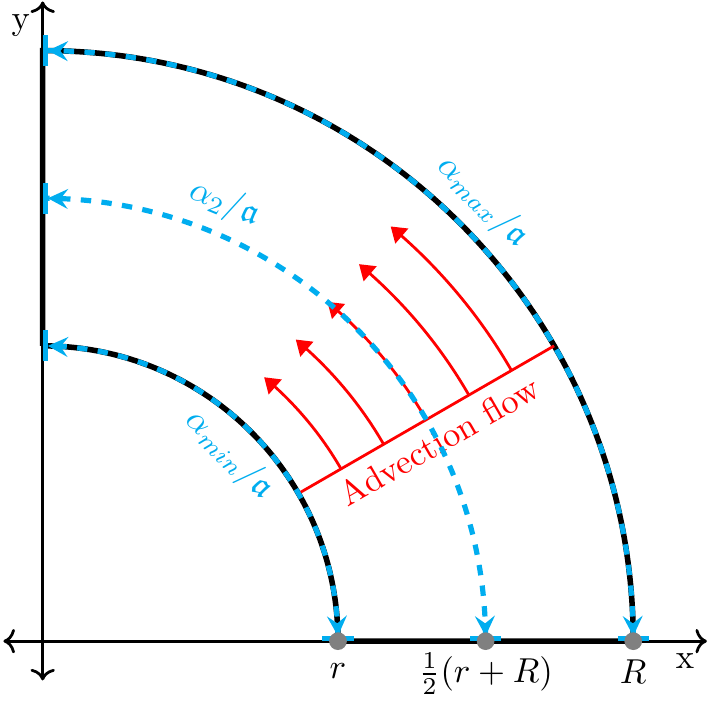}
    \caption{Case study 4. Illustration of the selection strategy of the shape parameters $\{\alpha_i: i=1,\ldots,\ell\}$ of the exponential functions to accommodate for different distances traveled by the advection flow for $\ell=3$ in \eqref{eq_ex_exp_annulus_tan}. The extremes of the chosen shape parameters $\alpha_1 = \alpha_{min}= \mathfrak{a}\,\frac{\pi}{2}r$ and $\alpha_3=\alpha_{max}= \mathfrak{a}\,\frac{\pi}{2}R$ are evidently related to the two extreme arc lengths of the domain, where $\mathfrak{a}$ is the modulus of the advection coefficient.}
    \label{fig_ex_adv_tan_description}
\end{figure} 

\begin{figure}[t!] \centering
\includegraphics[width=0.75\linewidth]{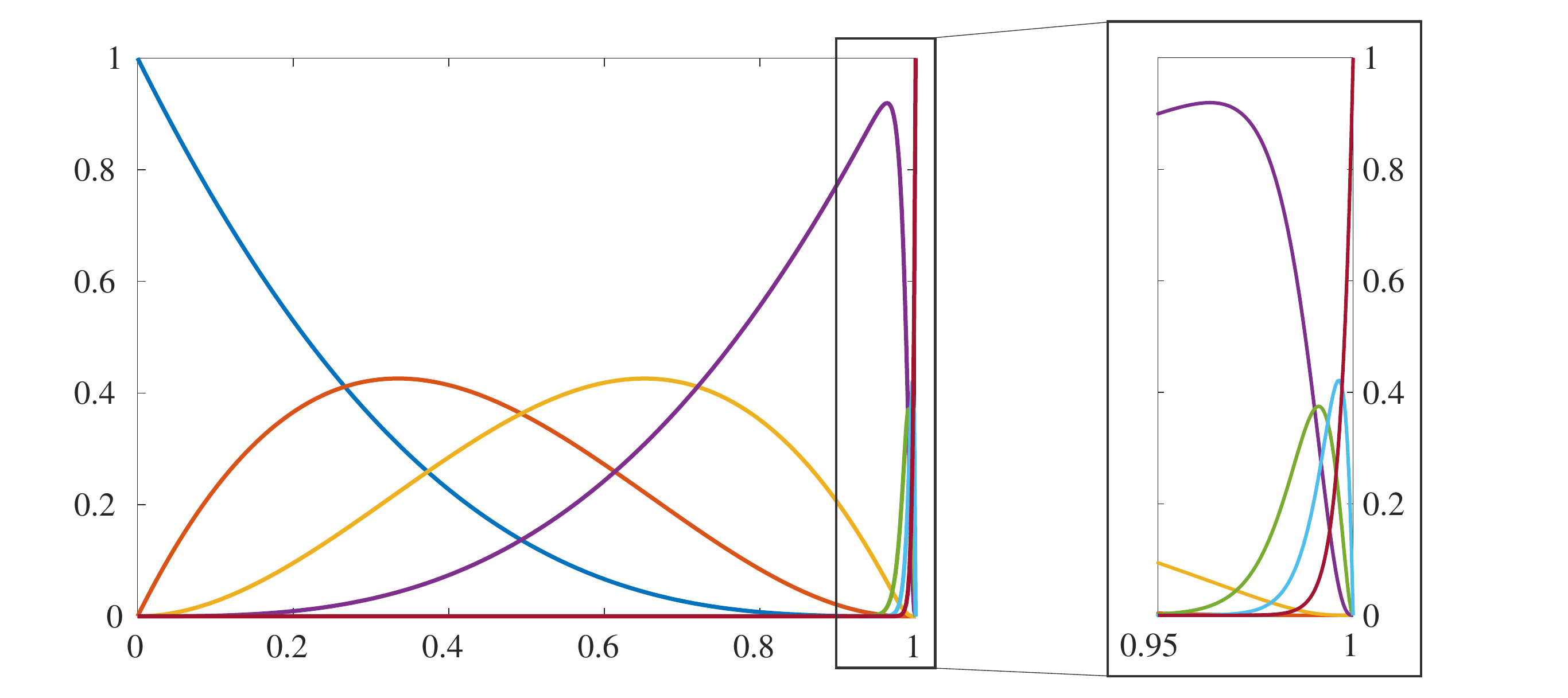}
    \caption{Case study 4. Tchebycheffian Bernstein basis of $\spaceT_6^{(\alpha_1,\alpha_2,\alpha_3,\myi \frac{\pi}{2},-\myi \frac{\pi}{2})}$ on the interval $[0,1]$, with a magnified view of the functions on the interval $[0.95,1]$.}
\label{fig_ex_adv_tabenstn} 
\end{figure}

To build suitable TB-spline spaces for the above problem, we need to consider both the geometry map in \eqref{eq_map_annulus} and the advection coefficient $\bm{a}$ in \eqref{eq_ex_adv_tan_data}. The geometry map requires trigonometric functions in the $t$-direction. Moreover, since we have a tangential flow along the $t$-direction, see \Cref{fig_ex_adv_tan_description}, it is natural to add exponential functions in that direction as well to accurately represent the sharp layer. The ECT-space in the $t$-direction can be chosen as
\begin{equation}\label{eq_tcheb_space_annulus}
\spaceT_{p_2}^{(\alpha_1,\alpha_2,\ldots,\alpha_\ell,\myi \frac{\pi}{2},-\myi \frac{\pi}{2})} = \Biggl\langle\, 1, t,\ldots, t^{p_2-2-\ell},e^{\alpha_1 t},\ldots e^{\alpha_\ell t}, \cos\biggl(\frac{\pi}{2} t\biggr), \sin\biggl(\frac{\pi}{2}t\biggr) \, \Biggr\rangle,\quad p_2\geq \ell+2\geq 3,
\end{equation}
where, similar to \eqref{eq_ex_exp_adv_radial}, for $i = 1,\dots,\ell$,
\begin{equation}\label{eq_ex_exp_annulus_tan} 
\alpha_i \coloneqq
\begin{cases} 
\frac{1}{2}(\alpha_{max}+\alpha_{min}), & \ell=1,
\\[0.1cm]
\alpha_{min} + \frac{i-1}{\ell-1}(\alpha_{max}-\alpha_{min}),& \ell\geq2,
\end{cases}
\end{equation}
and
$$ \alpha_{min}\coloneqq \mathfrak{a}\,\frac{\pi}{2}r,\quad \alpha_{max}\coloneqq \mathfrak{a}\,\frac{\pi}{2}R. $$
The governing factors in the selection of the shape parameters here are again the global \peclet number \eqref{eq_peclet}, which is $\bm{Pe}_g = \mathfrak{a}$, and the range of distances traveled by the advection flow along circular arcs. \Cref{fig_ex_adv_tan_description} illustrates this approach for the selection of the shape parameters by considering $\ell=3$ different distances along arcs traveled by the advection flow. 
The Tchebycheffian Bernstein basis of this ECT-space is illustrated in \Cref{fig_ex_adv_tabenstn} for $p_2=6$.
In the $s$-direction, we simply consider algebraic polynomials of degree $p_1=2$.
The resulting tensor-product space is 
$$\spaceT_{p_1} \otimes \spaceT_{p_2}^{(\alpha_1,\alpha_2,\dots,\alpha_{\ell},\myi \frac{\pi}{2},-\myi \frac{\pi}{2})}.$$ 

\begin{figure}[t!]
	\begin{subfigure}{0.5\textwidth}
		\includegraphics[width=\linewidth]{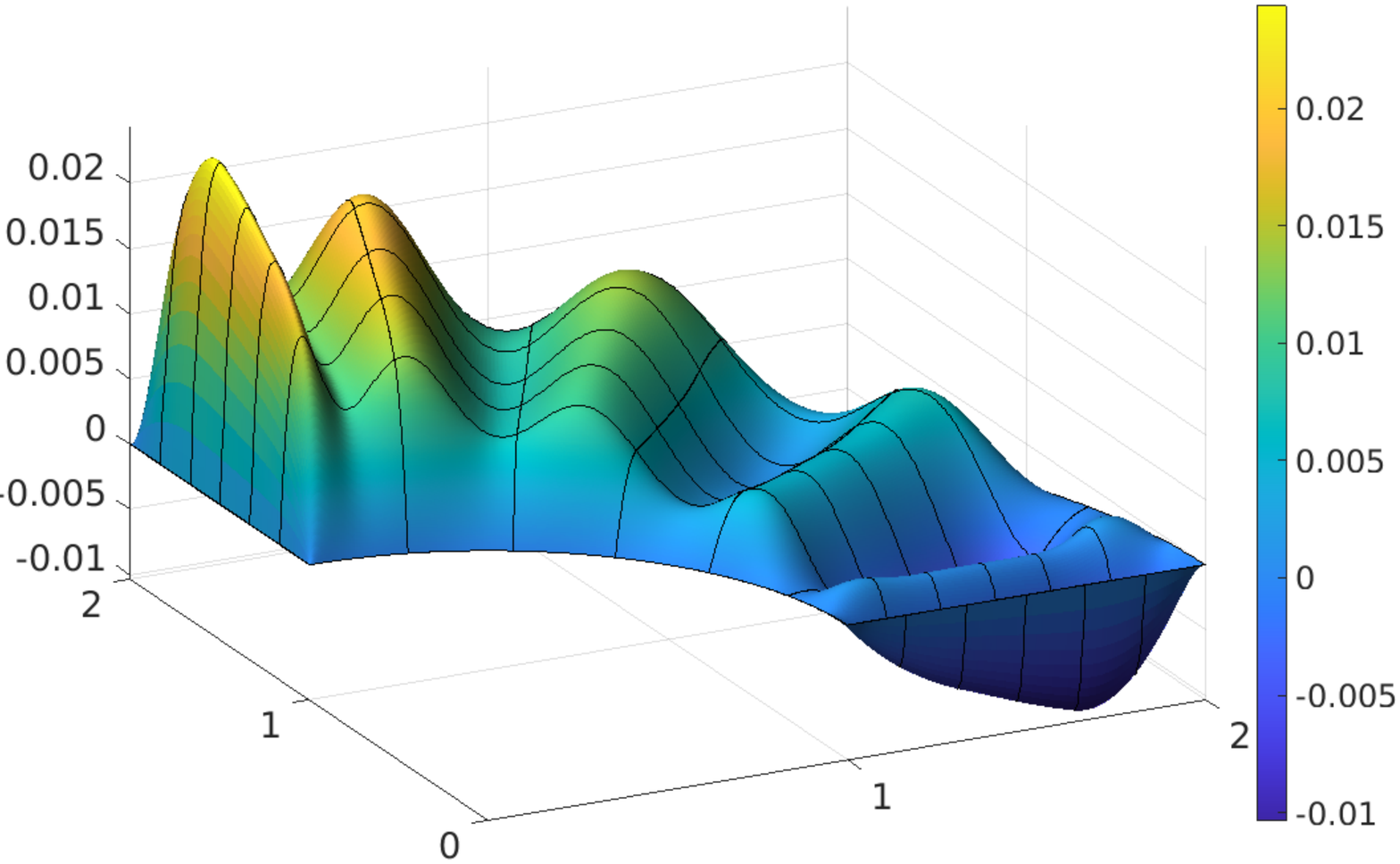}
		\caption{GB-spline space: $\xheightsub{\spaceT}{2} \otimes \spaceT_6^{(\myi \frac{\pi}{2},-\myi \frac{\pi}{2})}$} \label{fig_ex_annulus_gpoly}
	\end{subfigure}\hspace*{\fill}
	\begin{subfigure}{0.5\textwidth}
		\includegraphics[width=\linewidth]{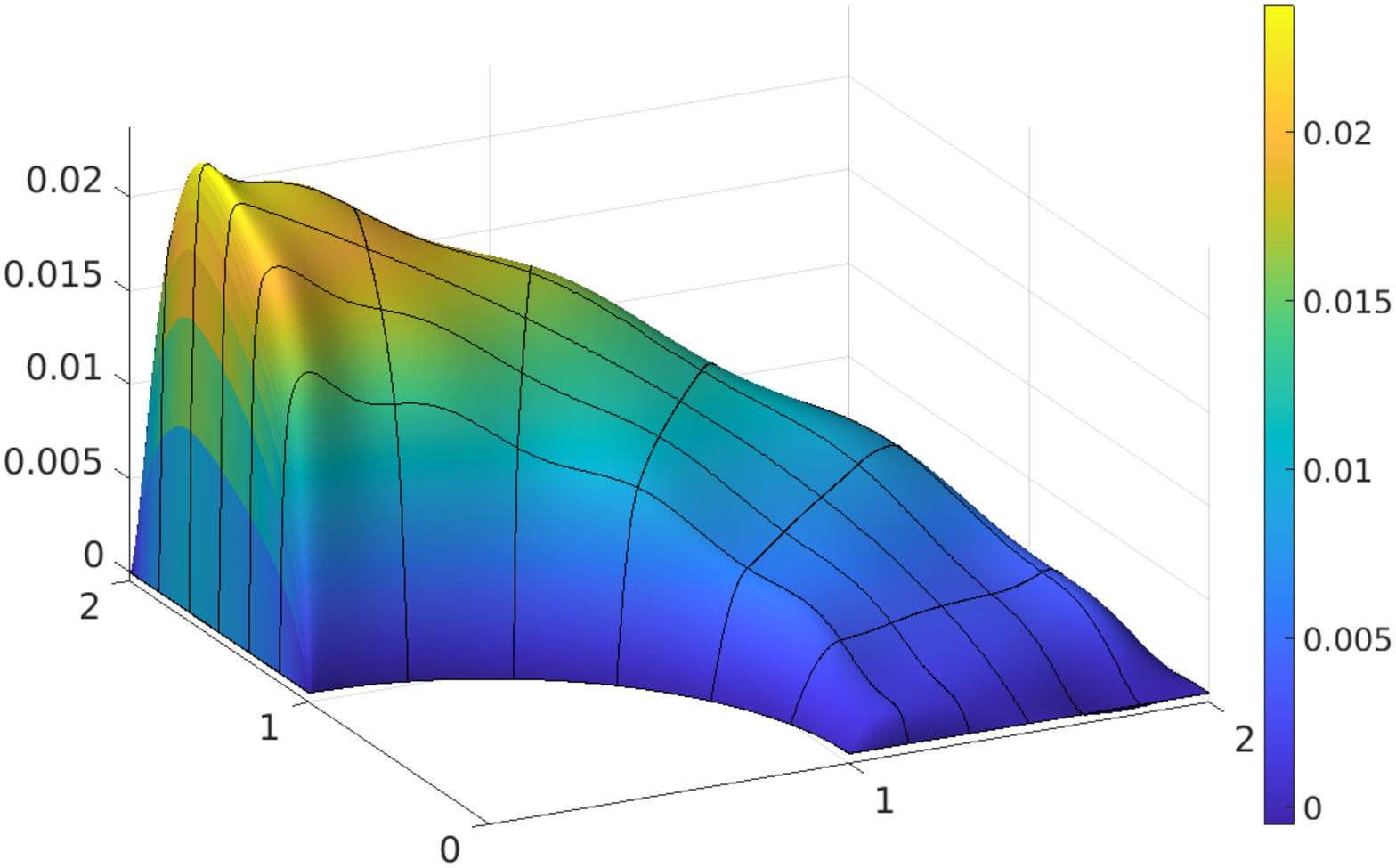}
		\caption{TB-spline space: $\xheightsub{\spaceT}{2} \otimes \spaceT_6^{(\alpha_1,\myi \frac{\pi}{2},-\myi \frac{\pi}{2})}$} \label{fig_ex_annulus_exp1}
	\end{subfigure}
	\medskip
	\begin{subfigure}{0.5\textwidth}
		\includegraphics[width=\linewidth]{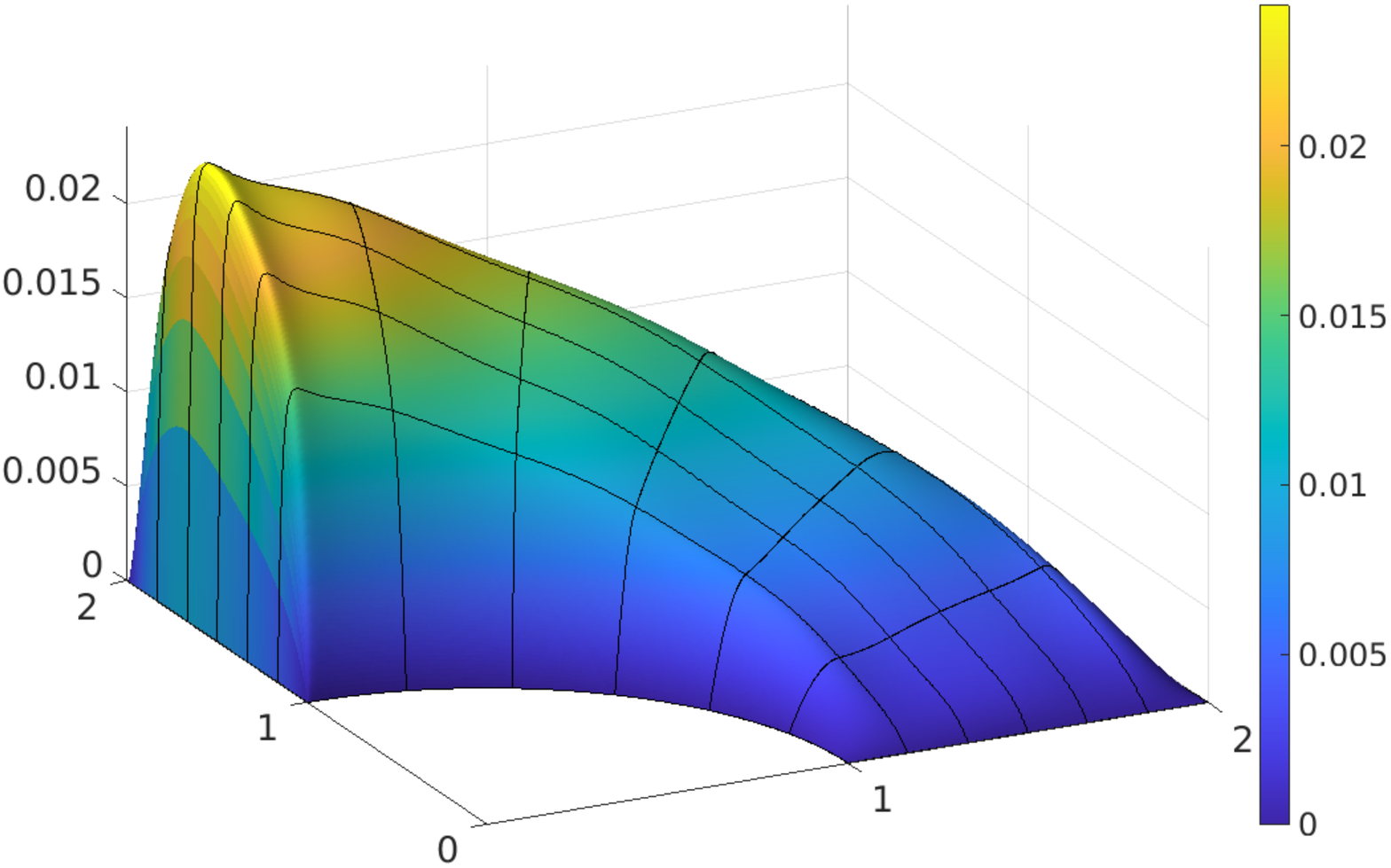}
		\caption{TB-spline space: $\xheightsub{\spaceT}{2} \otimes \spaceT_6^{(\alpha_1,\alpha_2,\myi \frac{\pi}{2},-\myi \frac{\pi}{2})}$} \label{fig_ex_annulus_exp2}
	\end{subfigure}\hspace*{\fill}
	\begin{subfigure}{0.5\textwidth}
		\includegraphics[width=\linewidth]{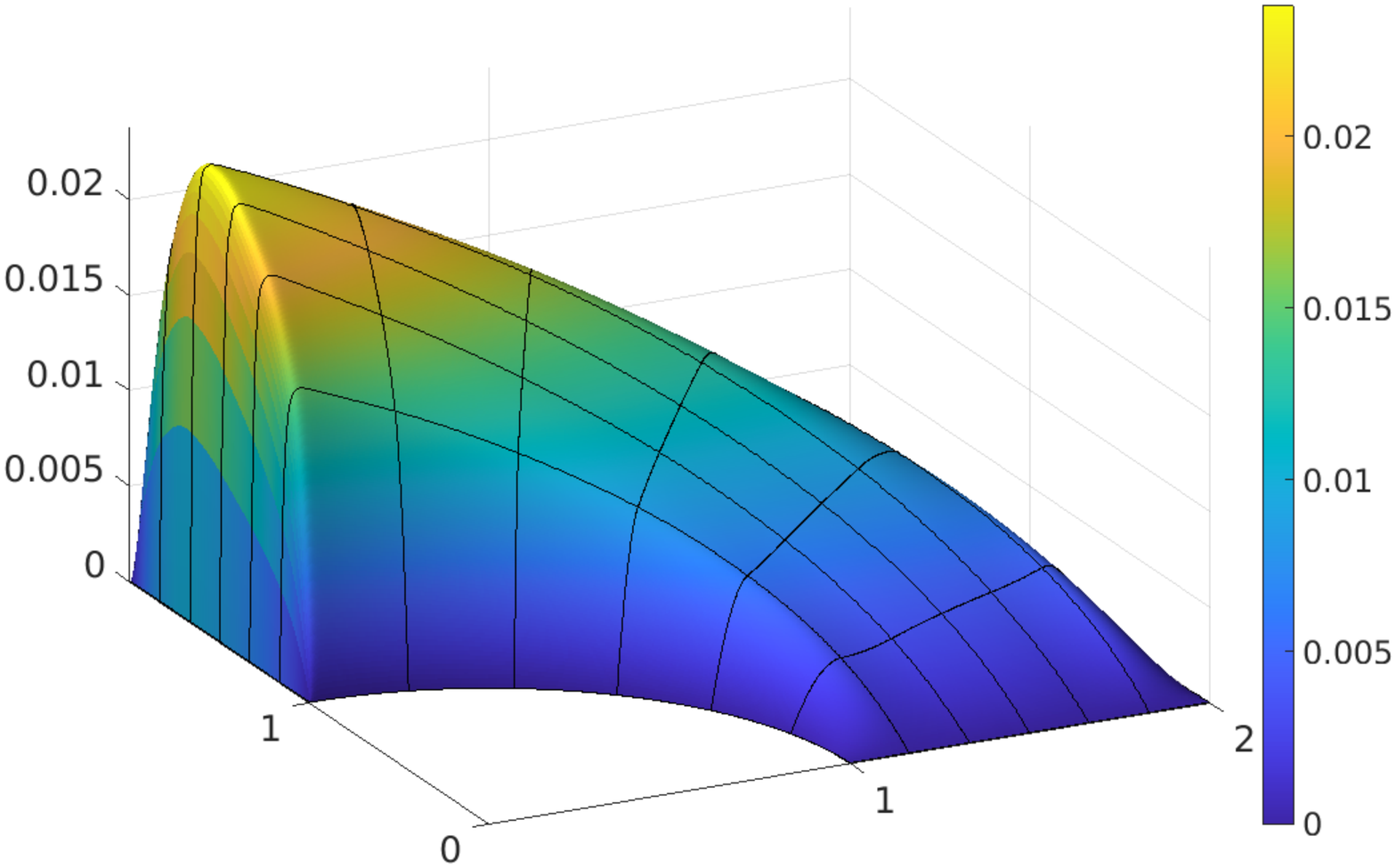}
		\caption{TB-spline space: $\xheightsub{\spaceT}{2} \otimes \spaceT_6^{(\alpha_1,\alpha_2,\alpha_3,\myi \frac{\pi}{2},-\myi \frac{\pi}{2})}$} \label{fig_ex_annulus_exp3}
	\end{subfigure}
	\caption{Case study 4. Plots of the approximate solution to the advection-diffusion problem with global \peclet number $\bm{Pe}_g = 100$, obtained by using tensor-product TB-splines with $p_1=2$, $p_2=6$, $ m = 6$, dof = $96$. Increasing the number of exponential functions with different shape parameters $\{\alpha_i: i=1,\ldots,\ell\}$ according to \eqref{eq_ex_exp_annulus_tan} in the TB-spline space results into a reduction of oscillations in the solution. } \label{fig_ex_annulus}
\end{figure}

For the sake of comparison, we also consider in the $t$-direction the generalized polynomial space obtained by adding to the polynomial space the pair of trigonometric functions necessary for an exact representation of the geometry, see \eqref{eq_gb_space_trig}, i.e.,
\begin{equation}\label{eq_gpoly_space_annulus}
\spaceT_{p_2}^{(\myi \frac{\pi}{2}, -\myi \frac{\pi}{2})} = \Biggl\langle\, 1, t,\dots, t^{{p_2}-2}, \cos\biggl(\frac{\pi}{2} t\biggr), \sin\biggl(\frac{\pi}{2}t\biggr)\, \Biggr\rangle,
\end{equation}
so that the resulting tensor-product space used for the comparison is
$$\spaceT_{p_1} \otimes \spaceT_{p_2}^{(\myi \frac{\pi}{2},-\myi \frac{\pi}{2})}.$$ 

Note that both the ECT-spaces in \eqref{eq_tcheb_space_annulus} and \eqref{eq_gpoly_space_annulus} are of the form \eqref{eq_spaces}. From \Cref{ex_TB_null-space} we know that the corresponding Tchebycheffian spline spaces admit a TB-spline basis on any partition of the interval $[0, 1]$.

In \Cref{fig_ex_annulus} we present some results obtained by using tensor-product TB-splines of degrees $p_1=2$ and $p_2=6$ on a mesh consisting of $6\times 6$ elements. 
To exhibit the effect of the presence of exponential functions in the TB-spline space for addressing the advection-dominant problem \eqref{eq_adv_diff} with \eqref{eq_ex_adv_tan_data}, we vary the number of shape parameters in the selection from \eqref{eq_ex_exp_annulus_tan}.
We can see in \Cref{fig_ex_annulus_gpoly} that the GB-spline space with only trigonometric functions leads to spurious oscillations in the solution.
However, taking just one exponential function in the TB-spline space already substantially reduces the oscillations 
as shown in \Cref{fig_ex_annulus_exp1}, and with two exponential functions all the over- and undershoots almost disappear as shown in \Cref{fig_ex_annulus_exp2}. The trend continues to follow and  increasing the number of exponential functions to $3$ ultimately gives us a solution with no oscillations (at least visually) as illustrated in \Cref{fig_ex_annulus_exp3}.

\subsection{Case study 5: Advection-diffusion problem with internal sharp layer on a square}
\label{sec_Example5}
In the previous two case studies the advection flow travels exactly along one of the parametric directions. Here we consider a classical benchmark problem with advection flow not parallel to any parametric direction \cite{hughes2005,manni2011tension,speleers2012}.

\begin{figure}[t!] \centering
\includegraphics[width=0.5\linewidth]{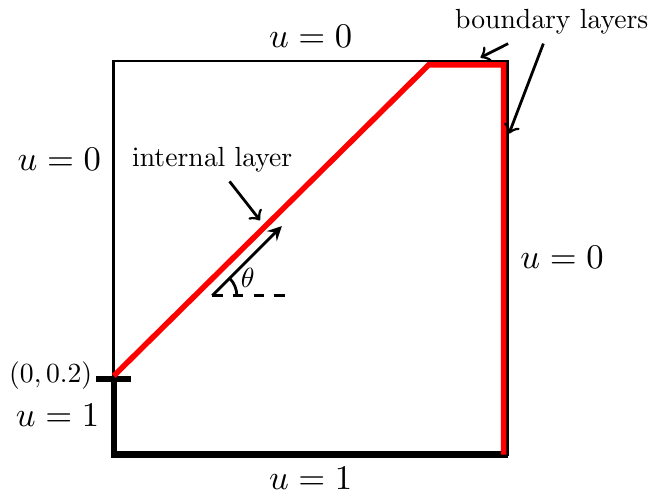}
\caption{Case study 5. The domain with the Dirichlet boundary conditions, the sharp internal layer (red) generated along the advection flow at angle $\theta$ and the boundary layers (red).}
\label{fig_ex_sqr_drchlt_bc}
\end{figure}

We want to solve the advection-diffusion problem \eqref{eq_adv_diff} on a square, hence the geometry map is just the identity, 
with
$$ \kappa = 1, \quad \mathbf{a} = \mathfrak{a}\, (\cos (\theta), \sin (\theta))^T,\quad \theta=\frac{\pi}{4},\quad \mathfrak{a}= 10^4, \quad {\rm f}(x,y) = 0,
$$
and the discontinuous Dirichlet boundary conditions as shown in \Cref{fig_ex_sqr_drchlt_bc}. 
The solution exhibits sharp boundary and internal layers depicted in red in the same figure.
In particular, the sudden jump in the boundary conditions at the point $(0, 0.2)$ generates an inner sharp layer aligned with the advection flow direction identified by
$(\cos (\theta), \sin (\theta))$.  
 
As already mentioned in \Cref{sec_Example3}, since the problem has a large global \peclet number, namely $\bm{Pe}_g = \mathfrak{a} = 10^4$, approximate solutions belonging to piecewise polynomial spaces exhibit spurious oscillations until the discretization is fine enough to resolve the sharp layers featured by the exact solution.
To overcome this issue, it is common to rely on stabilization methods such as the SUPG and GLS method; see \cite{brooks1982streamline,hughes1989new}. Typically, stabilization methods remove spurious oscillations but at the same time ``smooth out'' the layers featured by the exact solution. Their efficiency deeply depends on a careful choice of some parameters appearing in the various stabilization methods.

Tchebycheffian splines can offer a flexible alternative for the treatment of the above advection-dominated problem without the need for stabilization. To catch the sharp layers, it is natural to consider Tchebycheffian spline spaces that contain exponential functions with suitable parameters. Since the advection flow velocity is constant here, we can simply take its two components as parameters in the two (parameter) directions:
\begin{equation}\label{eq_space_square}
\spaceT_{p_1}^{(\mathfrak{a} \cos(\theta))} = \angleSpace{ 1, x,\dots, x^{p_1-1}, e^{\mathfrak{a}\cos(\theta) x}}, \quad
\spaceT_{p_2}^{(\mathfrak{a} \sin(\theta))} = \angleSpace{ 1, y,\dots, y^{p_2-1}, e^{\mathfrak{a}\sin(\theta) y}}.
\end{equation}
Both these spaces are of the form \eqref{eq_spaces} and, since they just contain polynomial and exponential functions, 
from \Cref{ex_TB_exp} we know that the corresponding Tchebycheffian spline spaces admit a TB-spline basis on any partition of the interval $[0, 1]$.
Note that, since the physical domain is simply a square, the parameter selection in \eqref{eq_space_square} is not governed by the geometry.

It is important to remark that a proper treatment of the boundary conditions is imperative due to its effect on the approximate solution. In order to avoid oscillations along the boundary, instead of using a least squares approximation as in the previous case studies, we have considered quasi-interpolation. More precisely, in order to exploit the shape-preserving property of the (T)B-spline representation, following the approach in \cite{hughes2005}, we approximate the boundary function through the Schoenberg operator (see \cite{lyche18,Schoenberg1967}) along each edge. In other words, we consider the linear combination of the boundary TB-splines whose coefficients are obtained by evaluating the boundary function at the corresponding Greville abscissae (under the assumptions that linear polynomials belong to the ECT-space of interest).

\begin{figure}[t!]
\vspace{-1cm}
  \begin{subfigure}{\textwidth}
  \centering
    \includegraphics[width=.45\linewidth]{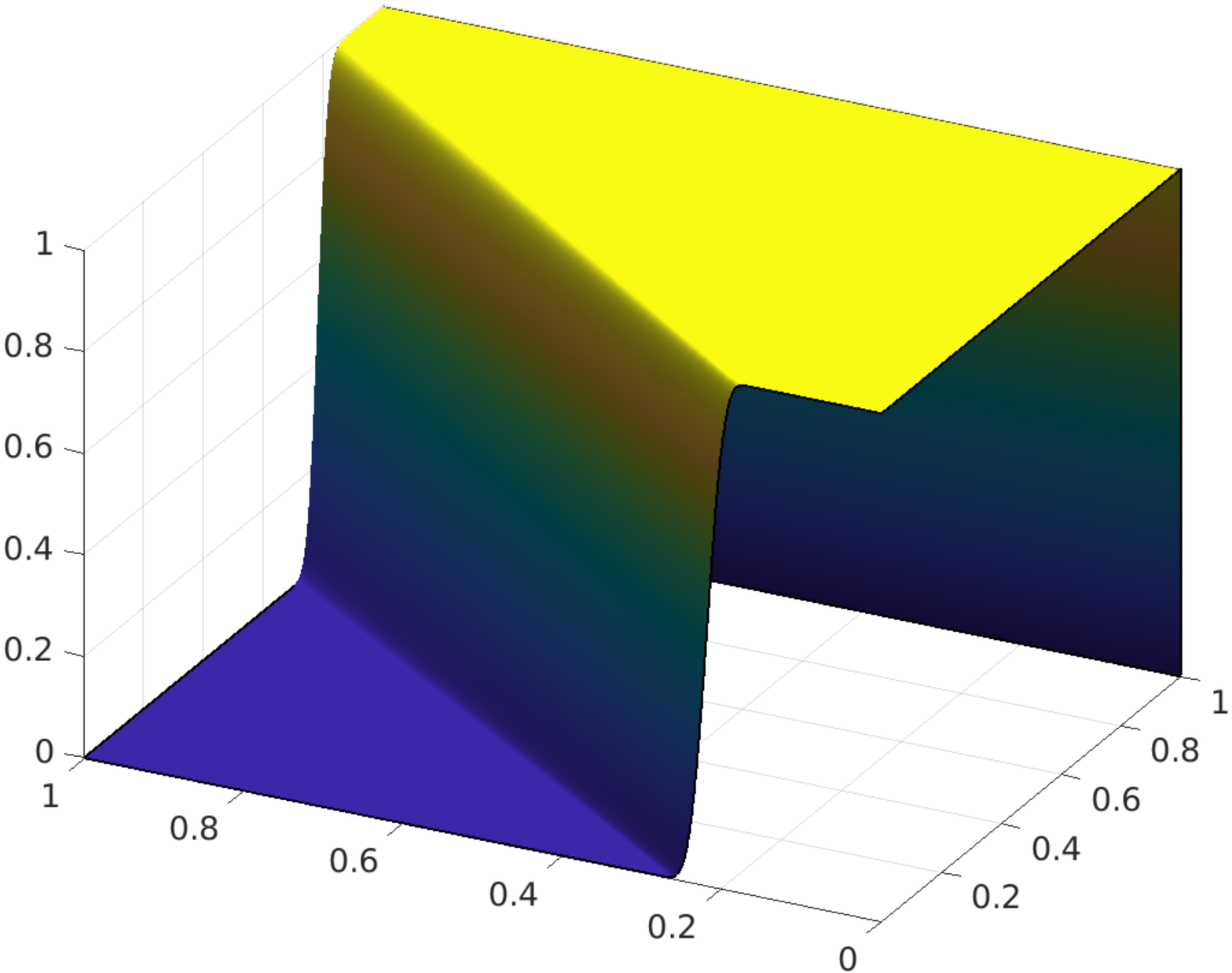}
    \hfill
    \includegraphics[width=.4\linewidth]{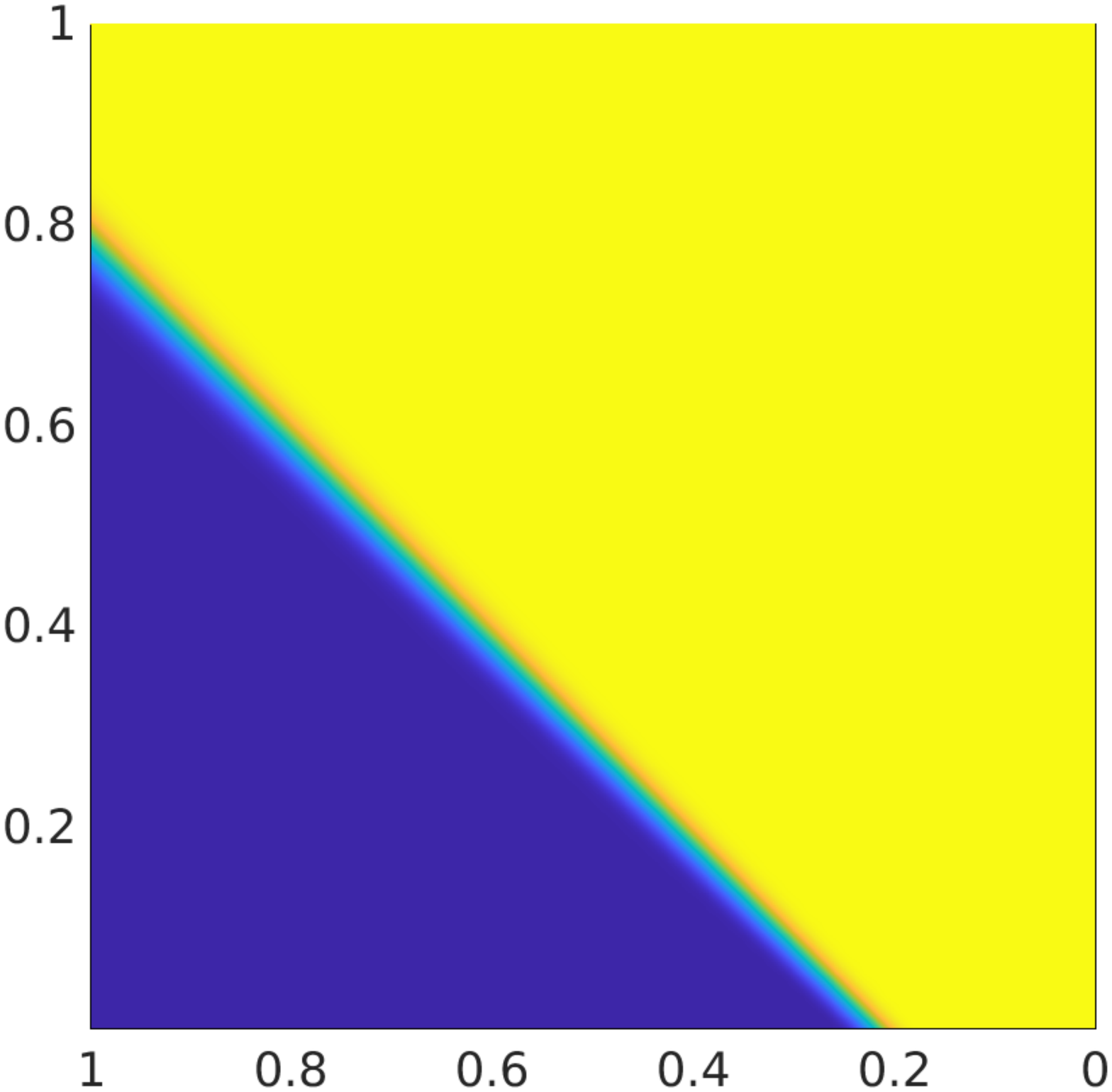}
    \caption{TB-spline space: $\spaceT_6^{(\mathfrak{a} \cos(\theta) )} \otimes \spaceT_6^{(\mathfrak{a} \sin(\theta) )}$ with $m=50$, dof $= 3136$.} \label{fig_adv_sqr_exp}
  \end{subfigure}
  \begin{subfigure}{\textwidth}
  \centering
    \includegraphics[width=.45\linewidth]{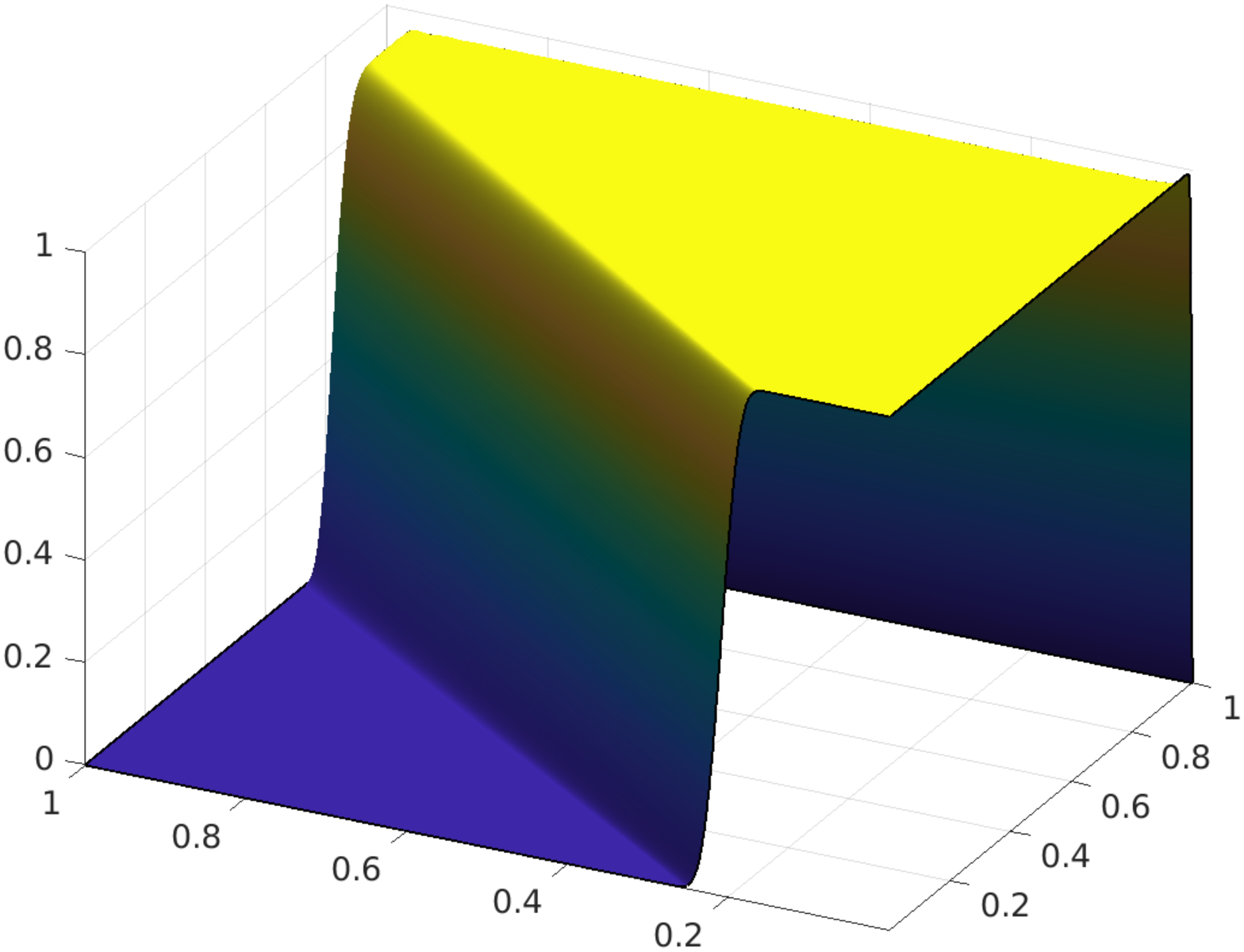}
    \hfill
    \includegraphics[width=.4\linewidth]{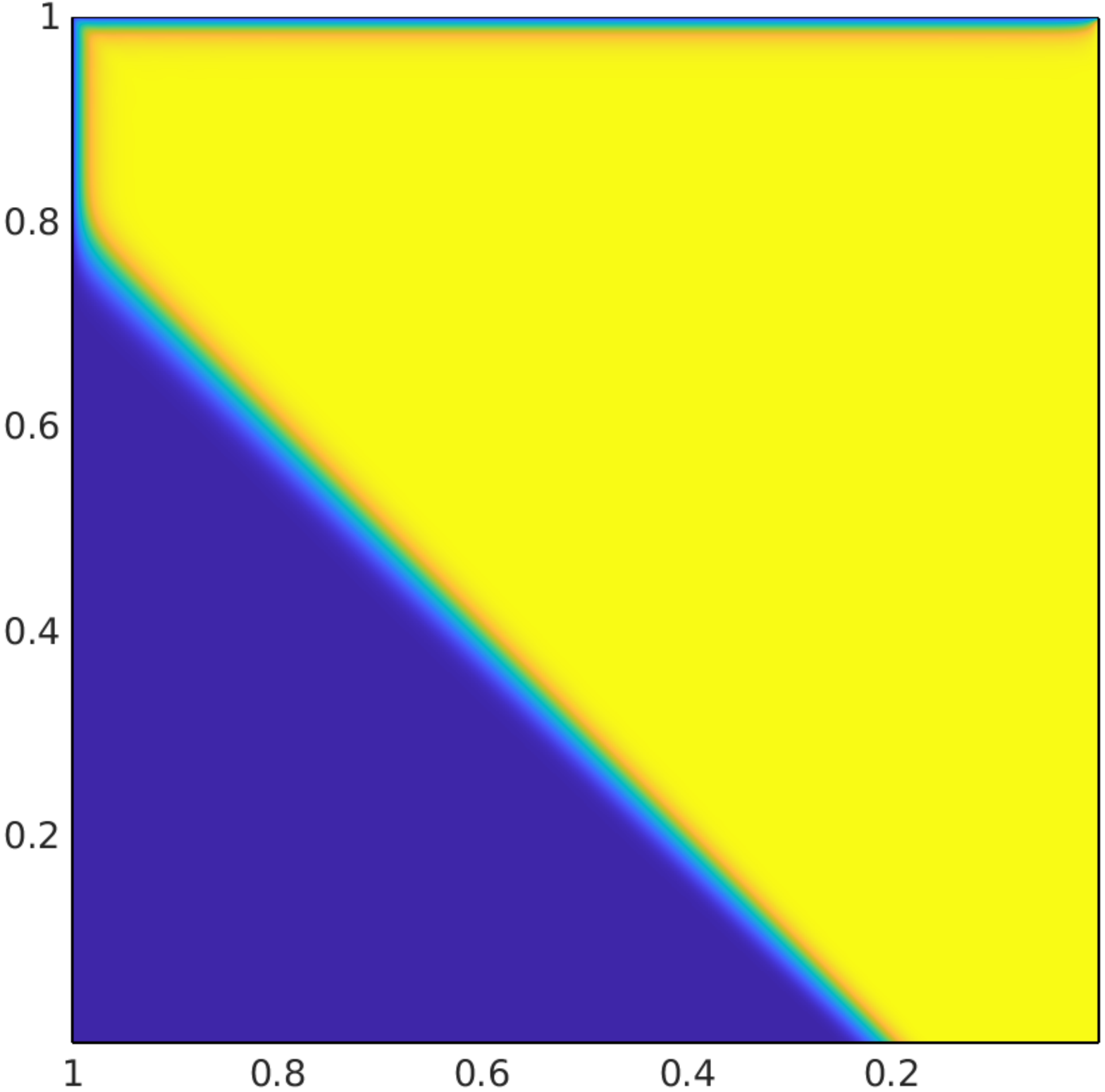}
    \caption{B-spline space: $\spaceT_6 \otimes \spaceT_6$ with SUPG stabilization and $m=50$, dof $= 3136$.} \label{fig_adv_sqr_SUPG}
  \end{subfigure}
  \begin{subfigure}{\textwidth}
  \centering
    \includegraphics[width=.45\linewidth]{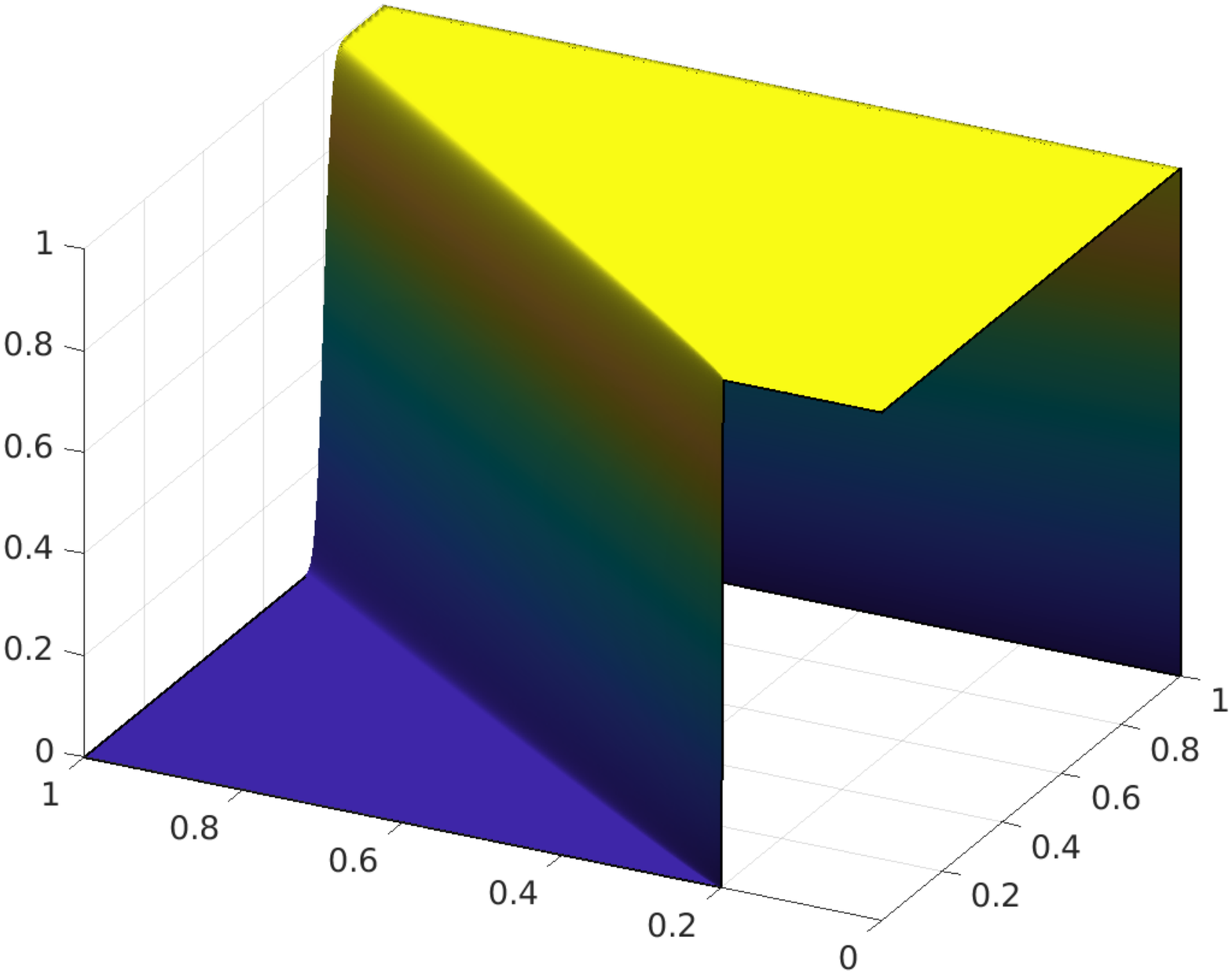}
    \hfill
    \includegraphics[width=.4\linewidth]{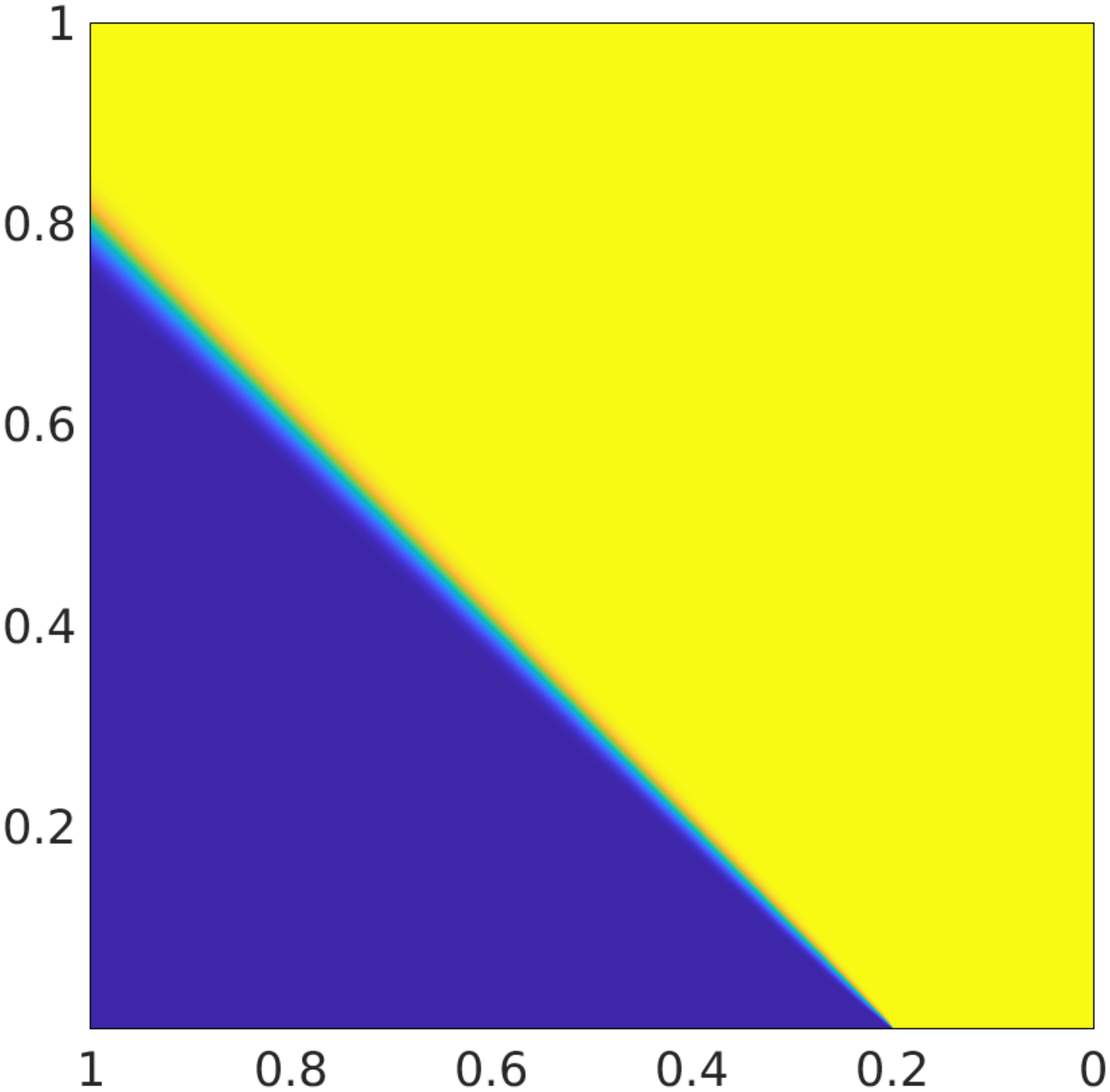}
    \caption{B-spline space: $\spaceT_6 \otimes \spaceT_6$  with $m=1500$, dof $= 2268036$.} \label{fig_adv_sqr_bspline}
  \end{subfigure}
  \begin{subfigure}{\textwidth}
  \centering
    \includegraphics[width=.45\linewidth]{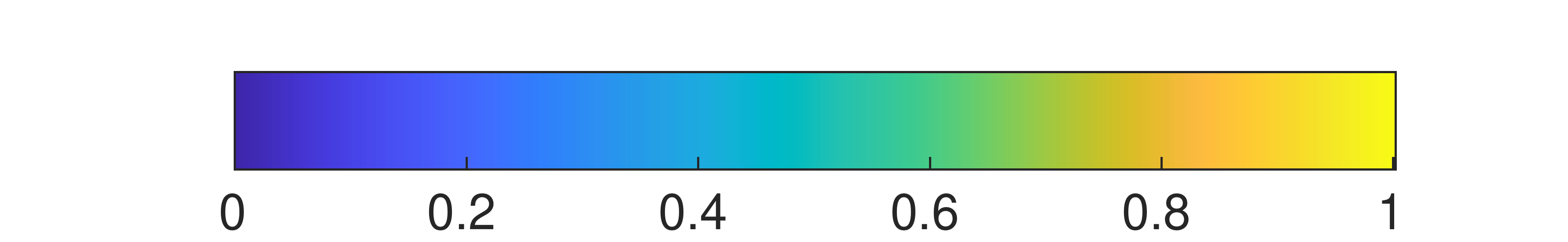}
  \end{subfigure}
\caption{Case study 5. Plots of the approximate solution to the advection-diffusion problem with advection skew at $\theta=45^\circ$ and global \peclet number $\bm{Pe}_g= 10^4 $, obtained by using tensor-product spline spaces with $p_1 = p_2 = 6$; the number of elements $m$ in each direction and the dof are mentioned in the subcaptions.} \label{fig_ex_square}
\end{figure}

In \Cref{fig_adv_sqr_exp} we illustrate the approximate solution obtained by using tensor-product TB-splines,  identified by \eqref{eq_space_square}, of degrees $p_1=p_2=6$ on a mesh consisting of $50 \times 50$ elements. If we use classical tensor-product B-splines of the same degrees on a mesh of only $50 \times 50$ elements, then the result will exhibit huge spurious oscillations near the layers. Hence, it makes sense to compare the stabilized solution of the B-splines against the TB-splines (without stabilization).  
\Cref{fig_adv_sqr_SUPG} shows the result obtained by using tensor-product B-splines with SUPG stabilization, where the stabilization constant is given by 
\begin{equation*}
\tau = \frac{h_{\mathbf{a}}}{2 \|\mathbf{a}\|},\quad \text{with }\quad h_{\mathbf{a}} = \frac{h}{\max\bigl\{\cos(\theta), \sin(\theta)\bigr\}}.
\end{equation*}
Note that $h_{\mathbf{a}}$ is the element length in the direction of the advection flow; see \cite{hughes2005} and \cite[Section~9.3.2]{IGAbook}. As a reference solution we are considering the one obtained by classical tensor-product B-splines of the same degrees on a very fine mesh of $1500 \times 1500$ elements, as shown in \Cref{fig_adv_sqr_bspline}.

\begin{table}[t!] 
\begin{center}
\renewcommand{\arraystretch}{1.25}
\begin{tabular} { >{\centering}p{3cm}  >{\centering}p{3.5cm}  >{\centering}p{1cm}  >{\centering}p{1.75cm} >{\centering}p{1.5cm} >{\centering\arraybackslash}p{2.5cm}}
\hline
TB-spline space   &SUPG stabilization    &$m$    &dof    &max  &min\\ 
\hline\hline
$\spaceT_6^{(\mathfrak{a} \cos(\theta))} \otimes \spaceT_6^{(\mathfrak{a} \sin(\theta) )}$   &No   &50   &3136    &1.0016   &$-1.6385\times10^{-3}$\\
$\spaceT_6 \otimes \spaceT_6$   &Yes   &50   &3136   &1.0009  &$-9.4671\times10^{-4}$\\
$\spaceT_6 \otimes \spaceT_6$   &No   &1500   &2268036  &1.0022   &$-1.8761\times10^{-6}$\\
\hline
\end{tabular}
\renewcommand{\arraystretch}{1}
\caption{Case study 5. Comparison of maximum and minimum of the TB-spline solution with exponential functions in the space against the polynomial B-spline solution with SUPG stabilization on the same mesh and the polynomial B-spline solution on a really high number of intervals with $p_1=p_2=6$ for $\mathfrak{a}=10^4$.}
\label{tab_adv_sqr}
\end{center}
\end{table} 

\Cref{tab_adv_sqr} collects the values of maximal over- and undershoot in the neighbourhood of the layers evaluated on a uniform grid of $3001$ points along each direction for the three considered spline spaces, together with the used number of elements and the total number of dof.  

When comparing the TB-spline solution with the stabilized B-spline solution, we observe that both approaches undoubtedly eliminate the spurious oscillations. However, stabilization produces a too smooth solution and the localization of the internal and boundary layers is less accurate than in the case of the TB-spline solution. On the other hand, the unstabilized B-spline solution on a very fine mesh is able to localize accurately the layers but at a very high cost (i.e., a large number of dof).
In conclusion, TB-splines on coarse meshes are able to predict the layers precisely, resulting in very few required dof while still giving stable results and eliminating the need of stabilization.

\section{Conclusion}
\label{sec-conclusion}
Tchebycheffian splines are smooth piecewise functions whose pieces are drawn from ECT-spaces, a natural generalization of algebraic polynomial spaces. Under suitable assumptions, they admit a basis of functions equipped with all the nice properties of classical polynomial B-splines, the so-called TB-splines. Theoretically, TB-splines are completely plug-to-plug compatible with polynomial B-splines.
Unfortunately, TB-splines in their most generality still lack practical algorithms for their evaluation and manipulation. For this reason, despite their great potential, they did not gain much attention in applications so far.

On the contrary, for the restricted class of TB-splines with pieces belonging to ECT-spaces that are null-spaces of constant-coefficient linear differential operators, efficient manipulation routines have been developed recently and made publicly available in a Matlab toolbox \cite{speleers2022algorithm}, so that they can be easily incorporated in any software library supporting polynomial B-splines to enrich its capability.
This subclass of TB-splines already provides a large variety of combinations of polynomial, exponential, and trigonometric functions equipped with a wide spectrum of shape parameters. They allow for an exact representation of profiles of interest in applications, they behave nicely with respect to differentiation and integration and, by construction, they include fundamental solutions of interesting differential operators. Therefore, they offer a valid alternative to classical polynomial B-splines and NURBS in Galerkin isogeometric methods. 

It turns out that TB-splines can outperform polynomial B-splines whenever
appropriate problem-driven selection strategies for the underlying ECT-spaces are applied. Typically, there are two aspects that can be considered to properly select and profit from a non-polynomial structure in the reference ECT-spaces:
the geometric aspect of the problem based on an exact geometry mapping by using suitable ECT-spaces, and the analytical behavior of the solution that can be better captured, avoiding spurious oscillations, by including in the discretization spaces functions related to fundamental solutions of the given differential operator.

The proposed setting is promising but also presents some theoretical, computational, and numerical issues that require more attention and deserve further investigation in future research.

\begin{itemize}
\item The existence of a TB-spline basis is ensured when the different pieces are taken from a $(p+1)$-dimensional  ECT-space on the entire interval $[a, b]$ which contains constants and its derivative space is a $p$-dimensional ECT-space on the same interval. Since our interest is confined to null-spaces of constant-coefficient linear differential operators, we can weaken this restriction and we only need to check that the critical length of the derivative space we are dealing with exceeds $\xi_{k+p}-\xi_{k+1}$ for $k=1,\ldots,n$. Besides the results concerning the critical length recalled in \Cref{sec-ECT-spaces,sec_null_spc}, we remark that a numerical procedure to estimate the critical length for a given ECT-space is provided in \cite{beccari2020critical}.
In the more general setting where the various pieces of the Tchebycheffian splines are drawn from different ECT-spaces, the existence of a TB-spline basis requires some theoretical assumptions. A numerical characterization for the existence can be found in \cite{beccari2019design}. Alternatively, for practical purposes one can simply use the Matlab toolbox in \cite{speleers2022algorithm} and directly check whether the produced set of functions possess the properties of interest (non-negativity and partition of unity). 

\item The Matlab toolbox in \cite{speleers2022algorithm} provides efficient evaluation routines for a wide range of TB-splines, with pieces belonging to ECT-spaces that are null-spaces of constant-coefficient linear differential operators, of interest in practical applications. 
However, as discussed and illustrated in \cite[Section~6.2]{speleers2022algorithm}, 
possible numerical instabilities can arise when extreme values for shape parameters (i.e., of the roots of the characteristic polynomial \eqref{eq_char_pol}), high degrees, high smoothness, and highly non-uniform
partitions are considered. A ``fidelity'' check for the obtained basis is recommended when dealing with extreme configurations. Such choices are anyway less relevant in most practical applications. Nevertheless, providing
more robust implementations is an interesting but challenging line of further research.

\item As for any Galerkin method, accurate quadrature rules are imperative for isogeometric Galerkin methods based on TB-splines. When using classical elementwise Gaussian quadrature, the same order of quadrature used for polynomial B-splines produces satisfactory results only for mild values of the shape parameters. For larger values of the shape parameters, however, the amount of quadrature points needs to be significantly increased (up to five times). The construction of tailored quadrature rules for TB-splines, in the spirit of \cite{Calabro2019} and references therein, is definitely worth to be investigated.
\end{itemize}

Finally, we observe that, due to the complete structural similarity with polynomial B-splines, the most popular local tensor-product structures supporting local refinement (such as T-splines \cite{SederbergZBN2003}, (truncated) hierarchical B-splines \cite{GiannelliJS2012}, and LR-splines \cite{DokkenLP2013}) can be easily extended to the Tchebycheffian setting; see \cite{BraccoLMRS2016-gb} for an example.

Although the large variety of ECT-spaces that are null-spaces of constant-coefficient linear differential operators already offers an extreme flexible environment for applications, an additional extension is provided by the so-called multi-degree Tchebycheffian B-splines (MDTB-splines); see \cite{hiemstra2020tchebycheffian,NurnbergerSSS1984,speleers2022algorithm}. MDTB-splines allow for ECT-spaces of different dimensions on different intervals of the given partitions; they are also supported by the Matlab toolbox in \cite{speleers2022algorithm}. However, finding (automatic) strategies to select both the dimension and the structure for different ECT-spaces seems to be a very challenging task.


\section*{Acknowledgements}
This work was supported 
by the MSCA-ITN-2019 program through the project GRAPES (contract n.~860843 -- CUP E54I19002020006) and 
by the MIUR Excellence Department Project awarded to the Department of Mathematics, University of Rome Tor Vergata (CUP E83C18000100006).
The last two authors are members of Gruppo Nazionale per il Calcolo Scientifico, Istituto Nazionale di Alta Matematica.

\bibliographystyle{amsplain}
\bibliography{sections/references}

\end{document}